\documentclass[preprint]{elsarticle}       
\usepackage{lipsum}
\makeatletter
\def\ps@pprintTitle{%
 \let\@oddhead\@empty
 \let\@evenhead\@empty
 \def\@oddfoot{}%
 \let\@evenfoot\@oddfoot}
\makeatother

\usepackage{lscape}
\usepackage{lineno,hyperref}
\usepackage{rotating}
\usepackage{makeidx}
\usepackage{float}
\usepackage{epsfig}
\usepackage{amsmath,amssymb}
\usepackage[latin1]{inputenc}
\usepackage{epic,eepic}
\usepackage{overpic}
\usepackage{color}
\usepackage{amsfonts}
\usepackage{graphicx}
\usepackage{mwe} 
\usepackage{subcaption}

\usepackage{algorithm}
\usepackage{algorithmicx}
\usepackage{algpseudocode}
\usepackage{enumitem}

\newenvironment{demo}{\smallskip\noindent{{\it Proof.}}\hskip \labelsep}%
            {\hfill\penalty10000\raisebox{-.09em}{\large\bf\rm $\blacksquare$}\par\medskip}

\newtheorem{theorem}{Theorem}[section]

\newtheorem{corollary}[theorem]{Corollary}

\newtheorem{remark}{Remark}

\DeclareMathOperator*{\argmin}{arg\,min}

\def\balpha{\boldsymbol{\alpha}}
\def\bbeta{\boldsymbol{\beta}}

\def\xx{\mathbf{x}}

\journal{}

\begin{document}

\small

\begin{frontmatter}

\title{Data dependent Moving Least Squares}\tnotetext[label1]{The  fourth author has been supported through project CIAICO/2021/227 (Proyecto financiado por la Conselleria de Innovaci\'on, Universidades, Ciencia y Sociedad digital de la Generalitat Valenciana), by grant PID2020-117211GB-I00 and by PID2023-146836NB-I00 funded by MCIN/AEI/10.13039/501100011033.}

\author[TAU]{David Levin}
\ead{levindd@gmail.com}
\author[UV]{Jos\'e M. Ram\'on}
\ead{Jose.Manuel.Ramon@uv.es}
\author[UPCT]{Juan Ruiz-\'Alvarez}
\ead{juan.ruiz@upct.es}
\author[UV]{Dionisio F. Y\'a\~nez}
\ead{Dionisio.Yanez@uv.es}

\date{Received: date / Accepted: date}

\address[TAU]{School of Mathematical Sciences. Tel-Aviv University, Tel-Aviv (Israel).}
\address[UV]{Departamento de Matem\'aticas. Universidad de Valencia, Valencia (Spain).}
\address[UPCT]{Departamento de Matem\'atica Aplicada y Estad\'istica. Universidad  Polit\'ecnica de Cartagena, Cartagena (Spain).}


\begin{abstract}

In this paper, we address a data dependent modification of the moving least squares (MLS) problem. We propose a novel approach by replacing the traditional weight functions with new functions that assign smaller weights to nodes that are close to discontinuities, while still assigning smaller weights to nodes that are far from the point of approximation. Through this adjustment, we are able to mitigate the undesirable Gibbs phenomenon that appears close to the discontinuities in the classical MLS approach, and reduce the smearing of discontinuities in the final approximation of the original data. The core of our method involves accurately identifying those nodes affected by the presence of discontinuities using smoothness indicators, a concept derived from the data-dependent WENO method. Our formulation results in a data-dependent weighted least squares problem where the weights depend on two factors: the distances between nodes and the point of approximation, and the smoothness of the data in a region of predetermined radius around the nodes. We explore the design of the new data-dependent approximant, analyze its properties including polynomial reproduction, accuracy, and smoothness, and study its impact on diffusion and the Gibbs phenomenon. Numerical experiments are conducted to validate the theoretical findings, and we conclude with some insights and potential directions for future research.

\end{abstract}
\begin{keyword}
WENO \sep high accuracy approximation \sep improved adaption to discontinuities \sep MLS    \sep 41A05 \sep  41A10 \sep 65D05 \sep 65M06 \sep 65N06
\end{keyword}
\end{frontmatter}

\section{Introduction and review}

The Moving Least Squares (MLS) method, that was originally proposed by Shepard in \cite{shepard1968} and further developed by Lancaster and Salkauskas in \cite{lancaster1981surfaces}, is a powerful mathematical tool for generating smooth surfaces from scattered data points and for mess free aproximation of data. The MLS method has been widely applied in various fields such as data approximation \cite{davidlevin}, image processing \cite{lee1997image}, and geometric modelling \cite{wendland}, among others.

In the classical MLS approach, the goal is to approximate a function $f$ given a set of scattered data points $\{(\xx_i, f_i)\}_{i=1}^N$. The more general form of the MLS approximation $\hat{f}(\xx)$ at a point $\xx$ is obtained by minimizing a weighted least squares error:
\begin{equation}\label{problema}
J(a(\xx)) = \sum_{i=1}^N \omega(\xx - \xx_i) \left( f_i - \sum_{j=0}^m a_j(\xx) \phi_j(\xx - \xx_i) \right)^2,
\end{equation}
where $\omega(\xx - \xx_i)$ is a weight function that decreases with the distance between $\xx$ and $\xx_i$, $\phi_j$ are basis functions, and $a_j(\xx)$ are the coefficients to be determined. The weight function $\omega(\xx - \xx_i)$ is designed so that points closer to $\xx$ have a larger influence on the approximation.

Despite its effectiveness in smooth regions, the classical MLS method tends to produce oscillations near jump discontinuities. This limitation arises because the method assumes a smooth underlying function, which is not valid in the presence of discontinuities. Various strategies have been proposed to address this issue, including modifications to the weight function and the introduction of data-dependent techniques \cite{tey2021moving, esfahani2023moving}.

In this work, we propose a data-dependent modification to the classical MLS method, inspired by the Weighted Essentially Non-Oscillatory (WENO) algorithm \cite{Liu, JiangShu}. Our purpose is to handle discontinuities more effectively including a data-dependent modification in the minimization problem (\ref{problema}). Our approach can be interpreted as an artificial adjustment of the distances of points near discontinuities, thereby reducing oscillations and improving the accuracy of the approximation. In the next subsection, we briefly introduce the WENO method, that will serve as inspiration for the data-dependent modification that we propose afterwards.

\subsection{The Weighted Essentially Non-Oscillatory (WENO) method}

The Weighted Essentially Non-Oscillatory (WENO) method was developed to solve hyperbolic partial differential equations with discontinuous solutions \cite{JiangShu,Liu}. The WENO algorithm builds upon the Essentially Non-Oscillatory (ENO) scheme \cite{harten1987uniformly}, aiming to achieve high-order accuracy in smooth regions while avoiding spurious oscillations near discontinuities in the solution of conservation laws.

The key idea of the WENO method is to construct a weighted combination of several candidate stencils. For a given point, the WENO scheme selects the smoothest stencil by assigning weights that diminish the influence of stencils crossing discontinuities. The smoothness indicators are used to measure the smoothness of the function within each stencil, and the data-dependent weights are computed based on these indicators.

Consider a function $f(x)$ that we want to approximate at the point $x_i$. The WENO reconstruction for the function at  $x_{i}$ is given by:

\begin{equation}
\hat{f}(x_{i})=\hat{f}_{i} = \sum_{k=0}^{r-1} \omega_k p_k(x_{i}),
\end{equation}
where $r$ is the number of stencils, $p(x_{i})$ are the polynomial approximations from each stencil, and $\omega_k$ are the data-dependent weights. The weights $\omega_k$ are computed as:
\begin{equation}
\omega_k = \frac{\alpha_k}{\sum_{l=0}^{r-1} \alpha_l},
\end{equation}
with
\begin{equation}
\alpha_k = \frac{C_k}{(\epsilon + I_k)^p},
\end{equation}
where $C_k$ are the linear weights, $I_k$ are the smoothness indicators, $\epsilon$ is a small positive number to avoid division by zero, and $p$ is a parameter typically set to 2 in order to obtain optimal accuracy at smooth zones. The smoothness indicators $I_k$ were designed in \cite{JiangShu} using a measure of the smoothness of the underlying data inspired in the total variation:
\begin{equation}
I_k = \sum_{l=1}^{r-1} \int_{x_{i-1/2}}^{x_{i+1/2}} \Delta x^{2l-1} \left( \frac{d^l}{dx^l} f_k(x) \right)^2 dx,
\end{equation}
where $\Delta x$ is the grid spacing. Given the nature of the problem we aim to solve, we will define the smoothness indicators differently. This approach is more suitable considering the scattered distribution of data points we are assuming. In the next subsection we explain the particular setting of our problem, including the description of the data and its distribution over the considered domain.

\subsection{Our setting}

We consider $\Omega\subseteq \mathbb{R}^n$ an open set, $\chi_N=\{\xx_i\in \Omega:i=1,\hdots,N\}$, that contains $N$ distinct nodes, and  $\mathcal{F}_N=\{f_i=f(\xx_i):i=1,\hdots,N\}$, that is the corresponding set of function values, where $f:\Omega\to\mathbb{R}$ is unknown. Throughout this paper, we assume that the nodes are quasi-equally spaced (or equally spaced), and we define the fill distance 
(see, e.g. \cite{FASSHAUER,wendland,wendland2002}) as:
\begin{equation}\label{filldistance}
h =\sup_{x\in\Omega}\min_{\xx_i\in\chi_N}\|\xx-\xx_i\|,
\end{equation}
we also choose a non-negative and compactly supported radial function $\omega:\Omega\to\mathbb{R}$, and define
$\omega_i(\xx)=\omega(\frac{\|\xx_i-\xx\|}{h}),$ where $\|\cdot\|$ is the Euclidean norm in $\mathbb{R}^n$ (but any other norm can be used).

\begin{figure}[H]
\begin{center}
		\begin{tabular}{ccc}
\hspace{-1cm}	Original& \hspace{-1cm}	   MLS & \hspace{-2cm}	MLS\\
\hspace{-1cm}		\includegraphics[width=6.3cm]{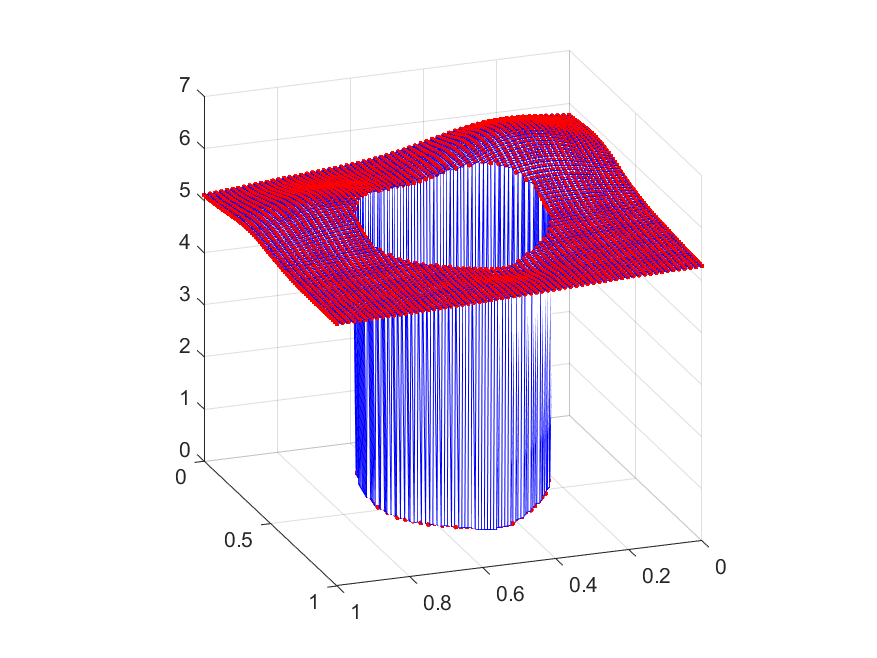} & \hspace{-1cm}	
			\includegraphics[width=6.3cm]{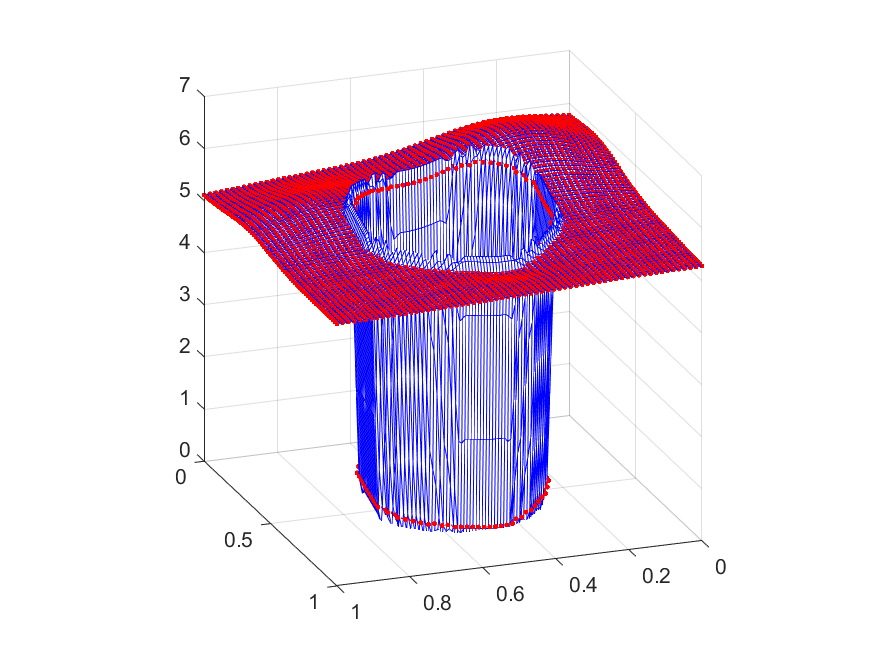} & \hspace{-1cm}
			\includegraphics[width=6.3cm]{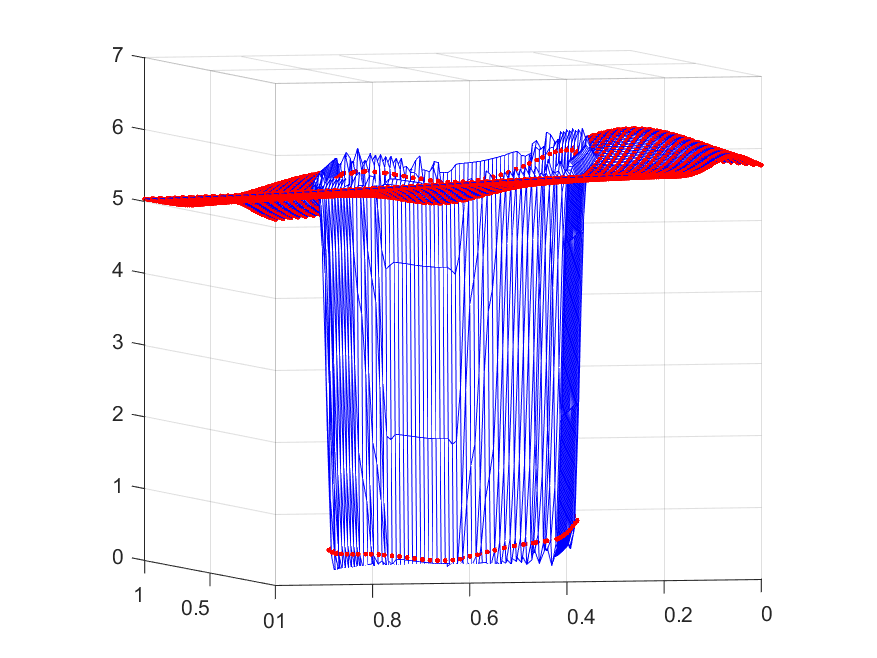}
		\end{tabular}
\end{center}
\caption{Approximation to the function $f$, Eq. \eqref{funcion1}, where $f_1$ is the Franke's function, $f_2$ is the Franke's function plus a constant, and the discontinuity curve is defined by the zero level of the level-set function $\gamma(x,y)=0.25^2-x^2-y^2$. We have used the MLS with $\omega(x)$ the $\mathcal{C}^2$ Wendland function \cite{wendland2002} and the class of polynomials $\Pi_2(\mathbb{R}^2)$. Red points: original function, blue lines: approximation.}
	\label{figure_ejemplo}
	\end{figure}

Let us consider a point $\xx_0\in \Omega$. We particularize the moving least squares problem described in equation ((\ref{problema})), to the polynomial approximation of the values $\{f(\xx_i)\}_{i=1}^N$. This involves calculating a polynomial of degree less than or equal to $d$ that closely approximates the given values $f(\xx_i)$ at the points $\xx_i$, and that satisfies:

\begin{equation}\label{MLSproblem}
\hat{p}_{\xx_0}=\argmin_{p \in \Pi_{d}(\mathbb{R}^n)} \sum_{i=1}^N (f(\xx_{i})-p(\xx_i))^2\omega_i(\xx_0).
\end{equation}
Then, it is evaluated at $\xx_0$, obtaining the approximation
$$f(\xx_0)\approx \hat{p}_{\xx_0}(\xx_0).$$
We can replace the condition that the function $\omega(r)$ is compactly supported with the requirement that it decreases rapidly as $r\to\infty$. In this case, interpolation is achieved if $\omega(0)=\infty$, \cite{davidlevin}.
There are many formulations for this problem in approximation theory (see e.g., \cite{BG1,BG2,BS,FASSHAUER,davidlevin}) and in statistics, where it is known as {\it local polynomial regression} (see e.g., \cite{FASSHAUER, LOADER}). The core idea of the method is to give more importance to the nodes near the point where we want to approximate. This method is effective because it can reproduce polynomials of degree $d$, which implies an accuracy order of $d+1$, and its smoothness depends on the chosen function $\omega$. However, if the data points have a strong gradient or come from a discontinuous function, some oscillations may appear near the discontinuity, as shown in Figure \ref{figure_ejemplo}.

In particular, this paper considers a curve \(\Gamma\) defined as the zero level set of a continuous level-set function \(\gamma: \Omega \to \mathbb{R}\). This function \(\gamma\) delineates two distinct sets within the domain \(\Omega\)
\begin{equation}
\begin{split}
\Omega^+&=\{\xx\in \Omega: \gamma(\xx)\geq 0\}, \\
\Omega^-&=\Omega\setminus \Omega^+.
\end{split}
\end{equation}
The unknown function $f$ is defined as:
\begin{equation}\label{funcion1}
f(\xx)=\left\{
      \begin{array}{ll}
        f_1(\xx), & \xx \in \Omega^+, \\
        f_2(\xx), & \xx \in \Omega \setminus \Omega^+,
      \end{array}
    \right.
\end{equation}
with $f_1\in \mathcal{C}^{d+1}(\overline{\Omega^+})$ and $f_2\in \mathcal{C}^{d+1}(\overline{\Omega^-})$.
In this paper, we replace the functions $\omega_i$ in (\ref{MLSproblem}), which determine the importance of each node, with a new function $\mathcal{W}_i$ that assigns a greater weight if the nodes are far from the discontinuity. In this way, we avoid the undesired effects produced by these nodes, and the Gibbs phenomenon is mitigated. The key to the method is to correctly detect the {\em infected} nodes using {\it smoothness indicators}, a concept defined in the context of data-dependent methods such as WENO (see e.g., \cite{doi:10.1137/070679065}). Therefore, our problem is a weighted least squares problem where the weights depend on different aspects, in this case, two: the distances between the nodes and the point to approximate, and the distance between the isolated discontinuity and the nodes. We could change these particular requirements to others, such as the density of points or the monotonicity of the data.

We divide the paper into four sections: We start by defining the ingredients to design the new data dependent approximator in Section \ref{nlmls}. The next section is devoted to analyzing the properties of the method, such as the reproduction of polynomials, the accuracy, and the smoothness presented by the new operator. Additionally, a study about the diffusion and the Gibbs phenomenon is shown. Some numerical experiments are performed to check the theoretical results in Section \ref{numericalexperiments}. Finally, some conclusions and future work are indicated in the last section.

\section{The data dependent MLS method}\label{nlmls}
The MLS problem that we consider, Eq. \eqref{MLSproblem}, can be reformulated from an algebraic point of view (see e.g., \cite{LOADER, lopezyanez}). We explain it for $n=2$, but it is similar for any $n$. Thus, we consider $\Omega \subset \mathbb{R}^2$, $\chi_N = \{\xx_i = (x_i, y_i) \in \Omega : i = 1, \hdots, N\}$, the set of polynomials of degree less than or equal to $d$
$$\Pi_d(\mathbb{R}^2)=\{\sum_{0\leq |\balpha|\leq d} a_{(\alpha_0,\alpha_1)}x^{\alpha_0} y^{\alpha_1} :a_{\balpha}\in\mathbb{R}, \balpha=(\alpha_0,\alpha_1) \in \mathbb{N}^2 \},$$
and a basis of $\Pi_d(\mathbb{R}^2)$ defined by:
$$\mathcal{B}^d_{\xx_0}=\{1,(x-x_0),(y-y_0),(x-x_0)^2,(y-y_0)^2,(x-x_0)(y-y_0),\hdots,(x-x_0)^d,\hdots,(y-y_0)^d\}.$$
Let $\xx_0=(x_0,y_0)\in \Omega$, we define  the matrices $\mathbb{X}(\xx_0) \in \mathbb{R}^{N\times \binom{d+2}{2}}$, $\mathbb{W}(\xx_0) \in \mathbb{R}^{N\times N}$ as:
$$
\mathbb{X}(\xx_0)=\begin{bmatrix}
1 & (x_1-x_0) & (y_1-y_0) & \hdots & (y_1-y_0)^d\\
1 & (x_2-x_0) & (y_2-y_0) & \hdots & (y_2-y_0)^d\\
1 & (x_3-x_0) & (y_3-y_0) & \hdots & (y_3-y_0)^d\\
\vdots & \vdots & \vdots & \ddots & \vdots\\
1 & (x_N-x_0) & (y_3-y_0) & \hdots & (y_N-y_0)^d
\end{bmatrix}, \quad
\mathbb{W}(\xx_0)=\begin{bmatrix}
\omega_1(\xx_0) & 0 & 0& \hdots & 0\\
0 & \omega_2(\xx_0) & 0 & \hdots & 0\\
0 & 0 & \omega_3(\xx_0) & \hdots & 0\\
\vdots & \vdots & \vdots & \ddots & \vdots\\
0 & 0 & 0 & \hdots & \omega_N(\xx_0)
\end{bmatrix}.
$$
If we write
$$\hat{p}_{\xx_0}(\xx)=\sum_{i=1}^{\binom{d+2}{2}} c_i(\xx_0) A_i(\xx),$$
with
$$A(\xx)=[1,(x-x_0),(y-y_0),(x-x_0)^2,(y-y_0)^2,(x-x_0)(y-y_0),\hdots,(x-x_0)^d,\hdots,(y-y_0)^d],$$
then the problem defined in Eq. \eqref{MLSproblem} can be expressed as:
\begin{equation}\label{equation111}
\mathbf{c}(\xx_0)=\argmin_{\bbeta \in \mathbb{R}^{\binom{d+2}{2}\times 1}}\|\mathbb{W}^\frac{1}{2}[f(\xx_1),\hdots,f(\xx_N)]^T-\mathbb{W}^\frac{1}{2}\mathbb{X} \bbeta||_2^2,
\end{equation}
whose solution is
$$\mathbf{c}(\xx_0)=(\mathbb{X}(\xx_0)^T\mathbb{W}(\xx_0)\mathbb{X}(\xx_0))^{-1}\mathbb{X}(\xx_0)^T\mathbb{W}(\xx_0)[f(\xx_1),\hdots,f(\xx_N)]^T.$$
So we get that
\begin{equation}\label{solucionsimple}
\mathcal{I}_{\text{MLS}}(\xx_0) \equiv \hat{p}_{\xx_0}(\xx_0)=c_0=[1,0,\hdots,0](\mathbb{X}(\xx_0)^T\mathbb{W}(\xx_0)\mathbb{X}(\xx_0))^{-1}\mathbb{X}(\xx_0)^T\mathbb{W}(\xx_0)[f(\xx_1),\hdots,f(\xx_N)]^T.
\end{equation}

\begin{remark}\label{nota1}
Note that if $\omega$ is compactly supported then some points such that $\omega_i(\xx_0)=0$ could exist. In these cases, we denote as:
$$\tilde{\chi}_{N_0}=\{\xx_{i}\in \chi_N: \omega_{i}(\xx_0)>0\}=\{\xx_{0_i}:i=1,\hdots,N_{0}\},$$
and construct the same problem replacing $\chi_N$ by $\chi_{N_0}$. In the rest of the paper, we consider that $\omega_i(\xx_0)>0$, for all $i=1,\hdots,N$
\end{remark}

In the solution of the MLS problem exposed before, we can see that the relevance of the matrix $\mathbb{W}$ is checked, such that weights are assigned depending on the distance of the nodes with respect to the point where we want to obtain the approximation. If the data points present an isolated discontinuity described by an unknown curve $\gamma$, then not only the distance but also the position of the nodes with respect to this curve are relevant. Therefore, in order to take into account these two variables, we construct a different problem by replacing the weight function with a data-dependent one. Thus, we define
\begin{equation}\label{si}
\tilde{\omega}_i(\xx)=\frac{\omega_i(\xx)}{(\epsilon+I_i)^t}=\frac{\omega(\frac{\|\xx-\xx_i\|}{h})}{(\epsilon+I_i)^t},
\end{equation}
where $\epsilon$ and $t$ are two parameters. We use $\epsilon$  to avoid zero values in the denominator. The purpose of $t$ is to reach the maximum order of accuracy of the approximation at smooth zones (see \cite{doi:10.1137/070679065}). Typically, we use $t=4$. Finally, the values $I_i$ with $i=1,\hdots,N$ are the smoothness indicators, i.e., the values which mark if one node is close to the discontinuity or not. In our case, we determine a ball centered in the node $\xx_i$ with a fixed radius $\delta$ denoted by
$$\mathcal{S}(\xx_i)=B(\xx_i,\delta)\cap \chi_N=\{\xx_{i_j}\in\chi_N:j=1,\hdots,N_i\},$$
and impose the following conditions (see e.g. \cite{ABM}):
\begin{enumerate}[label={\bfseries P\arabic*}]
\item\label{P1sm1d} The order of a smoothness indicator that is free of discontinuities is $h^2$, i.e.
$I_{i}=O(h^2) \,\, \text{if}\,\, f \,\, \text{is smooth in } \,\, \mathcal{S}_{i}.$

\item\label{P2sm1d} When the stencil $\mathcal{S}_{i}$ crosses a discontinuity, then
$I_{i} \nrightarrow 0 \,\, \text{as}\,\, h\to 0.$
\end{enumerate}


In this paper, we define the smoothness indicators in the following way: First we solve the linear least squares problem:
\begin{equation}\label{problemaparadetectar}
p_i=\argmin_{p\in\Pi_1(\mathbb{R}^2)} \sum_{i_j=1}^{N_i} (f(\xx_{i_j})-p(\xx_{i_j}))^2,
\end{equation}
and compute the smoothness indicators as:
\begin{equation}\label{indicador1}
I_i=\frac{1}{N_i}\sum_{j_i=1}^{N_i}|f(\xx_{j_i})-p_i(\xx_{j_i})|.
\end{equation}
With this definition, $I_i$ satisfies \ref{P1sm1d} and \ref{P2sm1d}. Now, with the new functions $\tilde{\omega}_i$, defined in (\ref{si}), we can pose the weighted least squares problem and find its solution. To do so, we define
$$
\tilde{\mathbb{W}}(\xx_0)=\begin{bmatrix}
\tilde{\omega}_1(\xx_0) & 0 & 0& \hdots & 0\\
0 & \tilde{\omega}_2(\xx_0) & 0 & \hdots & 0\\
0 & 0 & \tilde{\omega}_3(\xx_0) & \hdots & 0\\
\vdots & \vdots & \vdots & \ddots & \vdots\\
0 & 0 & 0 & \hdots & \tilde{\omega}_N(\xx_0)
\end{bmatrix},
$$
and get the new data dependent approximation DD-MLS as in Eq. \eqref{solucionsimple}:
\begin{equation}\label{solunueva}
\mathcal{I}_{\text{DD-MLS}}(\xx_0) \equiv
\tilde{\hat{p}}_{\xx_0}(\xx_0)=\tilde{c}_0=[1,0,\hdots,0](\mathbb{X}(\xx_0)^T\tilde{\mathbb{W}}(\xx_0)\mathbb{X}(\xx_0))^{-1}\mathbb{X}(\xx_0)^T\tilde{\mathbb{W}}(\xx_0)[f(\xx_1),\hdots,f(\xx_N)]^T.
\end{equation}
If we write $\lambda_i=(\epsilon+I_i)^t$, then we have that $\tilde{\omega}_i$ is a radial function since
 $$\tilde{\omega}_i(\xx)=\lambda_i\omega_i(\xx)=\lambda_i\omega\left(\frac{\|\xx-\xx_i\|}{h}\right),$$
 with compact support, which assigns small weights (or zero) when the point $\xx$ is far from $\xx_i$, but also, from the fact that $\lambda_i=O(1)$ or $\lambda_i=O(h^{-2t})$ depending on whether $\xx_i$ is close to or far from a discontinuity, the function $\tilde{\omega}_i$ assigns larger weights to {\it non-infected} points and smaller (close to zero) weights to {\it infected} ones {(meaning by {\it infected}, nodes that are at a distance smaller than $\delta$ from the discontinuity)}. In other words, we can interpret this data-dependent modification as a a change in the distance to the nodes that are close to discontinuities. All the nodes close to the discontinuity are considered to be far from any point and their importance is neglected in the final approximation. With these ingredients,  in the next section we prove some properties of the new approximation technique.
\section{Properties of the new method}\label{properties}
In this section, we prove some properties of the approximation technique described in the previous section. In particular, we focus on the smoothness, the reproduction of polynomials, the order of accuracy and the elimination of Gibbs phenomenon:
\begin{itemize}
\item Let us start with the smoothness: it is clear by Eq. \eqref{solunueva} that if $\mathbb{X}(\xx_0)$ has maximum rank, i.e. $\binom{d+2}{2}$, then the matrix
$\mathbb{X}(\xx_0)^T\tilde{\mathbb{W}}(\xx_0)\mathbb{X}(\xx_0)$ is non-singular, since all the entries of the matrix $\tilde{\mathbb{W}}(\xx_0)$
are positive by Remark \ref{nota1}. Therefore, the smoothness of the new approximator depends only on the function $\omega$ since $\tilde{\omega}_i=\lambda_i\omega$. We summarize this property in the next Theorem (Th. \ref{teo1})

\begin{theorem}\label{teo1}
Let $\nu\in\mathbb{N}$, $\Omega\subset \mathbb{R}^n$, $\omega$ be a function with $\omega\in\mathcal{C}^\nu(\Omega)$, then
the new approximation defined
in Eq. \eqref{solunueva} is $\mathcal{C}^\nu(\Omega)$.
\end{theorem}

In the literature there are many functions $\omega$ which are used in radial basis approximation. We summarize some of them in Table \ref{tabla1nucleos}  (see e.g. \cite{FASSHAUER,wendland2002}).


\begin{table}[h!]
\centering
\begin{tabular}{lll}
\hline
$\omega(  \,r)$ & RBF &                                                          \\ \hline
$e^{-  r^2}$ & Gaussian $C^\infty$ & G                                          \\
$\left(1 + r^2 \right)^{-1/2}$ & Inverse MultiQuadratic $C^\infty$ & IMQ  \\
$e^{-  r}$ & Mat\'ern $C^0$ & M0 \\
$e^{-  r} \left( 1 +   r \right)$ & Mat\'ern $C^2$ & M2 \\
$e^{-  r} \left( 3 + 3  r +   r^2 \right)$ & Mat\'ern $C^4$ & M4 \\
$(1 -   r)^2_+$ & Wendland $C^0$ & W0 \\
$(1 -   r)^4_+ \left( 4  r + 1 \right)$ & Wendland $C^2$ & W2 \\
$(1 -   r)^6_+ \left( 35  r^2 + 18  r + 3 \right)$ & Wendland $C^4$ & W4 \\ \hline
\end{tabular}
\caption{Examples of RBFs.}\label{tabla1nucleos}
\end{table}

\item The second property is the reproduction of polynomials. As $\mathbb{X}(\xx_0)^T\tilde{\mathbb{W}}(\xx_0)\mathbb{X}(\xx_0)$ is non-singular, the system has a unique solution. If the data are the discretization of a polynomial of degree less than or equal to $d$, then the solution of the problem Eq. \eqref{MLSproblem} with $\tilde{\omega}_i$ instead of ${\omega}_i$ is the same. Therefore, the interpolator reproduces $\Pi_d(\mathbb{R}^n)$. The next result collects this property.

\begin{theorem}\label{repropolis}
Let $\Omega\subset \mathbb{R}^n$ be an open set, $\xx_0\in \Omega$, $\chi_N=\{\xx_i\in \Omega:i=1,\hdots,N\}$ a set of $N$ distinct nodes and  $\mathcal{F}_N=\{f_i=f(\xx_i):i=1,\hdots,N\}$ a corresponding set of function values with $f\in\Pi_d(\mathbb{R}^n)$. Then the data dependent MLS approximation defined in Eq. \eqref{solunueva} satisfies
$$\tilde{\hat{p}}_{\xx_0}(\xx_0)=f(\xx_0).$$
\end{theorem}
A direct consequence of the last Theorem, Th. \ref{repropolis}, is the order of accuracy. If we consider the fill distance $h$, Eq. \eqref{filldistance}, then we can assure that if the function is $d+1$ continuous, the order of accuracy is $O(h^{d+1})$. Thus, we can enunciate the following corollary.
\begin{corollary}\label{cororder}
Let $\Omega\subset \mathbb{R}^n$. If $f\in\mathcal{C}^{d+1}(\overline{\Omega})$, $\chi_N=\{\xx_i\in \Omega:i=1,\hdots,N\}$ are quasi-uniformly distributed with fill distance $h$, the weight function $\omega$ is compactly supported with support size $c$, then the new approximation defined in Eq. \eqref{solunueva} fulfills
$$|f(\xx_0)-\tilde{\hat{p}}_{\xx_0}(\xx_0)|\leq Ch^{d+1}\max_{\boldsymbol{\xi}\in \overline{\Omega}}|D^{\boldsymbol{\alpha}} f(\boldsymbol{\xi})|, \quad|\boldsymbol{\alpha}|=d+1,\,\, \forall \, \xx_0 \in \Omega,$$
where $C$ is a constant independent of $h$.
\end{corollary}

\item Finally, we analyze the approximation when some points far enough from the discontinuities are used mixed with {\it infected} points. In this way, we will suppose that we have at least $\binom{d+2}{2}$ data points not marked as {\it infected}.
\end{itemize}

\begin{theorem}\label{repropolis2}
Let $\Omega\subset \mathbb{R}^2$ be an open set, $\chi_N=\{\xx_i\in \Omega:i=1,\hdots,N\}$ a set of $N$ distinct nodes with fill distance $h$, and  $\mathcal{F}_N=\{f_i=f(\xx_i):i=1,\hdots,N\}$ the corresponding set of function values with $f$ defined in Eq. \eqref{funcion1}.
Let $\xx_0\in \Omega^+$ and let
$$\tilde{\chi}_{N_0}=\{\xx_{i}\in \chi_N: \omega_{i}(\xx_0)>0\}=\{\xx_{0_i}\in \chi_N:i=1,\hdots,N_{0}\},$$
as the set of points used to calculate the approximation at $\xx_0$, with $N_0\ge \binom{d+2}{2}$. We assume there exist at least $\binom{d+2}{2}$ points belonging to $\tilde{\chi}_{N_0}$, denoted by $P_d$,  such that their smoothness indicators are of order $h^2$, all located within $\Omega^+$. We further assume that the interpolation operator in $\Pi_d(\mathbb{R}^2)$  defined using these points is well-defined, and that a constant independent of $h$ bounds its norm. Then the data dependent MLS approximation satisfies:
$$|f(\xx_0)-\mathcal{I}_{\text{DD-MLS}}(\xx_0|=O(h^{\min\{d+1,\kappa\}}),$$
with $\kappa\geq t$.
\end{theorem}
\begin{demo}
Let $\xx_0$ be a point in $\Omega$. We divide the set
$$\tilde{\chi}_{N_0}=\tilde{\chi}^f_{N_0}\cup \tilde{\chi}^d_{N_0},$$
with $\tilde{\chi}^f_{N_0}=\{\xx_{i}\in \chi_N: \omega_{i}(\xx_0)>0\,\, \wedge \,\, I_i=O(h^2)\}$ and $\tilde{\chi}^d_{N_0}=\tilde{\chi}_{N_0}\setminus\tilde{\chi}^f_{N_0}$. Let us define $\pi_1\in \Pi_d(\mathbb{R}^2)$ to be polynomial interpolating $f$ at the points $P_d$ and, as all the non-infected data are at one side of the discontinuity, then
$$f(\xx_i)=f_1(\xx_i)=\pi_1(\xx_i)+O(h^{d+1}), \quad \forall \xx_i \in \tilde{\chi}^f_{N_0}.$$
The MLS solution is defined by
\begin{equation*}
\tilde{\hat{p}}_{\xx_0}=\argmin_{p \in \Pi_{d}(\mathbb{R}^2)} \sum_{i=1}^{N_0} (f(\xx_{0_i})-p(\xx_{0_i}))^2\tilde{\omega}_{0_i}(\xx_0)=
\argmin_{p \in \Pi_{d}(\mathbb{R}^2)} \left(\sum_{\xx_{0_i}\in \tilde{\chi}^d_{N_0}} (f(\xx_{0_i})-p(\xx_{0_i}))^2\tilde{\omega}_{0_i}(\xx_0)+\sum_{\xx_{0_i}\in\tilde{\chi}^f_{N_0}} (f(\xx_{0_i})-p(\xx_{0_i}))^2\tilde{\omega}_{0_i}(\xx_0)\right).
\end{equation*}
Since $\pi_1\in \Pi_d(\mathbb{R}^2)$ and
\begin{equation}\label{omega}
\tilde{\omega}_{0_i}(\xx_0)=\begin{cases}
O(h^{-2t}),& \xx_{0_i}\in\tilde{\chi}^f_{N_0},\\
O(1),& \xx_{0_i}\in\tilde{\chi}^d_{N_0},\\
\end{cases}
\end{equation}
it follows that
\begin{equation*}
\begin{split}
\sum_{\xx_{0_i}\in \tilde{\chi}^d_{N_0}} (f(\xx_{0_i})-\tilde{\hat{p}}_{\xx_0}(\xx_{0_i}))^2\tilde{\omega}_{0_i}(\xx_0)&+\sum_{\xx_{0_i}\in\tilde{\chi}^f_{N_0}} (f_1(\xx_{0_i})-\tilde{\hat{p}}_{\xx_0}(\xx_{0_i}))^2\tilde{\omega}_{0_i}(\xx_0)
\leq\\
&\leq\sum_{\xx_{0_i}\in\tilde{\chi}^d_{N_0}} (f(\xx_{0_i})-\pi_1(\xx_{0_i}))^2\tilde{\omega}_{0_i}(\xx_0) +O(h^{2(d-t+1)}) \Rightarrow
\end{split}
\end{equation*}
\begin{equation*}
\begin{split}
\sum_{\xx_{0_i}\in \tilde{\chi}^f_{N_0}} (f_1(\xx_{0_i})-\tilde{\hat{p}}_{\xx_0}(\xx_{0_i}))^2\tilde{\omega}_{0_i}(\xx_0)
\leq &\sum_{\xx_{0_i}\in\tilde{\chi}^d_{N_0}} (f(\xx_{0_i})-\pi_1(\xx_{0_i}))^2\tilde{\omega}_{0_i}(\xx_0)\\&-\sum_{\xx_{0_i}\in\tilde{\chi}^d_{N_0}} (f(\xx_{0_i})-\tilde{\hat{p}}_{\xx_0}(\xx_{0_i}))^2\tilde{\omega}_{0_i}(\xx_0)+O(h^{2(d-t+1)})\\
=&\sum_{\xx_{0_i}\in\tilde{\chi}^d_{N_0}} (\pi_1(\xx_{0_i})-\tilde{\hat{p}}_{\xx_0}(\xx_{0_i}))(\pi_1(\xx_{0_i})-f(\xx_{0_i})+\tilde{\hat{p}}_{\xx_0}(\xx_{0_i})-f(\xx_{0_i}))\tilde{\omega}_{0_i}(\xx_0)\\
&+O(h^{2(d-t+1)}).
\end{split}
\end{equation*}
From \eqref{omega},
we know that the right-hand side of the inequality is bounded, then there exists a constant $C$ such that
\begin{equation}
\sum_{\xx_{0_i}\in\tilde{\chi}^d_{N_0}} (\pi_1(\xx_{0_i})-\tilde{\hat{p}}_{\xx_0}(\xx_{0_i}))(\pi_1(\xx_{0_i})-f(\xx_{0_i})+\tilde{\hat{p}}_{\xx_0}(\xx_{0_i})-f(\xx_{0_i}))\tilde{\omega}_{0_i}(\xx_0)\leq C,
\end{equation}
independent of $h$. By \eqref{omega}, $\tilde{\omega}_{0_i}(\xx_0)=O(h^{-2t})$ if $\xx_{0_i}\in\tilde{\chi}^f_{N_0}$, and this implies that
\begin{equation*}
\sum_{\xx_{0_i}\in \tilde{\chi}^f_{N_0}} (f_1(\xx_{0_i})-\tilde{\hat{p}}_{\xx_0}(\xx_{0_i}))^2\tilde{\omega}_{0_i}(\xx_0)=O(h^{-2t})\sum_{\xx_{0_i}\in \tilde{\chi}^f_{N_0}} (f_1(\xx_{0_i})-\tilde{\hat{p}}_{\xx_0}(\xx_{0_i}))^2
\leq C+O(h^{2(d-t+1)}),
\end{equation*}
then
\begin{equation*}
(f_1(\xx_{0_i})-\tilde{\hat{p}}_{\xx_0}(\xx_{0_i}))^2 = O(h^{2\kappa})\Rightarrow  |f_1(\xx_{0_i})-\tilde{\hat{p}}_{\xx_0}(\xx_{0_i})| = O(h^{\kappa}), \quad \forall\,\xx_{0_i}\in \tilde{\chi}^f_{N_0},
\end{equation*}
with $\kappa\geq t$. Therefore,
\begin{equation*}
\begin{split}
|\pi_1(\xx_{0_i})-\tilde{\hat{p}}_{\xx_0}(\xx_{0_i})|&=|\pi_1(\xx_{0_i})-f_1(\xx_{0_i})+f_1(\xx_{0_i})-\tilde{\hat{p}}_{\xx_0}(\xx_{0_i})|\\&\leq |\pi_1(\xx_{0_i})-f_1(\xx_{0_i})|+|f_1(\xx_{0_i})-\tilde{\hat{p}}_{\xx_0}(\xx_{0_i})| = O(h^{d+1})+O(h^{\kappa}), \quad \forall\,\xx_{0_i}\in \tilde{\chi}^f_{N_0}.
\end{split}
\end{equation*}
Since,
the interpolation operator is well-defined with respect to the points $P_d$
and its inverse is bounded by a constant independent of $h$  then
$$||\pi_1-\tilde{\hat{p}}_{\xx_0}||_{\infty,\Omega^+}=O(h^{d+1})+O(h^\kappa).$$
Therefore,
\begin{equation*}
\begin{split}
|f(\xx_{0})-\tilde{\hat{p}}_{\xx_0}(\xx_{0})|&=|f_1(\xx_{0})-\pi_1(\xx_0)+\pi_1(\xx_0)-\tilde{\hat{p}}_{\xx_0}(\xx_{0})|\\
&\leq|f_1(\xx_{0})-\pi_1(\xx_0)|+|\pi_1(\xx_0)-\tilde{\hat{p}}_{\xx_0}(\xx_{0})|=O(h^{d+1})+O(h^{d+1})+O(h^\kappa)\\
&=O(h^{\min\{d+1,\kappa\}}).
\end{split}
\end{equation*}
\end{demo}

In our case, we only have to take $t=d+1$ in order to get the maximum order of accuracy.

\begin{remark}{\bf Diffusion and oscillations in the MLS method.}
When applying the linear MLS method to discontinuous data with a radial weight function $\omega$ of support $c$, the resulting approximation exhibits significant diffusion and oscillations within a region of width
$\le ch$ near the discontinuity curve $\Gamma$.
\end{remark}

The following theorem demonstrates that the diffusion region is significantly reduced when using the data dependent MLS method, and, additionally, oscillations are eliminated.
We use the above definitions of $\Gamma$, $\Omega^+$ and $\Omega^-$ and of $f$ defined in Eq. \eqref{funcion1}.
We consider a weight function $\omega$ of support size $c$, $c\ge 2$, and such that $\omega(c-\epsilon_0)\ge C_0>0$. We further assume that $\Gamma$ is a smooth curve. We adopt below the assumptions in Theorem \ref{repropolis2} and the notation within its proof. Let $D_r$ denote the domain consisting of all points within a distance $\le r$ from $\Gamma$, and let
$\bar D_r$ denote its complement.

\begin{theorem}\label{teo3new}
Assume the data points $\tilde{\chi}^d_{N_0}$ are within $D_{kh}$ and let
$\xx_0\in\Omega^+\cap \bar D_{ch-kh}$. Then the data dependent MLS approximation of degree $d$ with $\omega$ of support size $c\ge d/2+k$ satisfies
$|\tilde{f}(\xx_0)- \mathcal{I}_{\text{DD-MLS}}(\xx_0)|=O(h^{d+1}).$
\end{theorem}

\begin{demo}
Define the support of data points respect to $\xx_0$ as $S(\xx_0)=\{\xx_i\in:\omega_i(\xx_0)>0\}$. For any point $\xx_i\in S(\xx_0)\cap \Omega^-$ the smoothness indicator $I_i$ is $O(1)$. Additionally, within a circular segment of width $dh$ there are points $\xx_i\in S(\xx_0)\cap \Omega^+$ for which the smoothness indicator $I_i$ is $O(h^2)$ as $h\to 0$.

Recalling that $h$ is the fill distance of the points $\chi_N$, it follows that there exist at least $\binom{d+2}{2}$ points points $\xx_i\in S(\xx_0)\cap \Omega^+$ for which the smoothness indicator $I_i$ is $O(h^2)$ as $h\to 0$. For these points, the interpolation operator, defined using them, is well-defined and its norm is bounded by a constant independent of $h$. Then, by Theorem \ref{repropolis2}, using $t=d+1$, the data dependent MLS approximation satisfies:
$$|f(\xx_0)-\mathcal{I}_{\text{DD-MLS}}(\xx_0)|=O(h^{d+1}).$$
\end{demo}

\begin{corollary}{\bf Reduced diffusion region and reduced oscillations}

The linear MLS provides full approximation order for points at a distance $\ge ch$ from $\Gamma$. By applying Theorem \ref{teo3new} we observe that the data dependent MLS extends this full approximation order to a larger region, specifically to points at a distance $\ge ch-kh$ from $\Gamma$.

The oscillations observed near
$\Gamma$ in the linear MLS approximation arise because it incorporates data from both sides of $\Gamma$
and applies both positive and negative weights. This occurs for points at a distance $< ch$ from $\Gamma$. As argued in Theorem \ref{teo3new}, for points within the range $ch-hk<dist<ch$, the data dependent MLS approximation uses data only from one side of $\Gamma$. Since the weights in the approximation sum to 1, and are applied exclusively to data from one side of $\Gamma$, no oscillations will occur within this range.

\end{corollary}

\section{Numerical experiments}\label{numericalexperiments}

In this section, we check numerically the theoretical properties shown in Section \ref{properties}: order of accuracy, reproduction of polynomials, reduction of Gibbs oscillations, and reduction of the smearing around discontinuities.  

The algorithm that we use is similar to the one designed in Chapter 19 of \cite{FASSHAUER}, we only introduce a detection of discontinuity points by calculating the smoothness indicator of each data point, i.e., we solve the problems in Eq. \eqref{problemaparadetectar} and evaluate Eq. \eqref{indicador1}. For this problem, we fix the following parameter as in \cite{cavo}:
$$\delta=\frac{\sqrt{2}}{\lfloor\frac{\sqrt{N}}{2}\rfloor}, \quad i=1,\hdots,N,$$
where $\lfloor x\rfloor$ is the floor function, which returns the greatest integer less than or equal to $x \in \mathbb{R}$.
This part is not very expensive computationally because it only involves solving a simple linear least squares problem. Therefore, the algorithms for the linear and data-dependent cases are very similar in terms of computational cost.

We will use the acronyms $\text{MLS}_{F}^{d}$ or $\text{DD-MLS}_{F}^{d}$ to call the linear and data dependent moving least squares methods, where
$F$ represents the weight function, chosen from Table \ref{tabla1nucleos}, and $d$ is the maximum degree of the polynomial used. We have selected the
Wendland $W2$, $W4$ functions and gaussian function $G$ to perform the experiments. We apply these functions in the following way,
$\omega_i(\xx)=\omega(\varepsilon \|\xx-\xx_i\|)$ with $\varepsilon \propto \frac{1}{h}$. In our experiments, we take:
$$ \varepsilon =\frac{1}{2}\left\lfloor\frac{\sqrt{N}}{2}\right\rfloor.$$
Also, in the Gaussian case, another condition is imposed: we only consider the values $\xx_i\in \chi_N$ such that $\omega_i(\xx)>10^{-10}$.
We divide our experiments in three subsections: we examine the order of accuracy in smooth zones for the values $d=0,1,2$. After that, we analyze the behaviour of the approximations close to the discontinuities, and, finally, we study intensively the smearing of discontinuities when $d=0$ (Shepard's method, \cite{shepard1968}) and $d=1$.

\subsection{Order of accuracy}

We start by analyzing the order of accuracy using the well-known Franke's function, defined as
\begin{equation}\label{frankefunction}
f(x,y)=\frac34 e^{-1/4((9x-2)^2+(9y-2)^2)}+\frac34 e^{-1/49((9x+1)^2-1/10(9y+1)^2)}+\frac12 e^{-1/4((9x-7)^2+(9y-3)^2)}-\frac15 e^{-((9x-4)^2-(9y-7)^2)},
\end{equation}
using as nodes two types of sets: a regular grid defined by \(\chi_{N=2^l+1}=\{(i/2^l,j/2^l):i,j=0,\hdots,2^l\}\), and a collection of \(N=(2^l+1)^2\) Halton's scattered data points, as described in \cite{halton}. We denote the fill distance by \(h_l\), and the associated errors by \(e^l_i=|f(\mathbf{z}_i)-\mathcal{I}^l(\mathbf{z}_i)|\), where \(\{\mathbf{z}_i:1\leq i\leq N_{\text{eval}}\}\) is a regular grid in $[0.025,0.975]^2$, and represents the set of evaluation points at which the function is approximated with $N_{eval}=120^2$ points. Finally, we define the maximum discrete \(\ell^2\) norms and their respective convergence rates as follows:
\begin{equation*}
\begin{split}
\text{MAE}_l=\max_{i=1,\hdots,N_{\text{eval}}}e^l_i, \quad \text{RMSE}_l=\left(\frac{1}{N_{\text{eval}}}\sum_{i=1}^{N_{\text{eval}}}(e_i^l)^2\right)^{\frac12}, \quad r_l^{\infty}=\frac{\log(\text{MAE}_{l-1}/\text{MAE}_{l})}{\log(h_{l-1}/h_{l})}, \quad r_l^2=\frac{\log(\text{RMSE}_{l-1}/\text{RMSE}_{l})}{\log(h_{l-1}/h_{l})}.
\end{split}
\end{equation*}
When $d=2$, we can see in Tables \ref{tabla:d2r} and \ref{tabla:d2h} that the order of accuracy reached is greater than 3. In fact, when the data points are placed on a regular grid, the order of accuracy is close to 4 in both linear and data-dependent cases. These results are shown in Table \ref{tabla:d2r}. We can see that the error is slightly small when the linear version of the MLS is employed. This effect can be reduced if we choose a value of the parameter $t$ in (\ref{si}) smaller, but then some oscillations can appear in the final approximation. When $d=1$, Tables \ref{tabla:d1r} and \ref{tabla:d1h}, the order of accuracy is 2 when using a regular grid, as it can be observed in Table \ref{tabla:d1r}, and also when using pseudo-random points, Table \ref{tabla:d1h}. Again, the error is slightly greater in the data-dependent case. Finally, when $d=0$ the order of accuracy with regular data points is larger than the one expected, as it is shown in Table \ref{tabla:d0r}, but smaller when Halton's data points are used, as shown in Table \ref{tabla:d0h}. It is clear that the results when the data dependent MLS method is used are very close to the linear ones, and that it works correctly for smooth functions and, therefore, in the smooth zones.

\begin{table}[!ht]
\begin{center}
\begin{tabular}{lrrrrrrrrrrr}
& \multicolumn{2}{c}{MLS$^2_{\text{W2}}$} & &   \multicolumn{2}{c}{DD-MLS$^2_{\text{W2}}$}& &\multicolumn{2}{c}{MLS$^2_{\text{W2}}$} & &   \multicolumn{2}{c}{DD-MLS$^2_{\text{W2}}$}\\ \cline{1-3} \cline{5-6} \cline{8-9} \cline{11-12}
$l$ & $\text{MAE}_l$ & $r_l^{\infty}$ & &$\text{MAE}_l$ & $r_l^{\infty}$ & &$\text{RMSE}_l$ & $r_l^2$ & &$\text{RMSE}_l$ & $r_l^2$       \\
\hline
$4$  &             2.9459e-02 &        &&   6.2989e-02 &        &&   4.5011e-03 &        &&   1.2700e-02 &                     \\
$5$  &             3.4607e-03 &  3.1049&&   9.5840e-03 &  2.7299&&   4.0810e-04 &  3.4805&&   1.5139e-03 &  3.0837         \\
$6$  &             2.5977e-04 &  3.7728&&   9.5542e-04 &  3.3594&&   2.7858e-05 &  3.9111&&   1.3019e-04 &  3.5747         \\
$7$  &             1.7035e-05 &  3.9115&&   5.1915e-05 &  4.1814&&   1.7626e-06 &  3.9629&&   5.6006e-06 &  4.5168         \\
\hline
\hline
& \multicolumn{2}{c}{MLS$^2_{\text{W4}}$} & &   \multicolumn{2}{c}{DD-MLS$^2_{\text{W4}}$}& &\multicolumn{2}{c}{MLS$^2_{\text{W4}}$} & &   \multicolumn{2}{c}{DD-MLS$^2_{\text{W4}}$}\\ \cline{1-3} \cline{5-6} \cline{8-9} \cline{11-12}
$l$ & $\text{MAE}_l$ & $r_l^{\infty}$ & &$\text{MAE}_l$ & $r_l^{\infty}$ & &$\text{RMSE}_l$ & $r_l^2$ & &$\text{RMSE}_l$ & $r_l^2$       \\
\hline
$4$   &          2.1519e-02 &        &&   4.7609e-02  &       &&   3.0906e-03 &        &&   9.0872e-03 &                 \\
$5$   &          2.2812e-03 &  3.2538&&   6.5667e-03  & 2.8721&&   2.6332e-04 &  3.5707&&   1.0050e-03 &  3.1924        \\
$6$   &          1.6727e-04 &  3.8069&&   6.3657e-04  & 3.4002&&   1.7560e-05 &  3.9451&&   8.5276e-05 &  3.5941        \\
$7$   &          1.0846e-05 &  3.9276&&   3.3648e-05  & 4.2210&&   1.1183e-06 &  3.9536&&   3.6034e-06 &  4.5425        \\
\hline
\hline
& \multicolumn{2}{c}{MLS$^2_{\text{G}}$} & &   \multicolumn{2}{c}{DD-MLS$^2_{\text{G}}$}& &\multicolumn{2}{c}{MLS$^2_{\text{G}}$} & &
\multicolumn{2}{c}{DD-MLS$^2_{\text{G}}$}\\ \cline{1-3} \cline{5-6} \cline{8-9} \cline{11-12}
$l$ & $\text{MAE}_l$ & $r_l^{\infty}$ & &$\text{MAE}_l$ & $r_l^{\infty}$ & &$\text{RMSE}_l$ & $r_l^2$ & &$\text{RMSE}_l$ & $r_l^2$       \\
 $4$   &      1.1701e-02  &       &&   3.0202e-02 &        &&   1.5423e-03 &        &&  4.9530e-03 &              \\
 $5$   &      1.0906e-03  & 3.4403&&   3.4672e-03 &  3.1383&&   1.2017e-04 &  3.7001&&  4.8617e-04 &  3.3654     \\
 $6$   &      7.8215e-05  & 3.8392&&   3.1008e-04 &  3.5176&&   8.0812e-06 &  3.9330&&  3.9902e-05 &  3.6427     \\
 $7$   &      5.2402e-06  & 3.8807&&   1.5806e-05 &  4.2732&&   5.6088e-07 &  3.8300&&  1.6869e-06 &  4.5417     \\                                                                                                          \hline
\end{tabular}
\end{center}
\caption{Errors and rates using linear and data dependent MLS methods for Franke's test function evaluated at grid points.}\label{tabla:d2r}
\end{table}

\begin{table}[!ht]
\begin{center}
\begin{tabular}{lrrrrrrrrrrr}
& \multicolumn{2}{c}{MLS$^2_{\text{W2}}$} & &   \multicolumn{2}{c}{DD-MLS$^2_{\text{W2}}$}& &\multicolumn{2}{c}{MLS$^2_{\text{W2}}$} & &   \multicolumn{2}{c}{DD-MLS$^2_{\text{W2}}$}\\ \cline{1-3} \cline{5-6} \cline{8-9} \cline{11-12}
$l$ & $\text{MAE}_l$ & $r_l^{\infty}$ & &$\text{MAE}_l$ & $r_l^{\infty}$ & &$\text{RMSE}_l$ & $r_l^2$ & &$\text{RMSE}_l$ & $r_l^2$       \\
\hline
$4$  &               3.0411e-02 &         &&  5.3557e-02 &        &&   4.5092e-03 &        &&   1.2182e-02 &                        \\
$5$  &               3.9097e-03 &  3.3987 &&  9.7501e-03 &  2.8224&&   5.1659e-04 &  3.5898&&   1.5725e-03 &  3.3920            \\
$6$  &               4.9130e-04 &  3.4450 &&  1.0768e-03 &  3.6595&&   5.4169e-05 &  3.7456&&   1.4011e-04 &  4.0161            \\
$7$  &               6.2966e-05 &  3.4057 &&  1.3020e-04 &  3.5022&&   6.7078e-06 &  3.4627&&   9.2877e-06 &  4.4986            \\
\hline
\hline
& \multicolumn{2}{c}{MLS$^2_{\text{W4}}$} & &   \multicolumn{2}{c}{DD-MLS$^2_{\text{W4}}$}& &\multicolumn{2}{c}{MLS$^2_{\text{W4}}$} & &   \multicolumn{2}{c}{DD-MLS$^2_{\text{W4}}$}\\ \cline{1-3} \cline{5-6} \cline{8-9} \cline{11-12}
$l$ & $\text{MAE}_l$ & $r_l^{\infty}$ & &$\text{MAE}_l$ & $r_l^{\infty}$ & &$\text{RMSE}_l$ & $r_l^2$ & &$\text{RMSE}_l$ & $r_l^2$       \\
\hline
$4$   &            2.3547e-02 &        &&   4.3721e-02 &        &&   3.3430e-03 &        &&   8.6294e-03 &              \\
$5$   &            2.7803e-03 &  3.5398&&   6.7029e-03 &  3.1071&&   3.8656e-04 &  3.5744&&   1.0595e-03 &  3.4750     \\
$6$   &            4.1415e-04 &  3.1625&&   7.9889e-04 &  3.5329&&   4.3740e-05 &  3.6192&&   9.6392e-05 &  3.9815     \\
$7$   &            6.4247e-05 &  3.0891&&   1.0343e-04 &  3.3890&&   5.8562e-06 &  3.3333&&   7.2394e-06 &  4.2916     \\
\hline
\hline
& \multicolumn{2}{c}{MLS$^2_{\text{G}}$} & &   \multicolumn{2}{c}{DD-MLS$^2_{\text{G}}$}& &\multicolumn{2}{c}{MLS$^2_{\text{G}}$} & &
\multicolumn{2}{c}{DD-MLS$^2_{\text{G}}$}\\ \cline{1-3} \cline{5-6} \cline{8-9} \cline{11-12}
$l$ & $\text{MAE}_l$ & $r_l^{\infty}$ & &$\text{MAE}_l$ & $r_l^{\infty}$ & &$\text{RMSE}_l$ & $r_l^2$ & &$\text{RMSE}_l$ & $r_l^2$       \\
 $4$   &       1.8162e-02 &        &&   3.1263e-02  &       &&   2.2047e-03 &        &&   4.8132e-03 &             \\
 $5$   &       2.2788e-03 &  3.4391&&   4.5951e-03  & 3.1769&&   2.7608e-04 &  3.4424&&   5.6946e-04 &  3.5364    \\
 $6$   &       3.8683e-04 &  2.9455&&   5.4298e-04  & 3.5472&&   3.3817e-05 &  3.4875&&   5.4933e-05 &  3.8842    \\
 $7$   &       5.6293e-05 &  3.1951&&   8.6890e-05  & 3.0376&&   4.8183e-06 &  3.2301&&   5.2895e-06 &  3.8797    \\                                                                                                          \hline
\end{tabular}
\end{center}
\caption{Errors and rates using linear and data dependent MLS methods for Franke's test function evaluated at Halton's points.}\label{tabla:d2h}
\end{table}

\begin{table}[!ht]
\begin{center}
\begin{tabular}{lrrrrrrrrrrr}
& \multicolumn{2}{c}{MLS$^1_{\text{W2}}$} & &   \multicolumn{2}{c}{DD-MLS$^1_{\text{W2}}$}& &\multicolumn{2}{c}{MLS$^1_{\text{W2}}$} & &   \multicolumn{2}{c}{DD-MLS$^1_{\text{W2}}$}\\ \cline{1-3} \cline{5-6} \cline{8-9} \cline{11-12}
$l$ & $\text{MAE}_l$ & $r_l^{\infty}$ & &$\text{MAE}_l$ & $r_l^{\infty}$ & &$\text{RMSE}_l$ & $r_l^2$ & &$\text{RMSE}_l$ & $r_l^2$       \\
\hline
$4$  &              1.1379e-01  &       &&   1.9511e-01 &        &&   2.9208e-02 &        &&   4.1874e-02  &              \\
$5$  &              3.1759e-02  & 1.8502&&   5.4100e-02 &  1.8597&&   8.2129e-03 &  1.8395&&   9.5939e-03  & 2.1364   \\
$6$  &              8.5522e-03  & 1.9115&&   1.1477e-02 &  2.2591&&   2.1242e-03 &  1.9703&&   2.1849e-03  & 2.1557   \\
$7$  &              2.2003e-03  & 1.9490&&   2.3391e-03 &  2.2835&&   5.3588e-04 &  1.9773&&   5.3824e-04  & 2.0114   \\
\hline
\hline
& \multicolumn{2}{c}{MLS$^1_{\text{W4}}$} & &   \multicolumn{2}{c}{DD-MLS$^1_{\text{W4}}$}& &\multicolumn{2}{c}{MLS$^1_{\text{W4}}$} & &   \multicolumn{2}{c}{DD-MLS$^1_{\text{W4}}$}\\ \cline{1-3} \cline{5-6} \cline{8-9} \cline{11-12}
$l$ & $\text{MAE}_l$ & $r_l^{\infty}$ & &$\text{MAE}_l$ & $r_l^{\infty}$ & &$\text{RMSE}_l$ & $r_l^2$ & &$\text{RMSE}_l$ & $r_l^2$       \\
\hline
$4$   &       8.9718e-02  &       &&   1.5742e-01 &        &&   2.3062e-02 &        &&   3.2179e-02 &             \\
$5$   &       2.4477e-02  & 1.8833&&   4.3323e-02 &  1.8707&&   6.3336e-03 &  1.8737&&   7.1698e-03 &  2.1769    \\
$6$   &       6.5848e-03  & 1.9130&&   8.3406e-03 &  2.4005&&   1.6265e-03 &  1.9807&&   1.6580e-03 &  2.1334    \\
$7$   &       1.6836e-03  & 1.9580&&   1.7678e-03 &  2.2273&&   4.0947e-04 &  1.9802&&   4.1084e-04 &  2.0030    \\
\hline
\hline
& \multicolumn{2}{c}{MLS$^1_{\text{G}}$} & &   \multicolumn{2}{c}{DD-MLS$^1_{\text{G}}$}& &\multicolumn{2}{c}{MLS$^1_{\text{G}}$} & &
\multicolumn{2}{c}{DD-MLS$^1_{\text{G}}$}\\ \cline{1-3} \cline{5-6} \cline{8-9} \cline{11-12}
$l$ & $\text{MAE}_l$ & $r_l^{\infty}$ & &$\text{MAE}_l$ & $r_l^{\infty}$ & &$\text{RMSE}_l$ & $r_l^2$ & &$\text{RMSE}_l$ & $r_l^2$       \\
 $4$   &        5.5298e-02 &        &&   9.7249e-02 &        &&   1.4210e-02 &        &&   1.8864e-02 &           \\
 $5$   &        1.4845e-02 &  1.9067&&   2.5712e-02 &  1.9288&&   3.7812e-03 &  1.9195&&   4.0917e-03 &  2.2158  \\
 $6$   &        3.9269e-03 &  1.9375&&   4.6121e-03 &  2.5035&&   9.6173e-04 &  1.9947&&   9.7048e-04 &  2.0965  \\
 $7$   &        9.9534e-04 &  1.9705&&   1.0290e-03 &  2.1536&&   2.4153e-04 &  1.9837&&   2.4203e-04 &  1.9938  \\
\hline
\end{tabular}
\end{center}
\caption{Errors and rates using linear and data dependent MLS methods for Franke's test function evaluated at grid points.}\label{tabla:d1r}
\end{table}

\begin{table}[!ht]
\begin{center}
\begin{tabular}{lrrrrrrrrrrr}
& \multicolumn{2}{c}{MLS$^1_{\text{W2}}$} & &   \multicolumn{2}{c}{DD-MLS$^1_{\text{W2}}$}& &\multicolumn{2}{c}{MLS$^1_{\text{W2}}$} & &   \multicolumn{2}{c}{DD-MLS$^1_{\text{W2}}$}\\ \cline{1-3} \cline{5-6} \cline{8-9} \cline{11-12}
$l$ & $\text{MAE}_l$ & $r_l^{\infty}$ & &$\text{MAE}_l$ & $r_l^{\infty}$ & &$\text{RMSE}_l$ & $r_l^2$ & &$\text{RMSE}_l$ & $r_l^2$       \\
\hline
$4$  &              1.0485e-01  &         & &   1.9775e-01  &         & &   2.8246e-02  &         & &   4.2458e-02  &                 \\
$5$  &              3.4265e-02  &  1.8530 & &   4.9345e-02  &  2.3000 & &   8.2717e-03  &  2.0348 & &   9.6733e-03  &  2.4507     \\
$6$  &              9.7104e-03  &  2.0943 & &   1.1974e-02  &  2.3520 & &   2.1446e-03  &  2.2420 & &   2.2221e-03  &  2.4431     \\
$7$  &              2.4667e-03  &  2.2716 & &   2.4842e-03  &  2.6073 & &   5.4304e-04  &  2.2769 & &   5.4552e-04  &  2.3282     \\
\hline
\hline
& \multicolumn{2}{c}{MLS$^1_{\text{W4}}$} & &   \multicolumn{2}{c}{DD-MLS$^1_{\text{W4}}$}& &\multicolumn{2}{c}{MLS$^1_{\text{W4}}$} & &   \multicolumn{2}{c}{DD-MLS$^1_{\text{W4}}$}\\ \cline{1-3} \cline{5-6} \cline{8-9} \cline{11-12}
$l$ & $\text{MAE}_l$ & $r_l^{\infty}$ & &$\text{MAE}_l$ & $r_l^{\infty}$ & &$\text{RMSE}_l$ & $r_l^2$ & &$\text{RMSE}_l$ & $r_l^2$       \\
\hline
$4$   &              8.3978e-02 &        &&   1.5293e-01  &       &&  2.2556e-02 &        &&   3.2350e-02 &                \\
$5$   &              2.8296e-02 &  1.8024&&   3.7951e-02  & 2.3092&&  6.4277e-03 &  2.0800&&   7.2640e-03 &  2.4748       \\
$6$   &              8.1573e-03 &  2.0659&&   8.9301e-03  & 2.4031&&  1.6564e-03 &  2.2522&&   1.6989e-03 &  2.4132       \\
$7$   &              2.1030e-03 &  2.2471&&   2.0241e-03  & 2.4606&&  4.1929e-04 &  2.2774&&   4.1958e-04 &  2.3183       \\
\hline
\hline
& \multicolumn{2}{c}{MLS$^1_{\text{G}}$} & &   \multicolumn{2}{c}{DD-MLS$^2_{\text{G}}$}& &\multicolumn{2}{c}{MLS$^1_{\text{G}}$} & &
\multicolumn{2}{c}{DD-MLS$^1_{\text{G}}$}\\ \cline{1-3} \cline{5-6} \cline{8-9} \cline{11-12}
$l$ & $\text{MAE}_l$ & $r_l^{\infty}$ & &$\text{MAE}_l$ & $r_l^{\infty}$ & &$\text{RMSE}_l$ & $r_l^2$ & &$\text{RMSE}_l$ & $r_l^2$       \\
 $4$   &       5.7777e-02 &        &&   8.1826e-02 &        &&  1.4432e-02 &        &&  1.8700e-02 &           \\
 $5$   &       2.1624e-02 &  1.6284&&   2.5200e-02 &  1.9514&&  3.9649e-03 &  2.1406&&  4.2732e-03 &  2.4458  \\
 $6$   &       6.1526e-03 &  2.0876&&   6.3276e-03 &  2.2953&&  1.0226e-03 &  2.2508&&  1.0340e-03 &  2.3568  \\
 $7$   &       1.5647e-03 &  2.2698&&   1.5688e-03 &  2.3118&&  2.5971e-04 &  2.2719&&  2.5859e-04 &  2.2975  \\
\hline
\end{tabular}
\end{center}
\caption{Errors and rates using linear and data dependent MLS methods for Franke's test function evaluated at Halton's points.}\label{tabla:d1h}
\end{table}

\begin{table}[!ht]
\begin{center}
\begin{tabular}{lrrrrrrrrrrr}
& \multicolumn{2}{c}{MLS$^0_{\text{W2}}$} & &   \multicolumn{2}{c}{DD-MLS$^0_{\text{W2}}$}& &\multicolumn{2}{c}{MLS$^0_{\text{W2}}$} & &     \multicolumn{2}{c}{DD-MLS$^0_{\text{W2}}$}\\ \cline{1-3} \cline{5-6} \cline{8-9} \cline{11-12}
$l$ & $\text{MAE}_l$ & $r_l^{\infty}$ & &$\text{MAE}_l$ & $r_l^{\infty}$ & &$\text{RMSE}_l$ & $r_l^2$ & &$\text{RMSE}_l$ & $r_l^2$       \\
\hline
$4$  &            1.1379e-01  &       &&   2.1215e-01 &        &&   2.9355e-02  &       &&   7.9470e-02 &               \\
$5$  &            3.1754e-02  & 1.8504&&   9.1316e-02 &  1.2222&&   8.3385e-03  & 1.8247&&   3.0115e-02 &  1.4069   \\
$6$  &            8.5514e-03  & 1.9115&&   3.1215e-02 &  1.5640&&   2.1290e-03  & 1.9891&&   8.2379e-03 &  1.8887   \\
$7$  &            2.2009e-03  & 1.9485&&   5.6983e-03 &  2.4417&&   5.3596e-04  & 1.9803&&   1.0663e-03 &  2.9352   \\
\hline
\hline
& \multicolumn{2}{c}{MLS$^0_{\text{W4}}$} & &   \multicolumn{2}{c}{DD-MLS$^0_{\text{W4}}$}& &\multicolumn{2}{c}{MLS$^0_{\text{W4}}$} & &     \multicolumn{2}{c}{DD-MLS$^0_{\text{W4}}$}\\ \cline{1-3} \cline{5-6} \cline{8-9} \cline{11-12}
$l$ & $\text{MAE}_l$ & $r_l^{\infty}$ & &$\text{MAE}_l$ & $r_l^{\infty}$ & &$\text{RMSE}_l$ & $r_l^2$ & &$\text{RMSE}_l$ & $r_l^2$       \\
\hline
$4$   &             8.9721e-02 &        &&   1.8258e-01 &        &&   2.3203e-02 &        &&   6.5426e-02 &                \\
$5$   &             2.4474e-02 &  1.8835&&   7.5082e-02 &  1.2884&&   6.3935e-03 &  1.8688&&   2.4098e-02 &  1.4481       \\
$6$   &             6.5832e-03 &  1.9132&&   2.4913e-02 &  1.6074&&   1.6274e-03 &  1.9936&&   6.4179e-03 &  1.9277       \\
$7$   &             1.6845e-03 &  1.9569&&   4.4363e-03 &  2.4773&&   4.0959e-04 &  1.9806&&   8.1827e-04 &  2.9570       \\
\hline
\hline
& \multicolumn{2}{c}{MLS$^0_{\text{G}}$} & &   \multicolumn{2}{c}{DD-MLS$^0_{\text{G}}$}& &\multicolumn{2}{c}{MLS$^0_{\text{G}}$} & &   \multicolumn{2}{c}{DD-MLS$^0_{\text{G}}$}\\ \cline{1-3} \cline{5-6} \cline{8-9} \cline{11-12}
$l$ & $\text{MAE}_l$ & $r_l^{\infty}$ & &$\text{MAE}_l$ & $r_l^{\infty}$ & &$\text{RMSE}_l$ & $r_l^2$ & &$\text{RMSE}_l$ & $r_l^2$       \\
\hline
$4$   &         5.5301e-02 &         &  &   1.3307e-01 &          & & 1.4256e-02 &         & &    4.4167e-02  &                \\
$5$   &         1.4844e-02 &   1.9068&  &   5.0079e-02 &   1.4169 & &  3.7889e-03 &   1.9212& &    1.5365e-02  &  1.5308       \\
$6$   &         3.9248e-03 &   1.9382&  &   1.5705e-02 &   1.6896 & &  9.6178e-04 &   1.9976& &    3.8797e-03  &  2.0053       \\
$7$   &         9.9655e-04 &   1.9679&  &   2.6784e-03 &   2.5393 & &  2.4169e-04 &   1.9828& &    4.8548e-04  &  2.9838       \\
\hline
\end{tabular}
\end{center}
\caption{Errors and rates using linear and data dependent MLS methods for Franke's test function evaluated at  grid points.}\label{tabla:d0r}
\end{table}

\begin{table}[!ht]
\begin{center}
\begin{tabular}{lrrrrrrrrrrr}
& \multicolumn{2}{c}{MLS$^0_{\text{W2}}$} & &   \multicolumn{2}{c}{DD-MLS$^0_{\text{W2}}$}& &\multicolumn{2}{c}{MLS$^0_{\text{W2}}$} & &   \multicolumn{2}{c}{DD-MLS$^0_{\text{W2}}$}\\ \cline{1-3} \cline{5-6} \cline{8-9} \cline{11-12}
$l$ & $\text{MAE}_l$ & $r_l^{\infty}$ & &$\text{MAE}_l$ & $r_l^{\infty}$ & &$\text{RMSE}_l$ & $r_l^2$ & &$\text{RMSE}_l$ & $r_l^2$       \\
\hline
$4$  &          1.0586e-01 &        &&   2.4268e-01  &       & &  2.8270e-02&         &&   8.4870e-02 &               \\
$5$  &          3.4724e-02 &  1.8469&&   9.8871e-02  & 1.4877& &  9.3858e-03&   1.8269&&   3.1069e-02 &  1.6650   \\
$6$  &          1.3981e-02 &  1.5109&&   3.3478e-02  & 1.7986& &  2.8053e-03&   2.0059&&   9.1451e-03 &  2.0313   \\
$7$  &          6.1264e-03 &  1.3678&&   9.2646e-03  & 2.1296& &  1.0559e-03&   1.6198&&   1.6218e-03 &  2.8673   \\
\hline
\hline
& \multicolumn{2}{c}{MLS$^0_{\text{W4}}$} & &   \multicolumn{2}{c}{DD-MLS$^0_{\text{W4}}$}& &\multicolumn{2}{c}{MLS$^0_{\text{W4}}$} & &   \multicolumn{2}{c}{DD-MLS$^0_{\text{W4}}$}\\ \cline{1-3} \cline{5-6} \cline{8-9} \cline{11-12}
$l$ & $\text{MAE}_l$ & $r_l^{\infty}$ & &$\text{MAE}_l$ & $r_l^{\infty}$ & &$\text{RMSE}_l$ & $r_l^2$ & &$\text{RMSE}_l$ & $r_l^2$       \\
\hline
$4$   &         8.2583e-02 &        &&   2.0157e-01 &        &&   2.2883e-02 &        &&   7.0341e-02 &              \\
$5$   &         3.3019e-02 &  1.5189&&   8.3508e-02 &  1.4600&&   7.8197e-03 &  1.7791&&   2.5145e-02 &  1.7044     \\
$6$   &         1.5863e-02 &  1.2176&&   3.1235e-02 &  1.6334&&   2.6272e-03 &  1.8116&&   7.4538e-03 &  2.0196     \\
$7$   &         6.7009e-03 &  1.4285&&   9.0153e-03 &  2.0599&&   1.1405e-03 &  1.3833&&   1.5213e-03 &  2.6344     \\
\hline
\hline
& \multicolumn{2}{c}{MLS$^0_{\text{G}}$} & &   \multicolumn{2}{c}{DD-MLS$^0_{\text{G}}$}& &\multicolumn{2}{c}{MLS$^0_{\text{G}}$} & &   \multicolumn{2}{c}{DD-MLS$^0_{\text{G}}$}\\ \cline{1-3} \cline{5-6} \cline{8-9} \cline{11-12}
$l$ & $\text{MAE}_l$ & $r_l^{\infty}$ & &$\text{MAE}_l$ & $r_l^{\infty}$ & &$\text{RMSE}_l$ & $r_l^2$ & &$\text{RMSE}_l$ & $r_l^2$       \\
\hline
$4$   &          7.5331e-02  &       &&   1.4497e-01  &       &&   1.6087e-02 &        &&   4.8153e-02 &              \\
$5$   &          3.2711e-02  & 1.3821&&   6.3294e-02  & 1.3731&&   6.3405e-03 &  1.5427&&   1.6858e-02 &  1.7389     \\
$6$   &          1.8551e-02  & 0.9421&&   2.8400e-02  & 1.3311&&   2.7718e-03 &  1.3743&&   5.3489e-03 &  1.9067     \\
$7$   &          8.8940e-03  & 1.2187&&   1.0852e-02  & 1.5948&&   1.4143e-03 &  1.1155&&   1.5756e-03 &  2.0262     \\
\hline
\end{tabular}
\end{center}
\caption{Errors and rates using linear and data dependent MLS methods for Franke's test function evaluated at Halton's points.}\label{tabla:d0h}
\end{table}

\subsection{Avoiding oscillations}

\begin{figure}[!ht]
\begin{center}
		\begin{tabular}{cc}
 \hspace{-1cm}				\includegraphics[width=7cm, height=5cm]{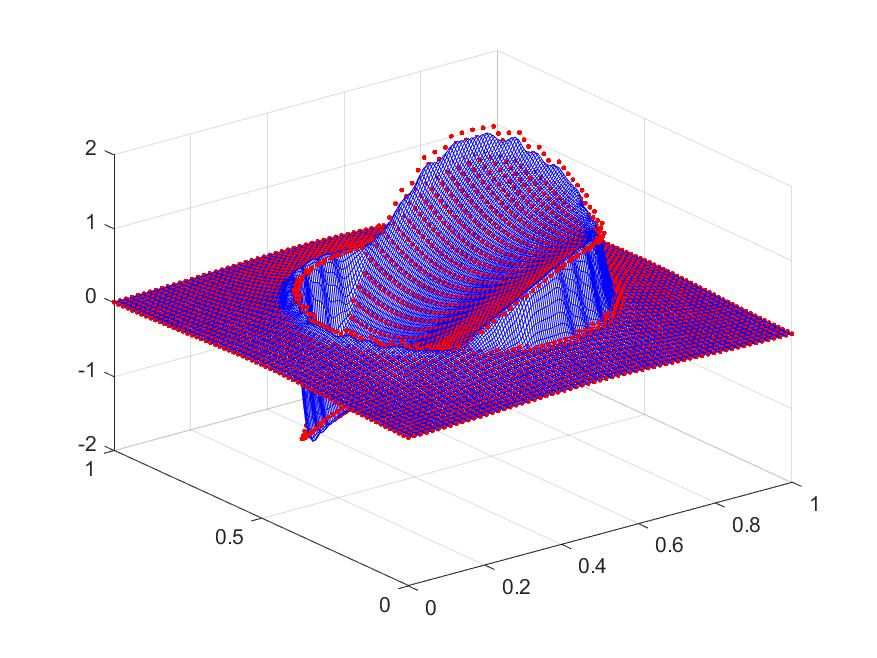} & \hspace{-1cm}
			\includegraphics[width=7cm, height=5cm]{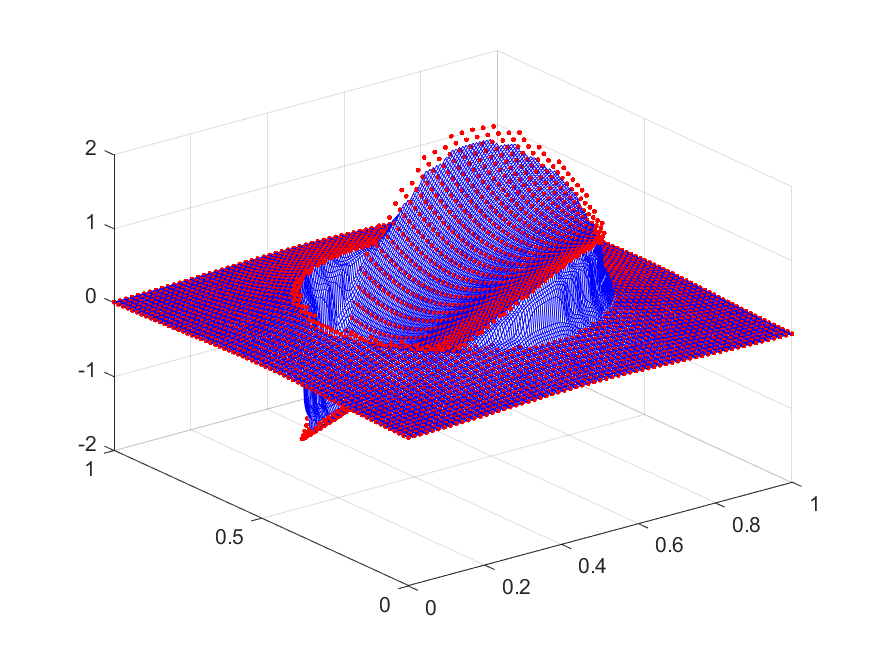}\\
 \hspace{-1cm}	
			\includegraphics[width=7cm, height=5cm]{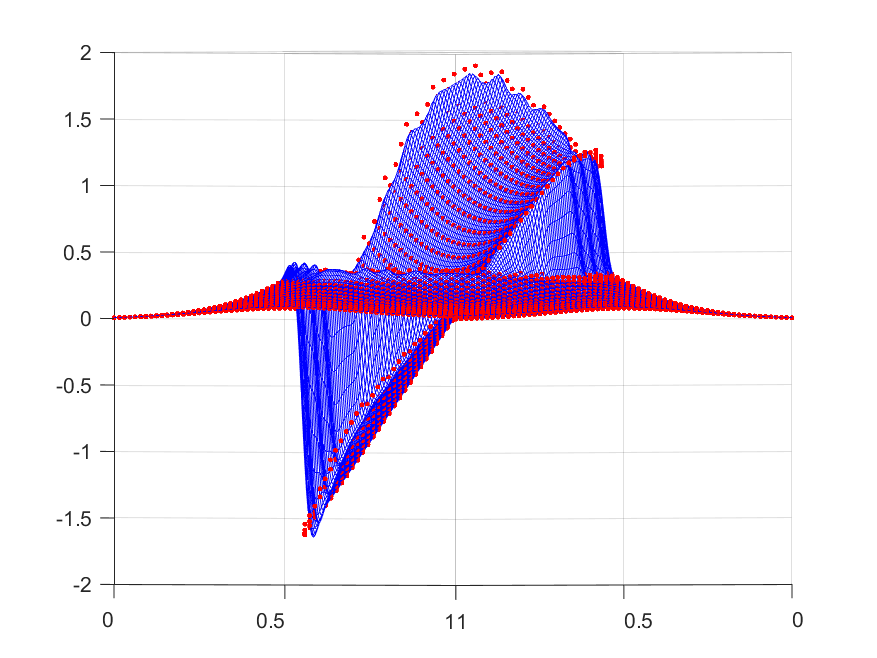} & \hspace{-1cm}
			\includegraphics[width=7cm, height=5cm]{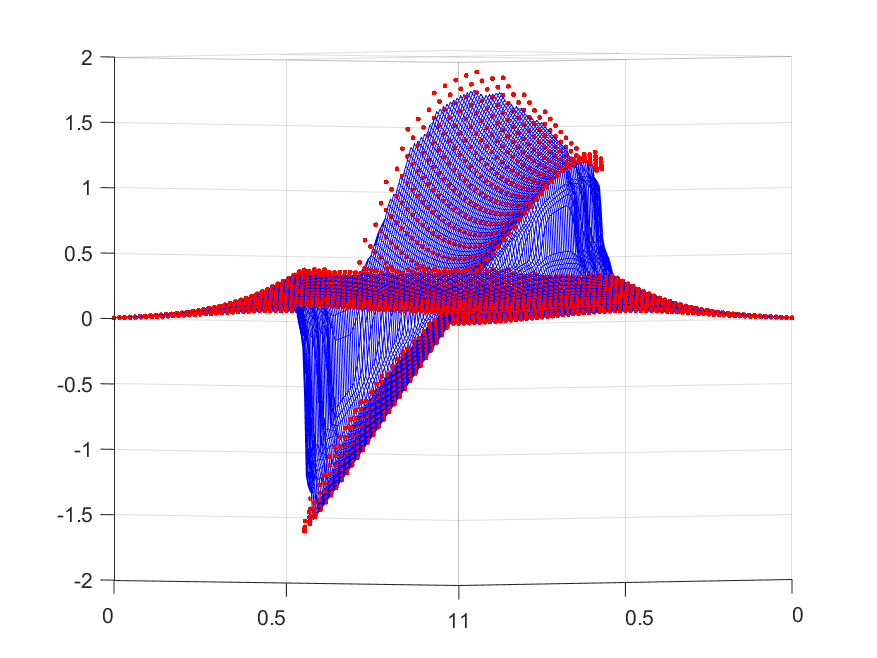} \\
	a)  MLS$^2_{\text{W2}}$ & \hspace{-2cm}	b) DD-MLS$^2_{\text{W2}}$\\
		 \hspace{-1cm}	\includegraphics[width=7cm, height=5cm]{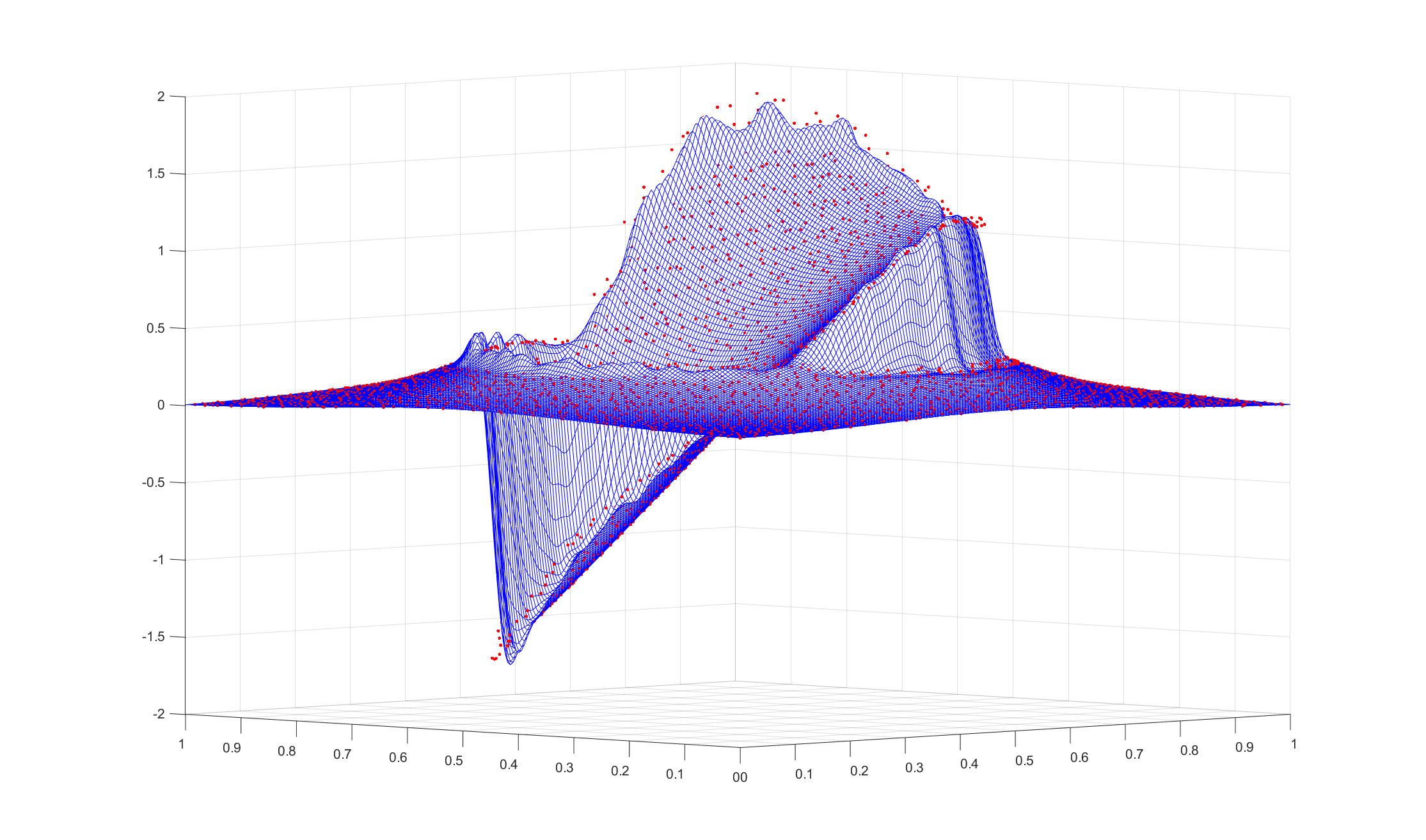} & \hspace{-1cm}	\includegraphics[width=7cm, height=5cm]{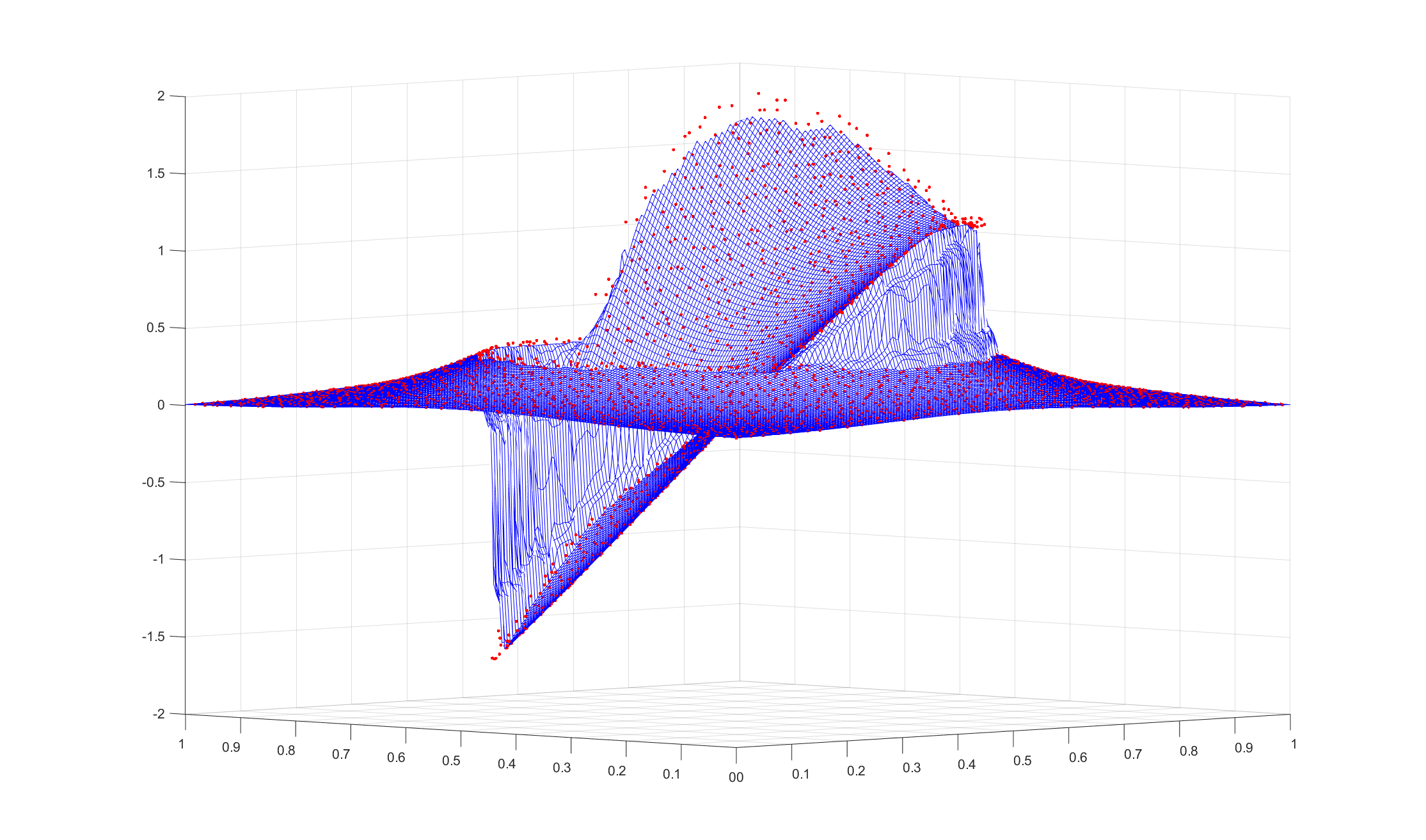}	\\
c) {MLS$^2_{\text{W2}}$} & d) {DD-MLS$^2_{\text{W2}}$}\\
\hspace{-1cm}\includegraphics[width=7cm, height=5cm]{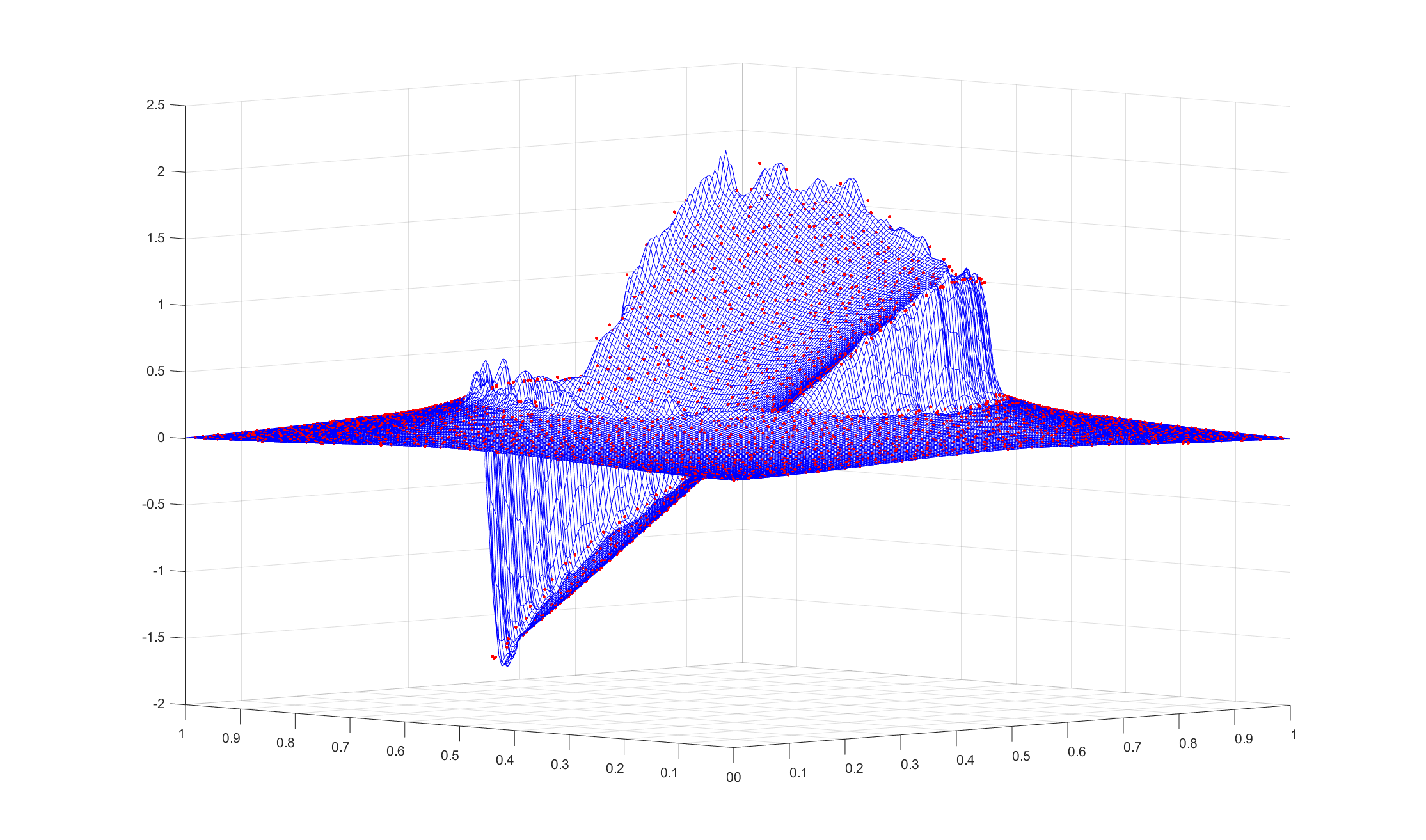} 	  &\hspace{-1cm} \includegraphics[width=7cm, height=5cm]{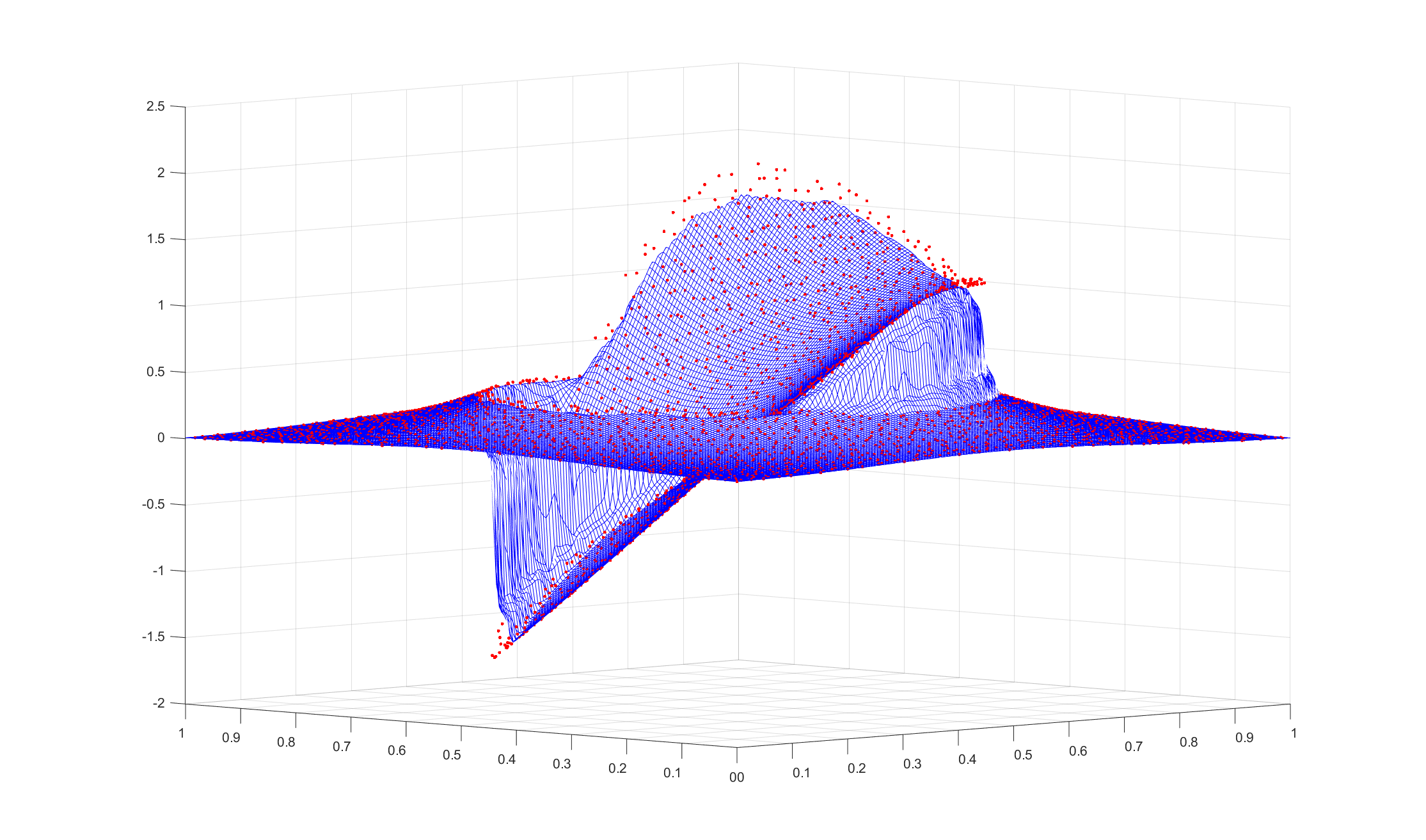} \\ e) MLS$^2_{\text{G}}$ & f) DD-MLS$^2_{\text{G}}$ 	
		\end{tabular}
\end{center}
\caption{Approximation to the function $g$, Eq. \eqref{ejemplolevin}, using linear and data dependent MLS with different $\omega(x)$ functions \cite{wendland2002} and the class of polynomials $\Pi_2(\mathbb{R}^2)$. Red points: original function, blue lines: approximation.}
	\label{figureejemplot1}
	\end{figure}

In this subsection, we approximate some functions with discontinuities as we have defined in Eq. \ref{funcion1}.  We start by approximating the function, $g$ on $[0,1]^2$, defined in \cite{AmirLevin}:
\begin{equation}\label{ejemplolevin}
g(x, y) =
\begin{cases}
-(x + y + 1) \cos(4x) + \sin(4(x + y)), & (x - 0.5)^2 + (y - 0.5)^2 \geq 0.1, \\
\exp\left(-10\left((x - 0.5)^2 + (y - 0.5)^2\right)\right), & (x - 0.5)^2 + (y - 0.5)^2 < 0.1,
\end{cases}
\end{equation}
using $N=65^2$ grid data points, Figure \ref{figureejemplot1}.a) {using the linear method}, and Figure \ref{figureejemplot1}.b), using W2. When a polynomial of degree $d\geq 2$ is used, we can see
that some oscillations appear close to the discontinuities, Figure \ref{figureejemplot1}.a), in the linear case.
This phenomenon is not avoided even if we refine the mesh. We can observe that these non-desired oscillations disappear when the data-dependent
method is employed, Figure \ref{figureejemplot1}.b). This result is very similar when the data points are pseudorandom, Figures \ref{figureejemplot1}.c)
and \ref{figureejemplot1}.d). If the approximator is $C^\infty$, Figures \ref{figureejemplot1}.e) and \ref{figureejemplot1}.f), then the result is similar, but in the DD-MLS we can observe that {some smearing of the discontinuities appears.}


%

In order to analyze the behaviour when the discontinuity is more  pronounced, we perform the example designed in \cite{amatetalnuma} with the following function
\begin{equation}\label{ejemplonuma}
z(x, y) =
\begin{cases}
\cos(xy), & (x - 0.5)^2 + (y - 0.5)^2 \geq 0.25^2, \\
\sin(xy), & (x - 0.5)^2 + (y - 0.5)^2 < 0.25^2.
\end{cases}
\end{equation}
In this experiment, Figures \ref{figureejemplonuma} and \ref{figureejemplonuma2}, we show the result using $N=65^2$ and $N=33^2$ data points, and the error between the original function and the approximated one. We can observe that the non-desired oscillations disappear for W2 and G cases.

\begin{landscape}
	\begin{figure}[!ht]
\begin{center}
		\begin{tabular}{cccc}
   \multicolumn{2}{c}{MLS$^2_{\text{W2}}$}& 	\multicolumn{2}{c}{DD-MLS$^2_{\text{W2}}$}\\
 \hspace{-1cm}				\includegraphics[width=5.6cm]{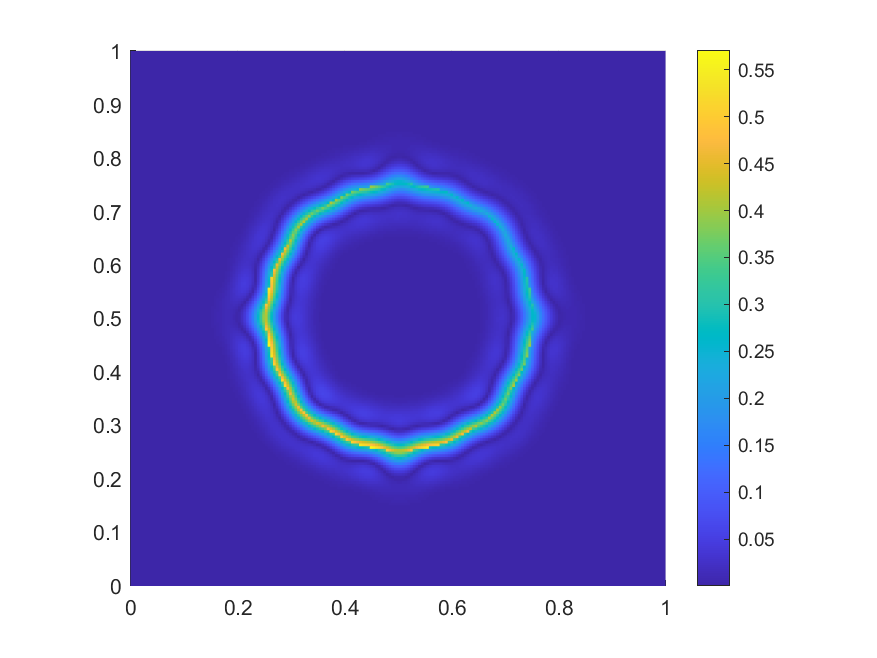} & \hspace{-1cm} \includegraphics[width=5.6cm]{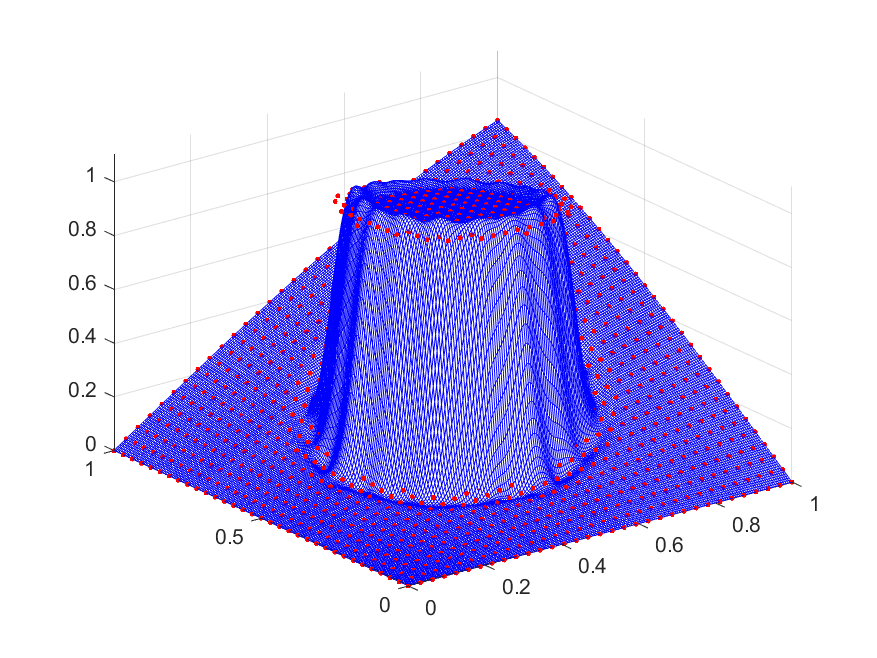}&
 \hspace{-1cm}				\includegraphics[width=5.6cm]{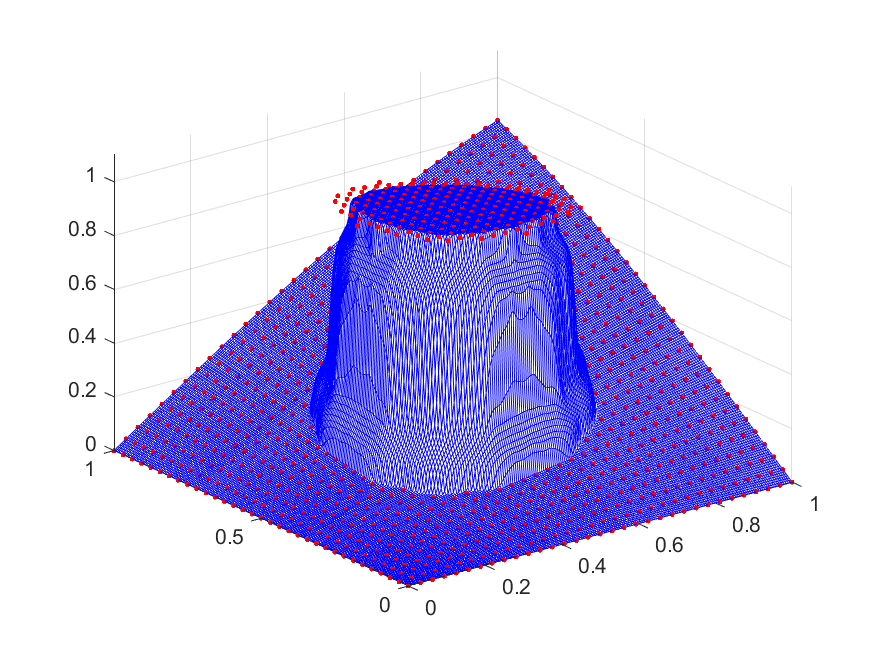} & \hspace{-1cm} \includegraphics[width=5.6cm]{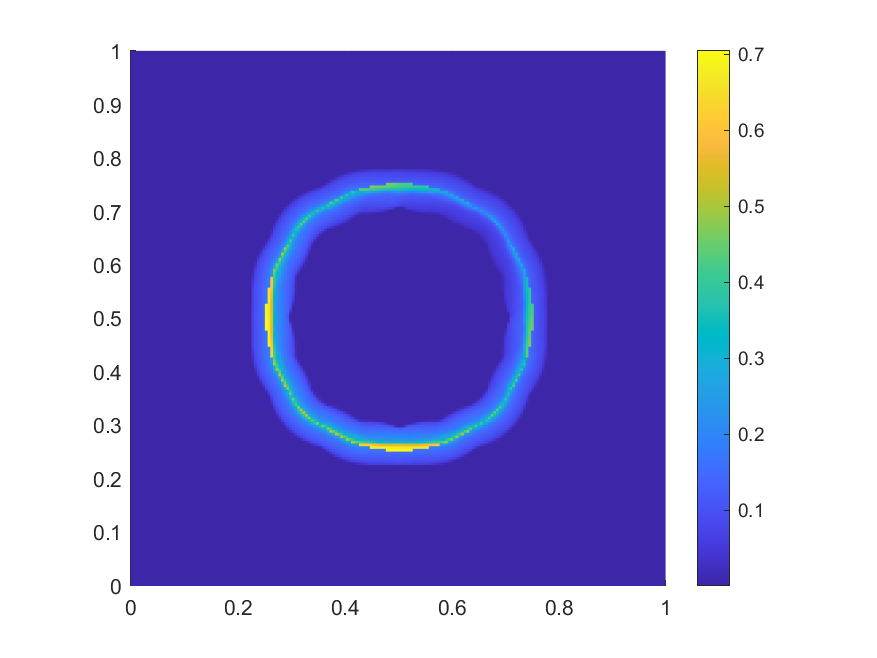}\\
 \hspace{-1cm}			\includegraphics[width=5.6cm]{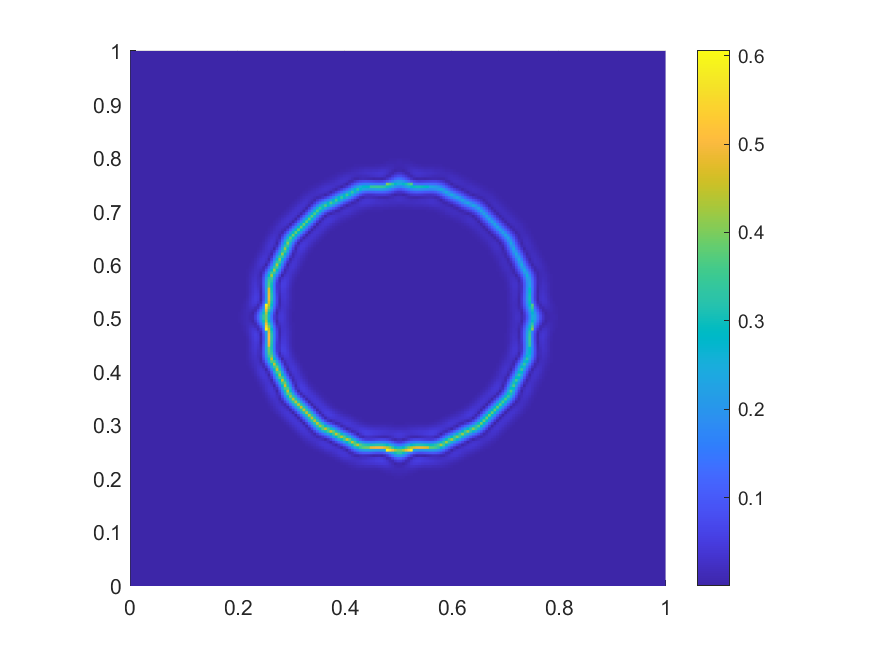} & \hspace{-1cm}			\includegraphics[width=5.6cm]{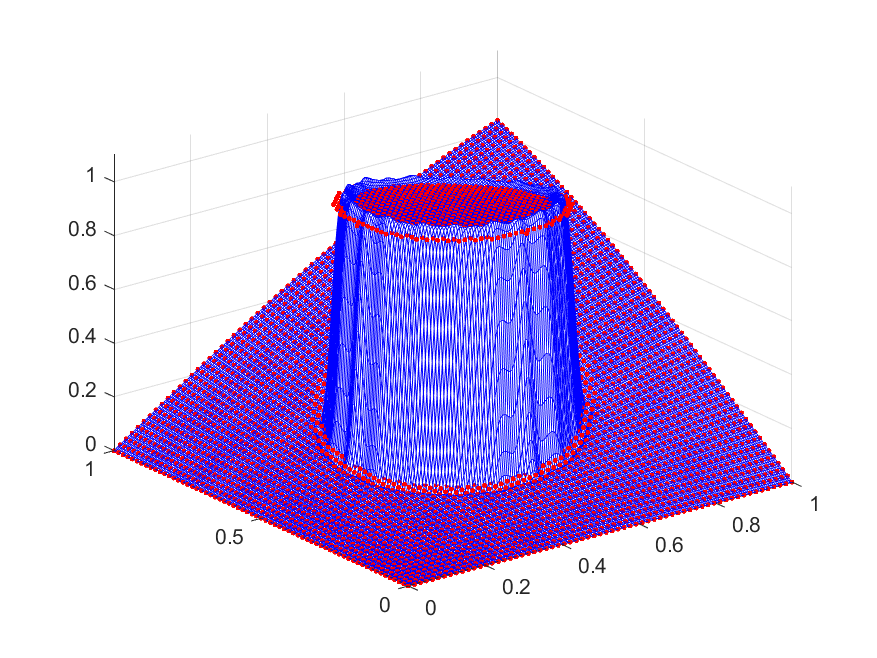} & \hspace{-1cm}			\includegraphics[width=5.6cm]{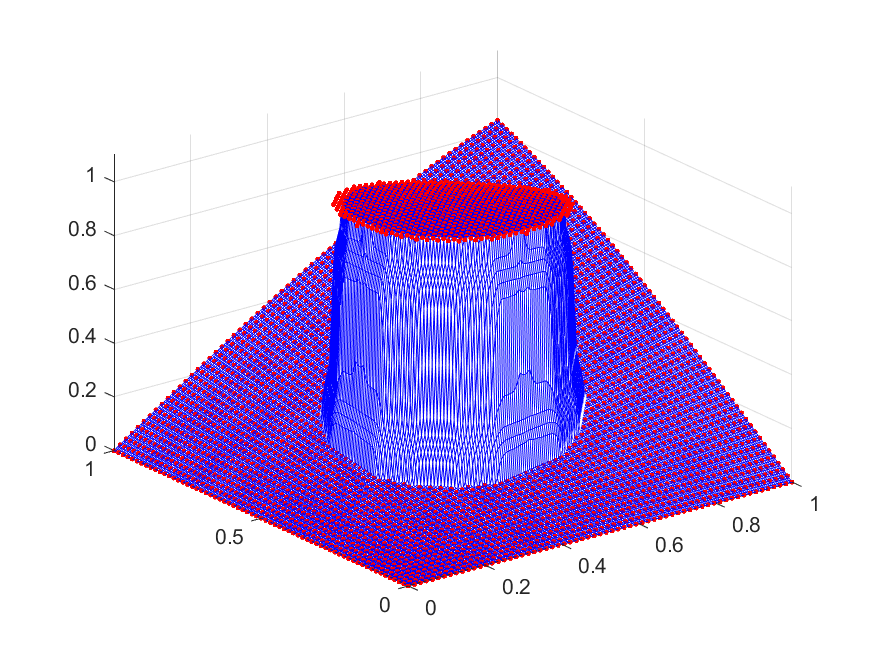} & \hspace{-1cm} \includegraphics[width=5.6cm]{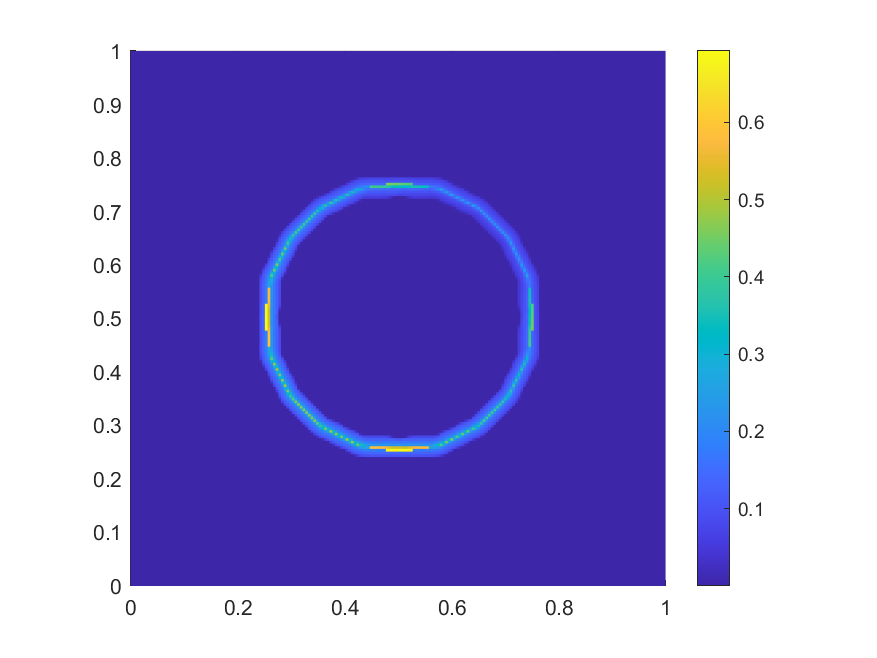}\\
\hspace{-1cm}				\includegraphics[width=5.6cm]{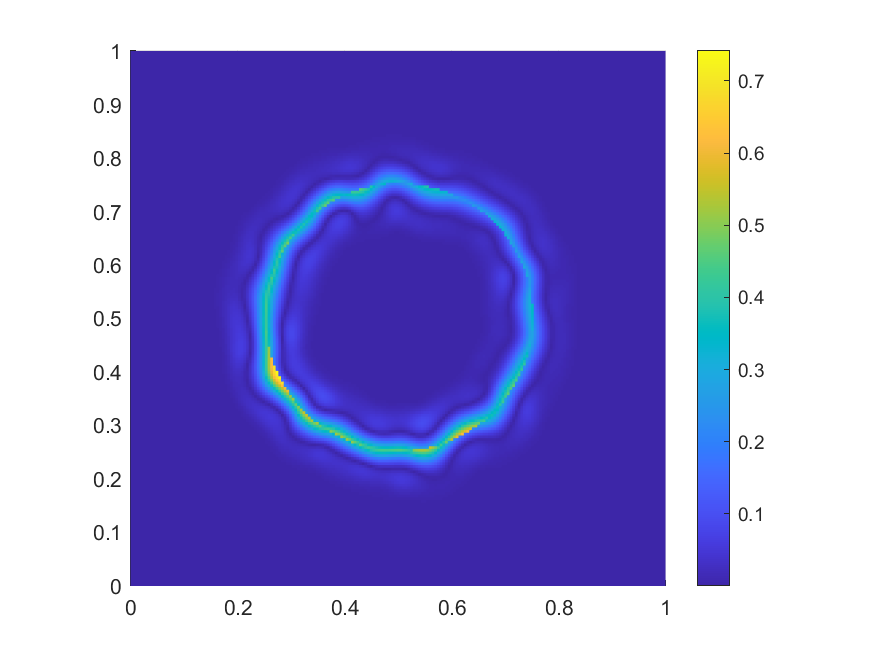} & \hspace{-1cm}\includegraphics[width=5.6cm]{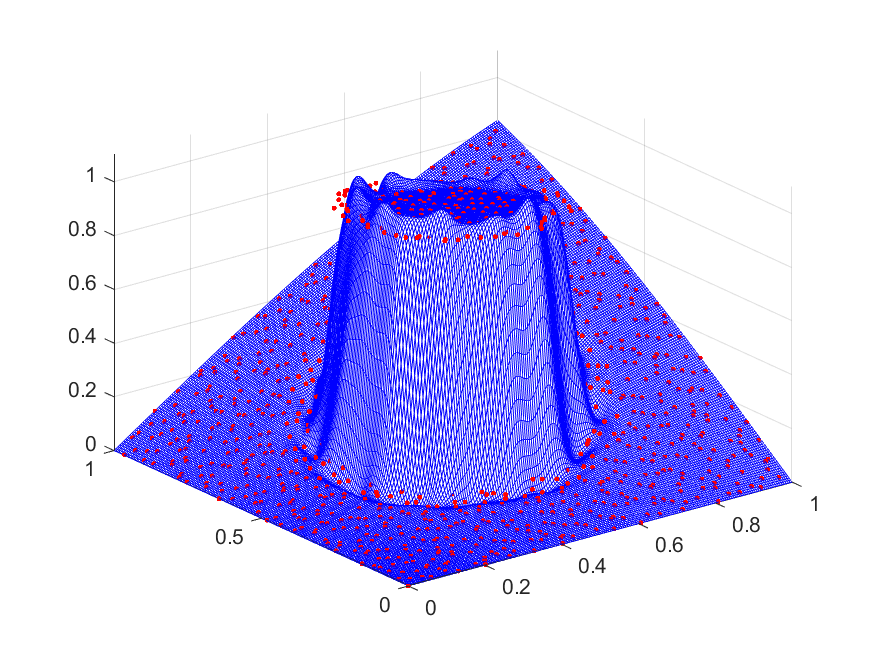} & \hspace{-1cm}			\includegraphics[width=5.6cm]{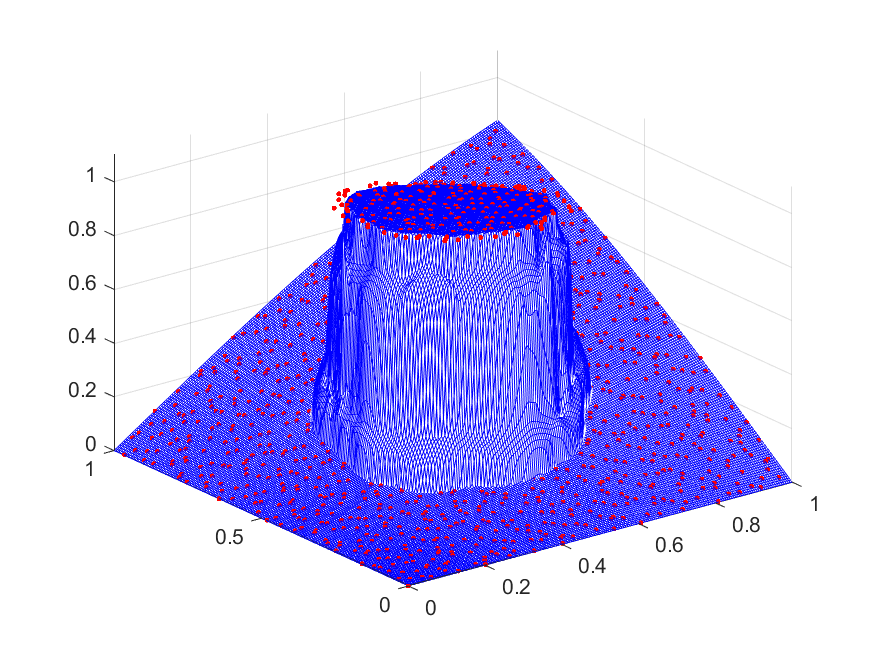} & \hspace{-1cm} \includegraphics[width=5.6cm]{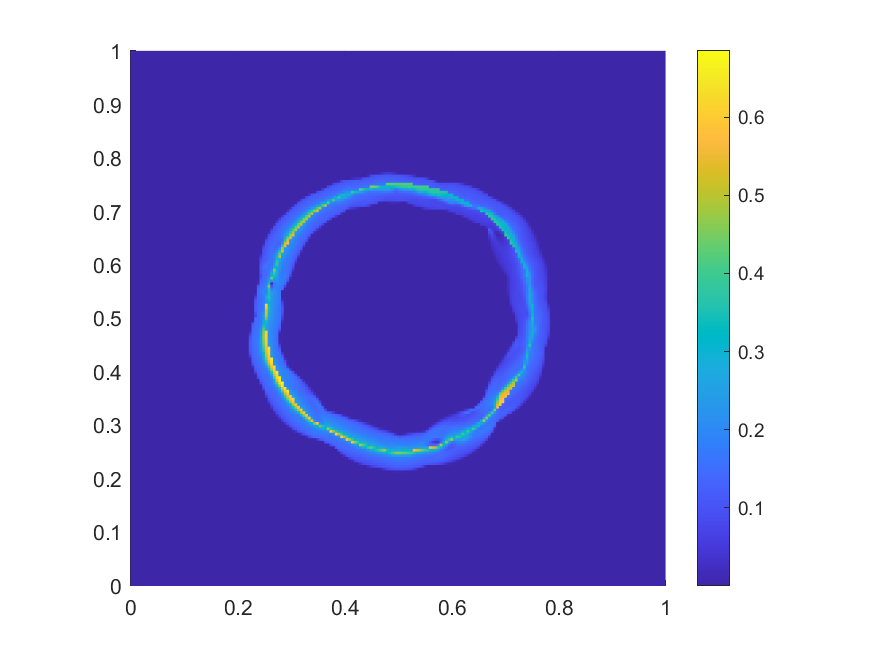}\\
\hspace{-1cm}	\includegraphics[width=5.6cm]{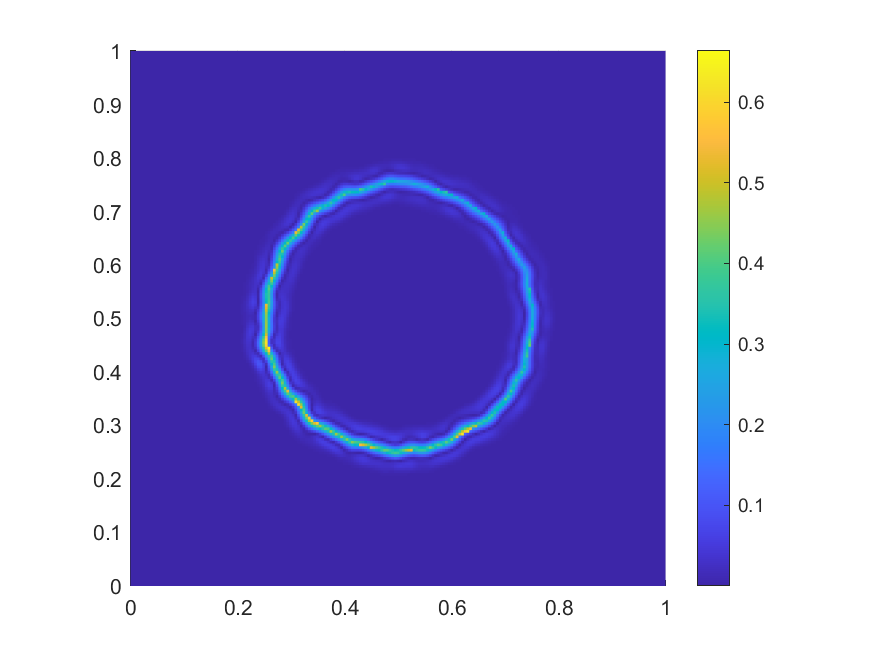} &  \hspace{-1cm}	\includegraphics[width=5.6cm]{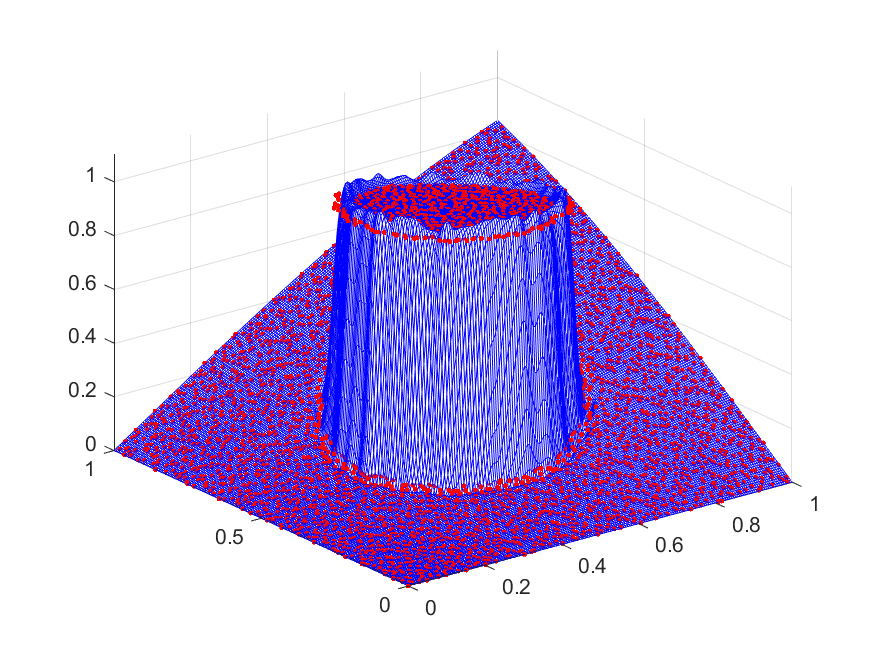} & \hspace{-1cm}
			\includegraphics[width=5.6cm]{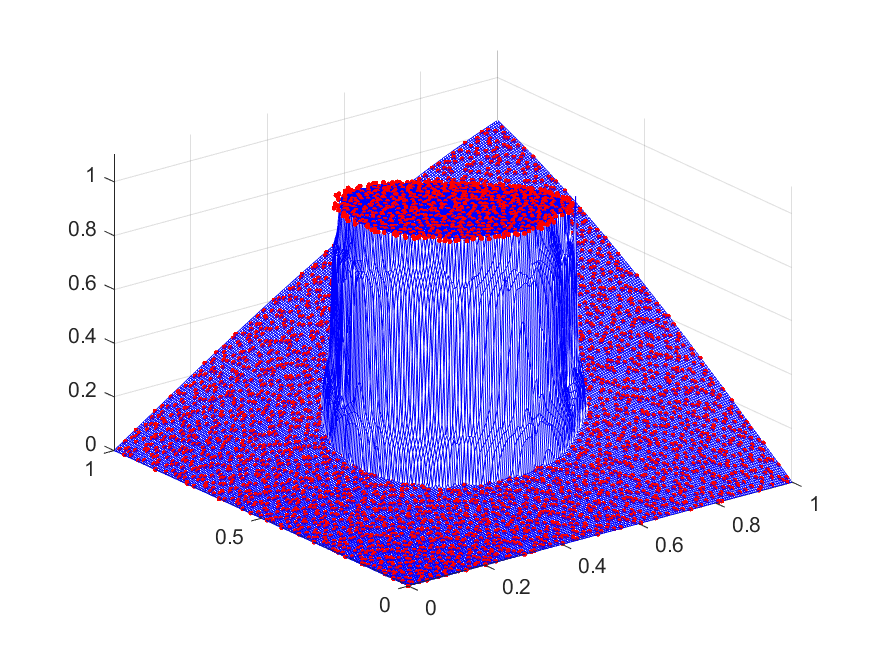} & \hspace{-1cm} \includegraphics[width=5.6cm]{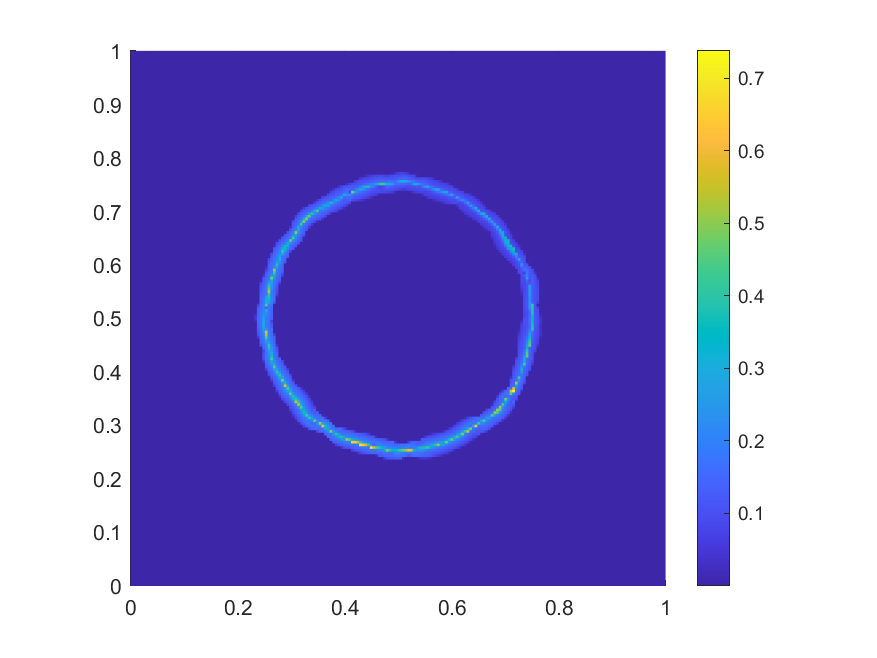}
		\end{tabular}
\end{center}
\caption{Approximation to the function $z$, Eq. \eqref{ejemplonuma}, using linear and data dependent MLS with $\omega(x)$ the $\mathcal{C}^2$ Wendland function \cite{wendland2002} and the class of polynomials $\Pi_2(\mathbb{R}^2)$. Red points: original function, blue lines: approximation.}
	\label{figureejemplonuma}
	\end{figure}

	\begin{figure}[!ht]
\begin{center}
		\begin{tabular}{cccc}
   \multicolumn{2}{c}{MLS$^2_{\text{G}}$}& 	\multicolumn{2}{c}{DD-MLS$^2_{\text{G}}$}\\
 \hspace{-1cm}				\includegraphics[width=5.6cm]{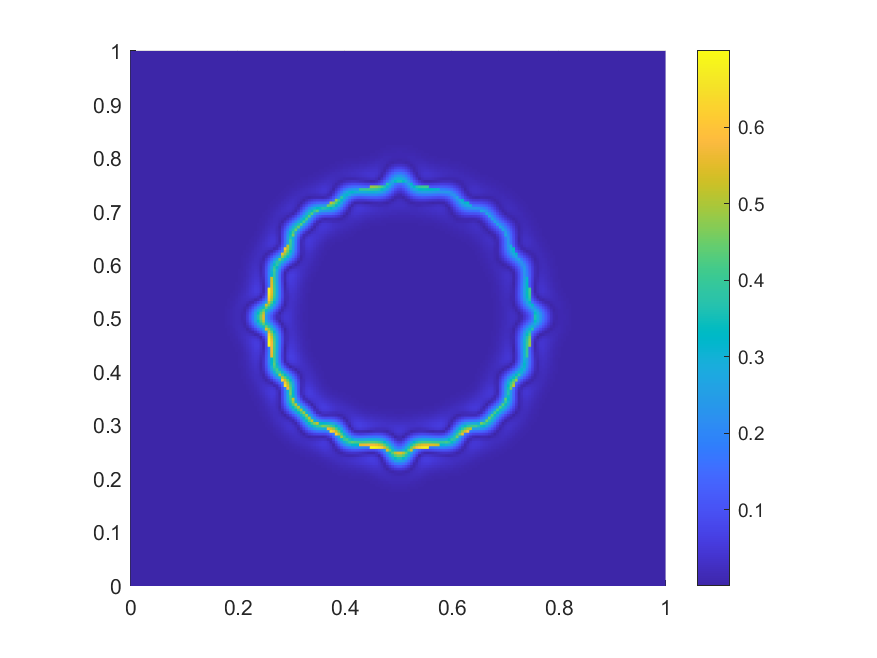} & \hspace{-1cm} \includegraphics[width=5.6cm]{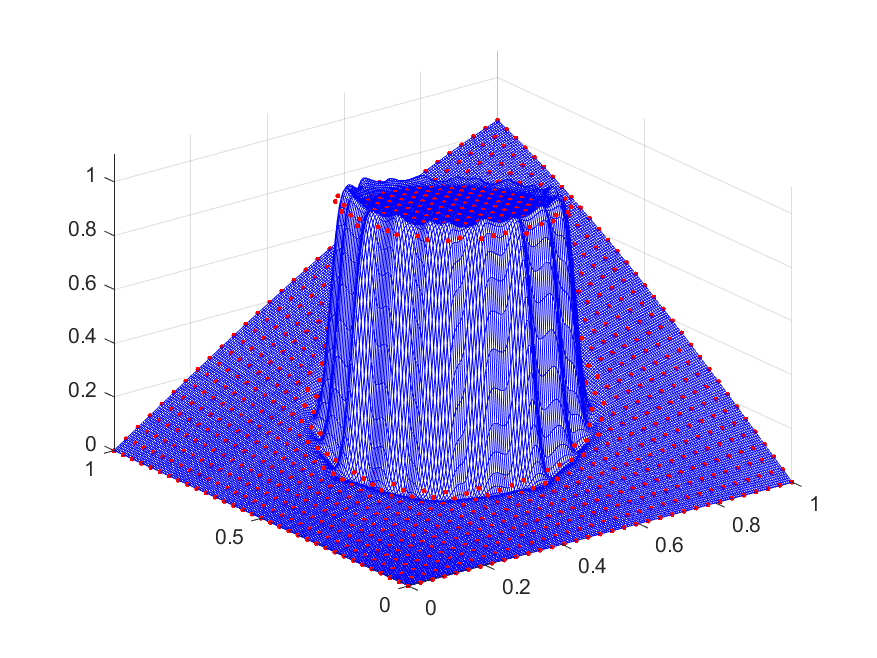}&
 \hspace{-1cm}				\includegraphics[width=5.6cm]{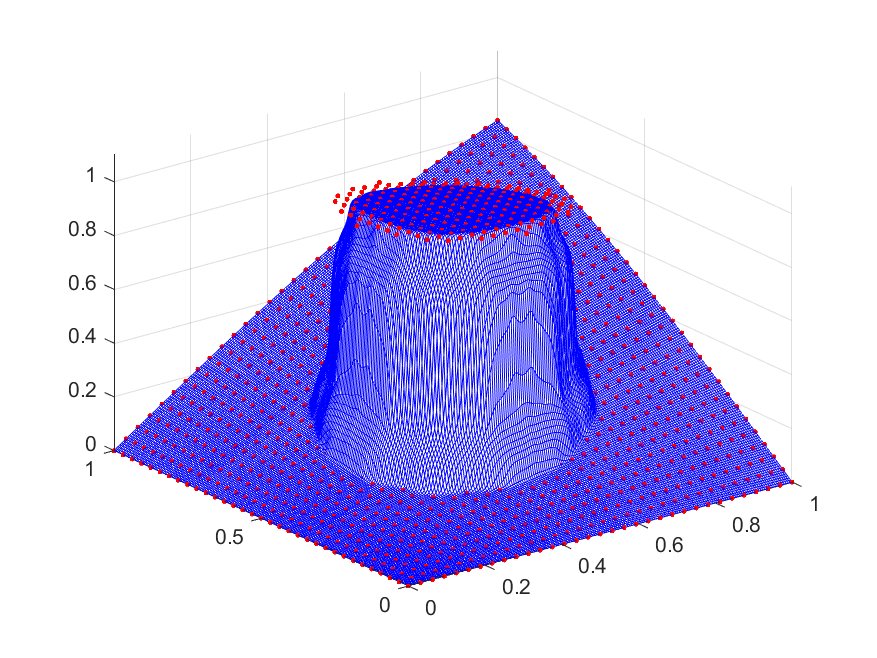} & \hspace{-1cm} \includegraphics[width=5.6cm]{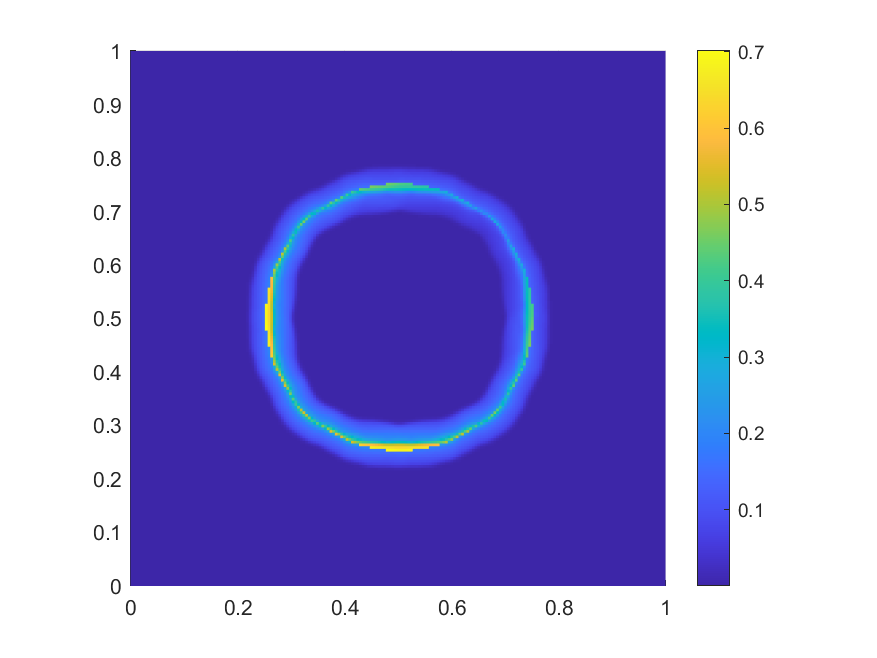}\\
 \hspace{-1cm}			\includegraphics[width=5.6cm]{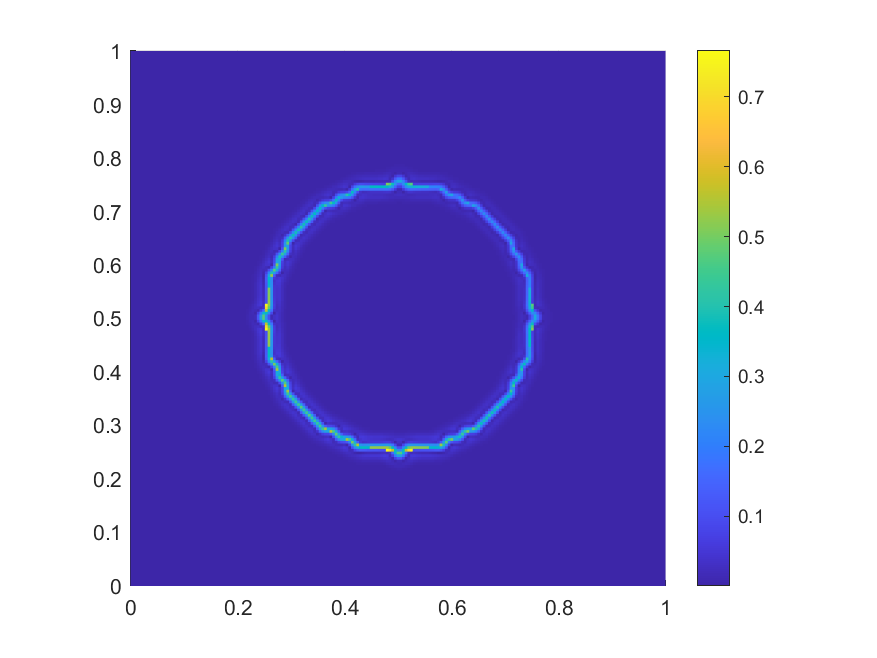} & \hspace{-1cm}			\includegraphics[width=5.6cm]{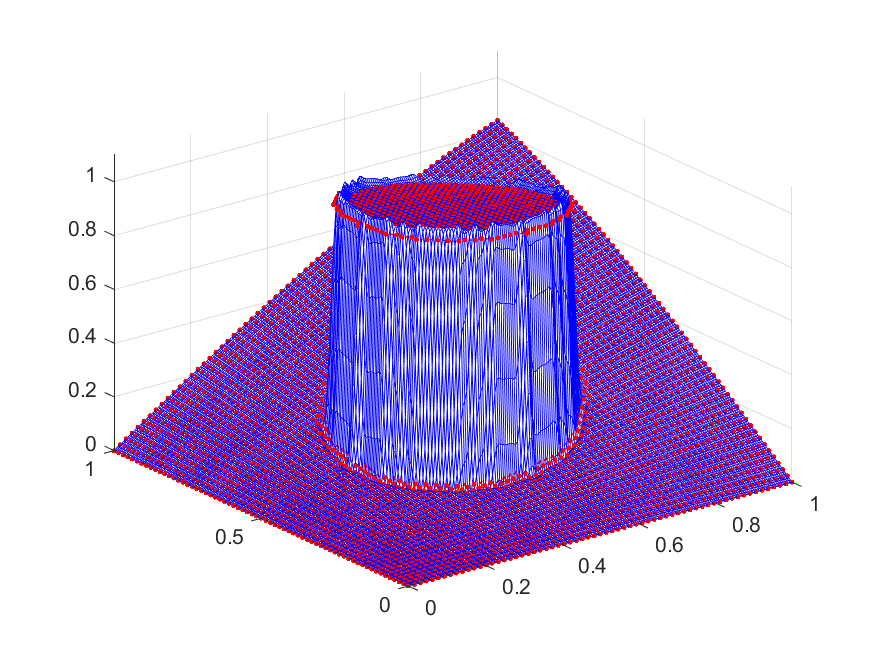} & \hspace{-1cm}			\includegraphics[width=5.6cm]{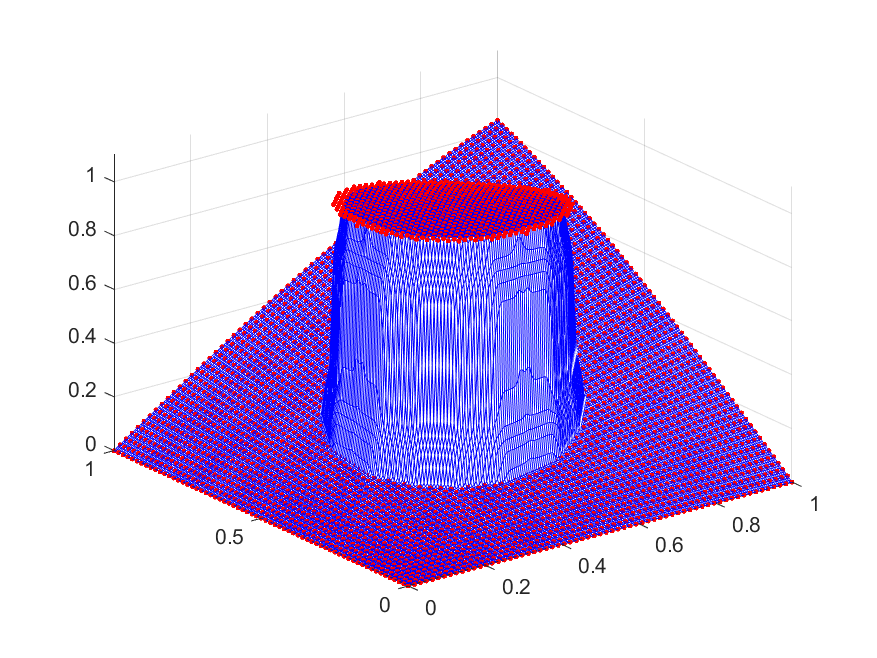} & \hspace{-1cm} \includegraphics[width=5.6cm]{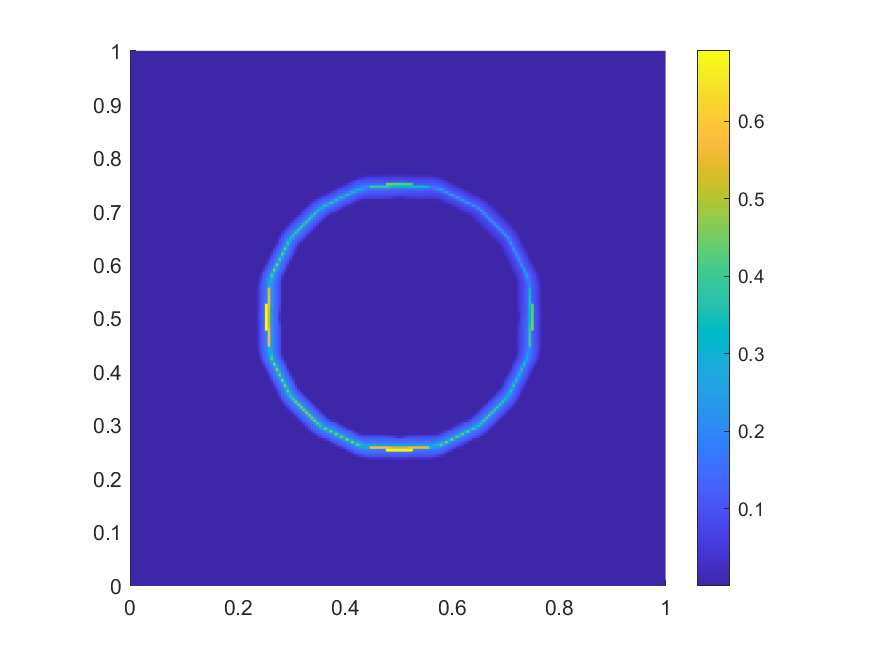}\\
\hspace{-1cm}				\includegraphics[width=5.6cm]{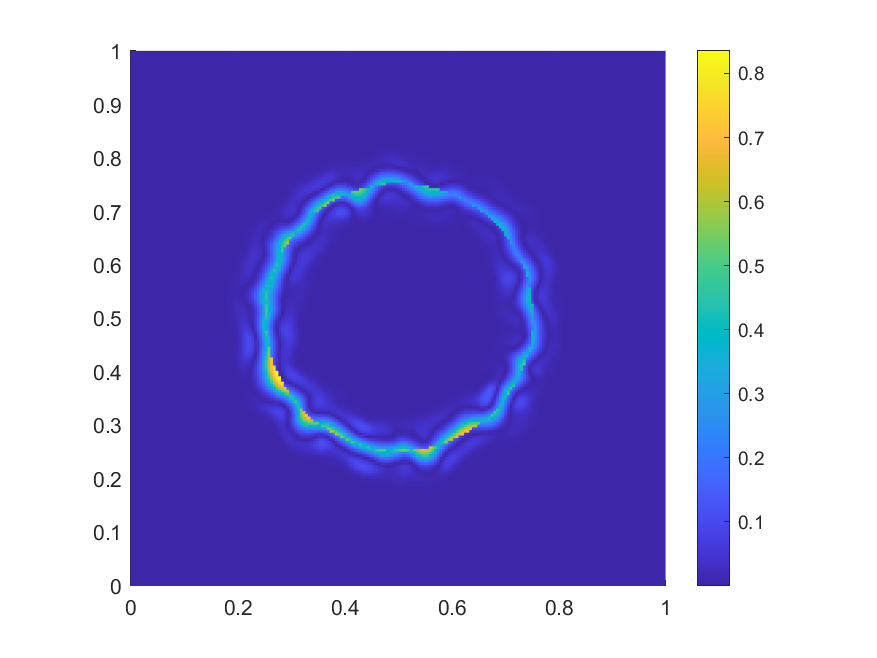} & \hspace{-1cm}\includegraphics[width=5.6cm]{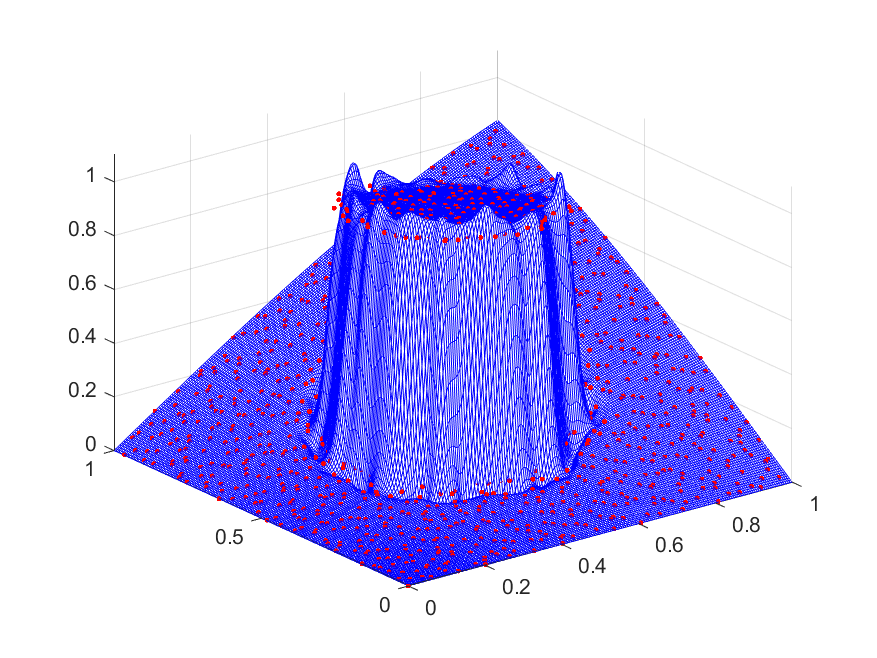} & \hspace{-1cm}			\includegraphics[width=5.6cm]{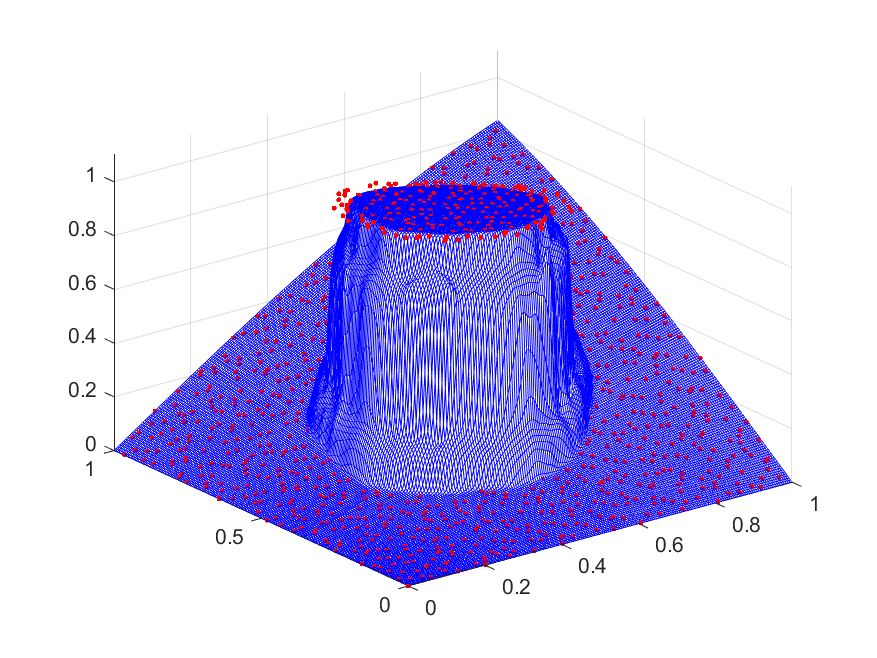} & \hspace{-1cm} \includegraphics[width=5.6cm]{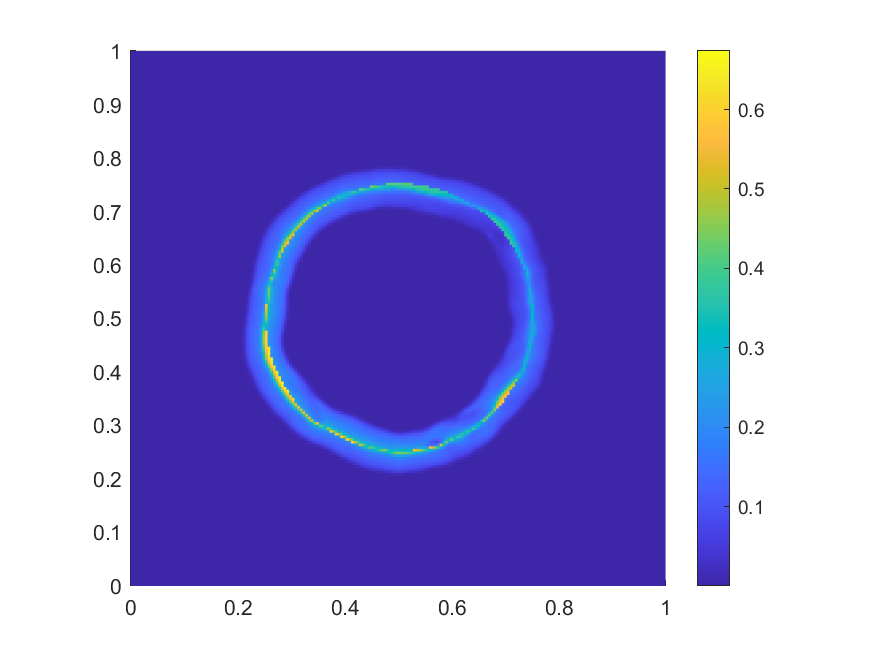}\\
\hspace{-1cm}	\includegraphics[width=5.6cm]{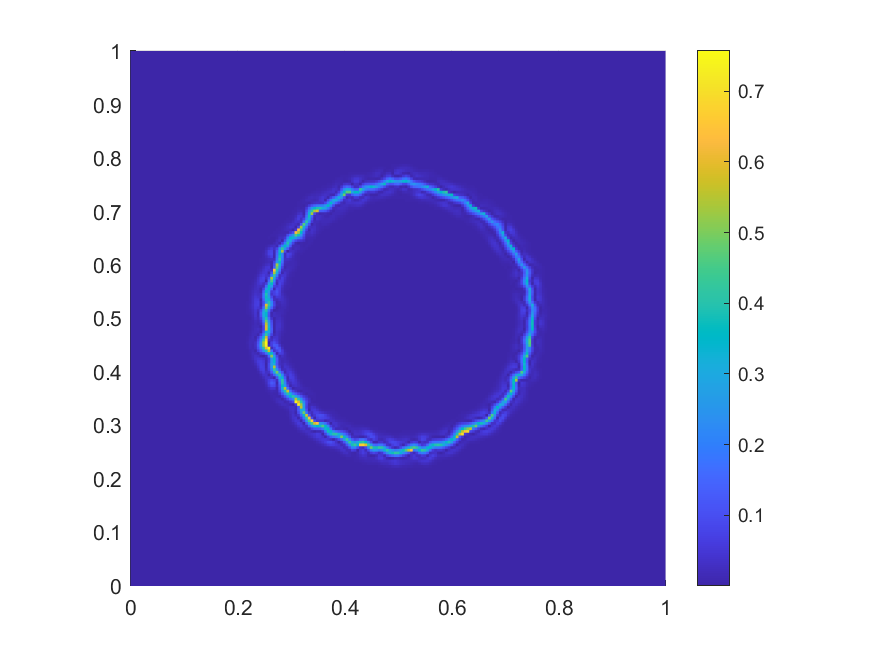} &  \hspace{-1cm}	\includegraphics[width=5.6cm]{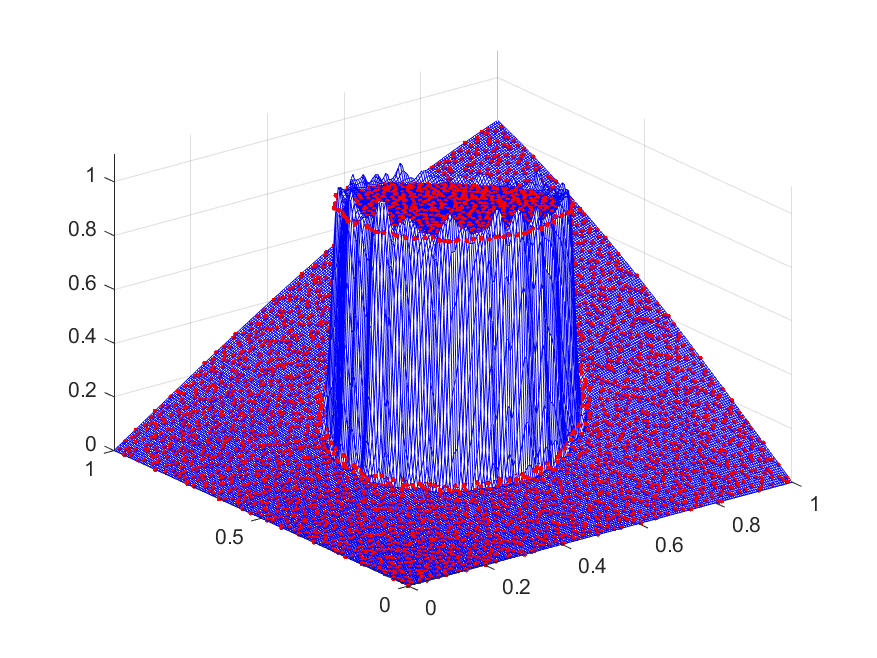} & \hspace{-1cm}
			\includegraphics[width=5.6cm]{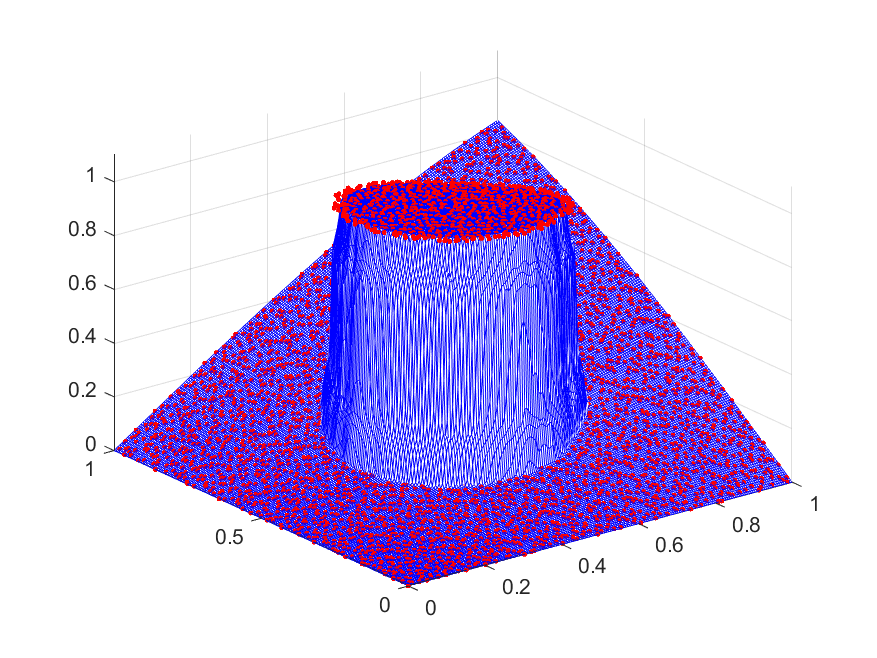} & \hspace{-1cm} \includegraphics[width=5.6cm]{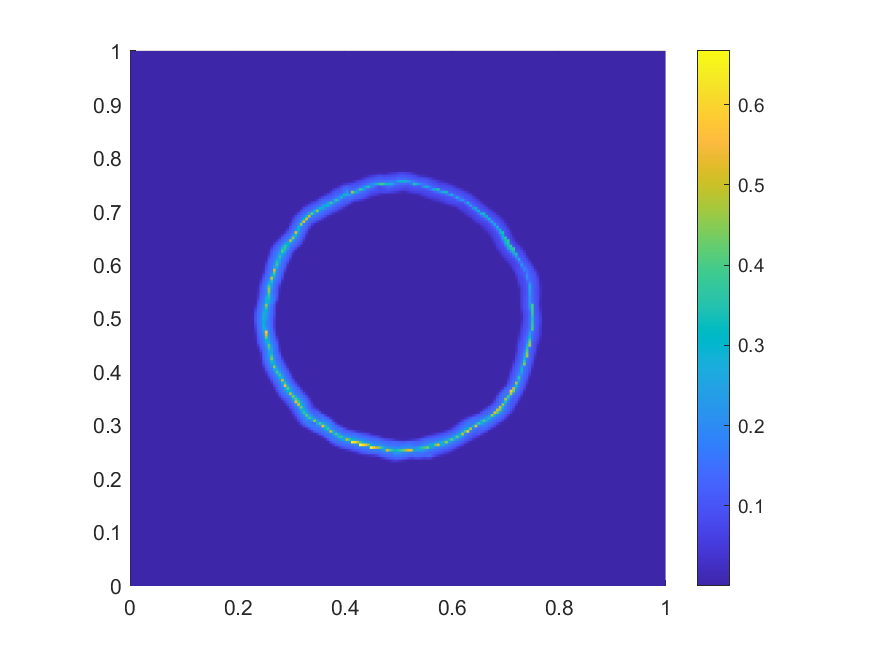}
		\end{tabular}
\end{center}
\caption{Approximation to the function $z$, Eq. \eqref{ejemplonuma}, using linear and data dependent MLS with $\omega(x)$ the $\mathcal{C}^\infty$ Gauss function \cite{wendland2002} and the class of polynomials $\Pi_2(\mathbb{R}^2)$. Red points: original function, blue lines: approximation. }
	\label{figureejemplonuma2}
	\end{figure}

\end{landscape}

\subsection{Reduction of the smearing zone across the discontinuity}

\begin{figure}[!ht]
\begin{center}
		\begin{tabular}{cccc}
\multicolumn{3}{c}{MLS$^0_{\text{W2}}$} \\
 \hspace{-1cm} 2D-plot& \hspace{-1cm} Approximation &  \hspace{-1cm} Error\\
\hspace{-1cm}			\includegraphics[width=6.3cm]{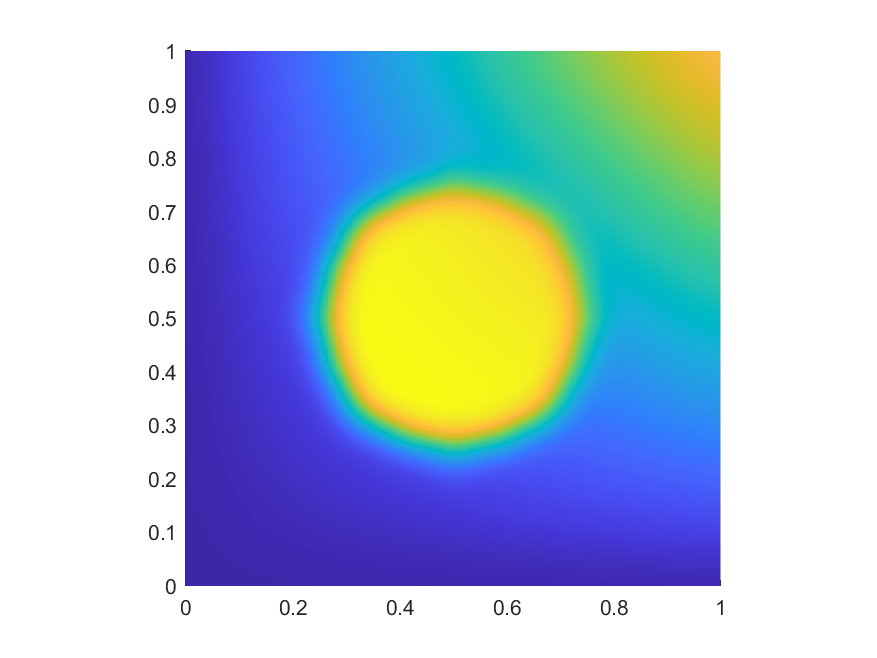}  &
\hspace{-1cm}		    \includegraphics[width=6.3cm]{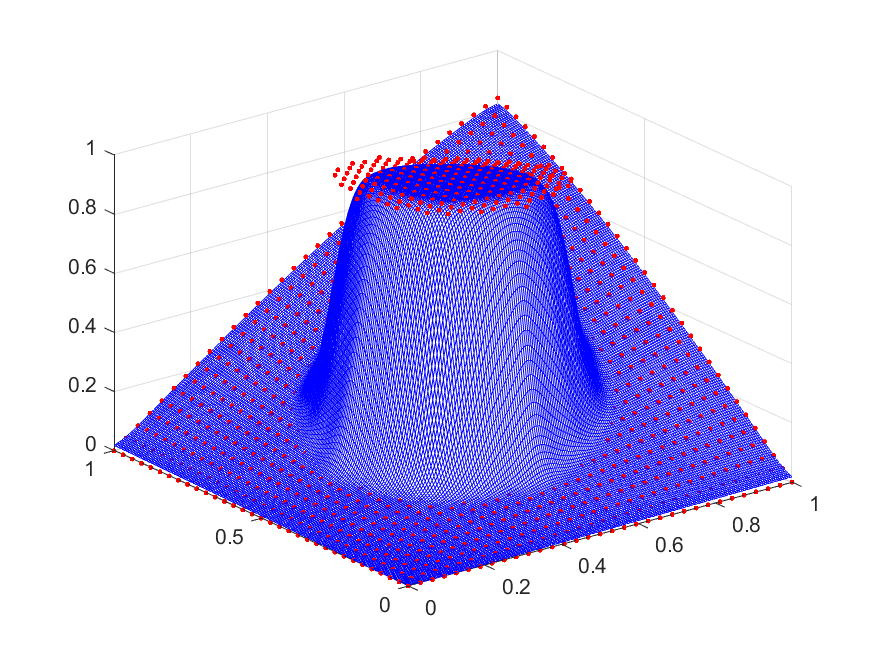} &	
\hspace{-1cm}            \includegraphics[width=6.3cm]{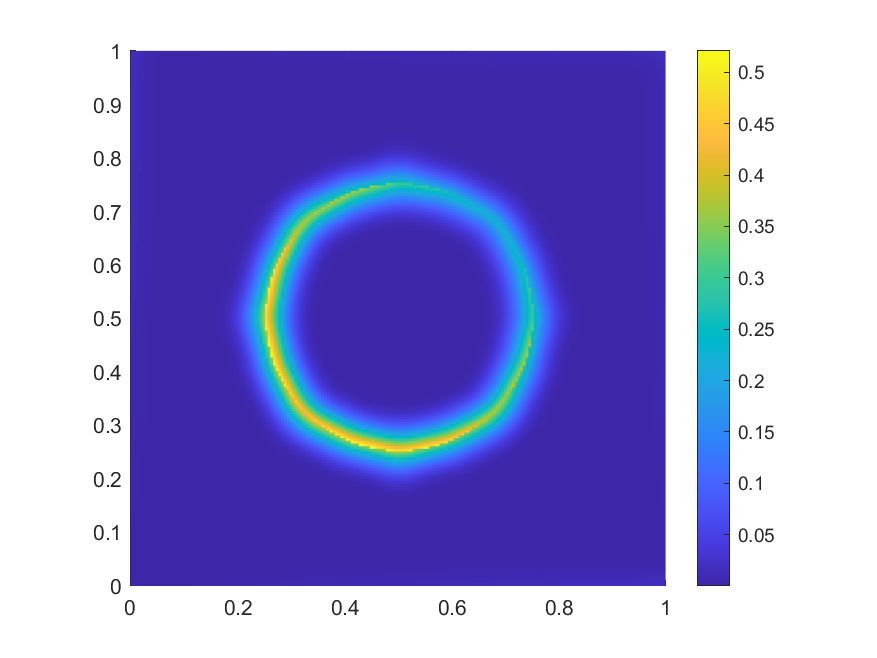} & 	
			\\
\hspace{-1cm}			\includegraphics[width=6.3cm]{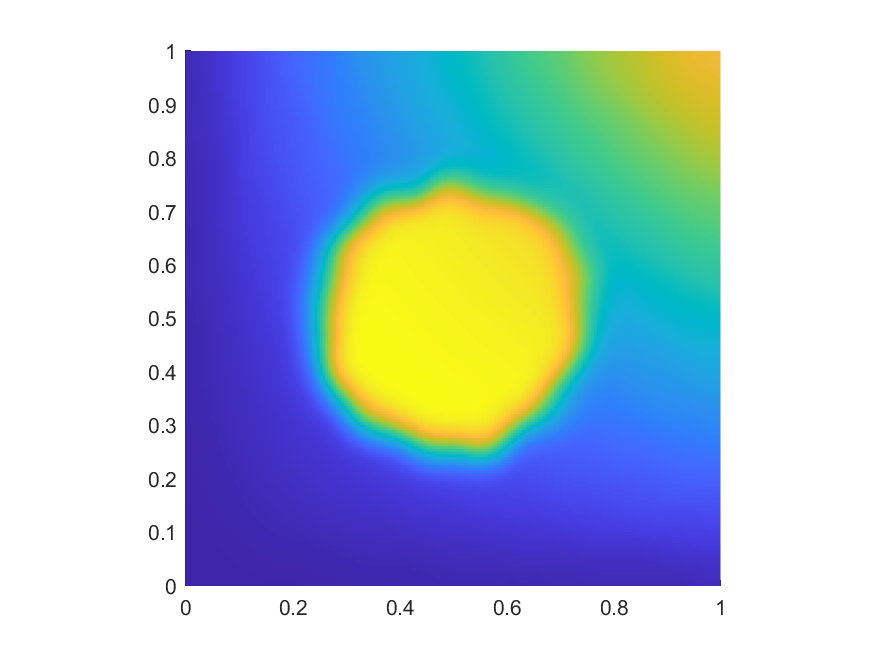}  &
\hspace{-1cm}		    \includegraphics[width=6.3cm]{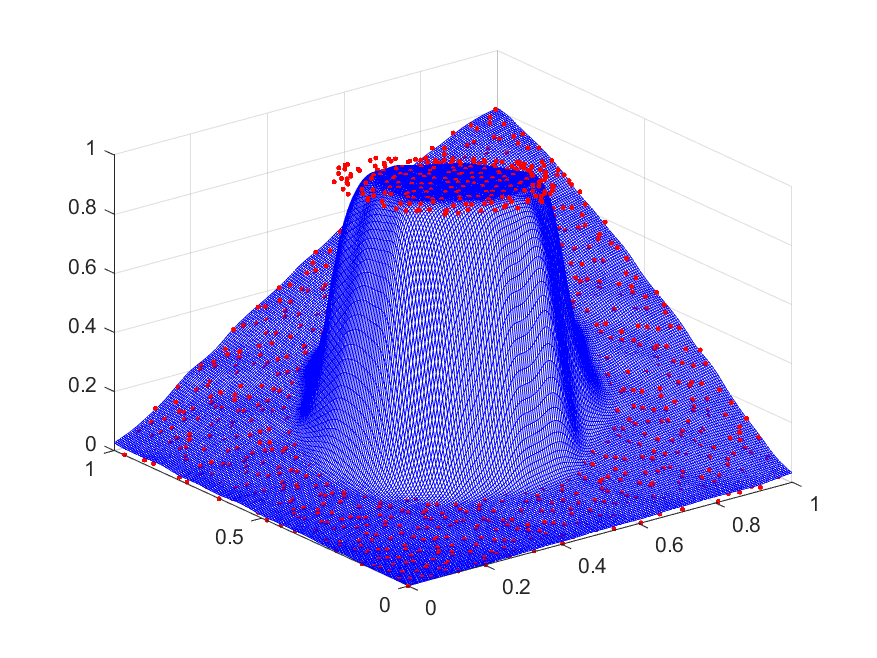} &	
\hspace{-1cm}            \includegraphics[width=6.3cm]{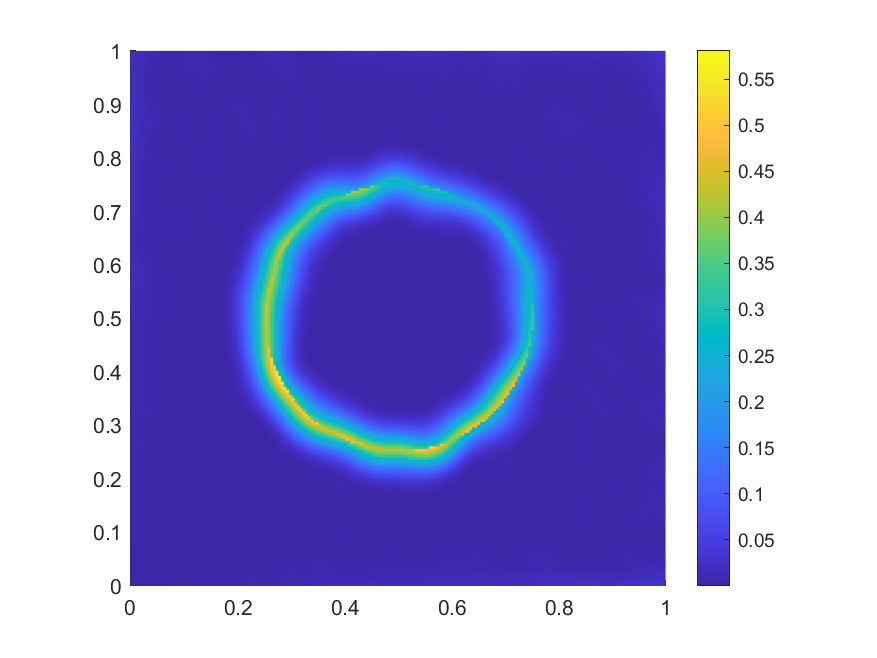} & 	
			\\
 \multicolumn{3}{c}{DD-MLS$^0_{\text{W2}}$}\\
\hspace{-1cm}			\includegraphics[width=6.3cm]{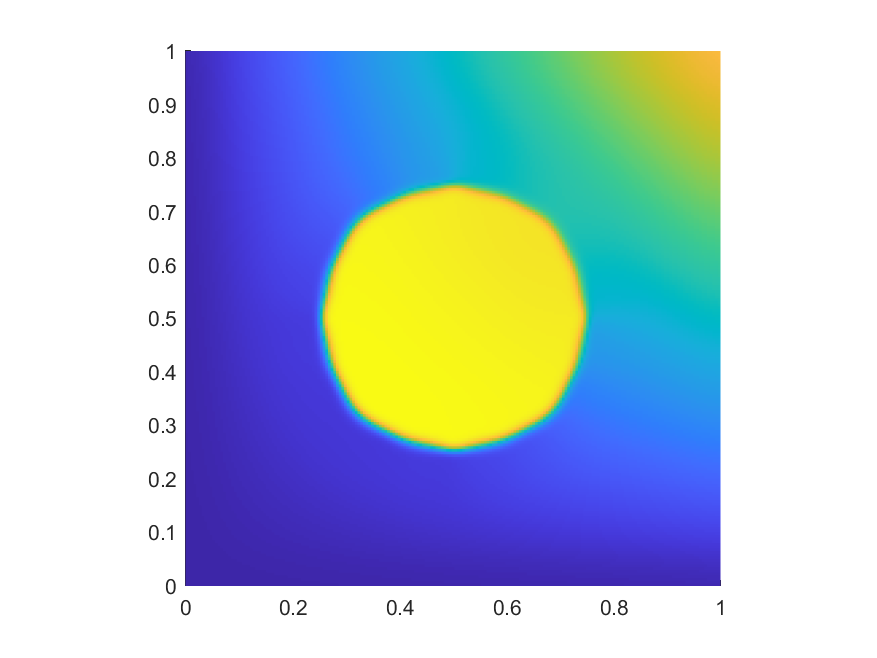}  &
\hspace{-1cm}			\includegraphics[width=6.3cm]{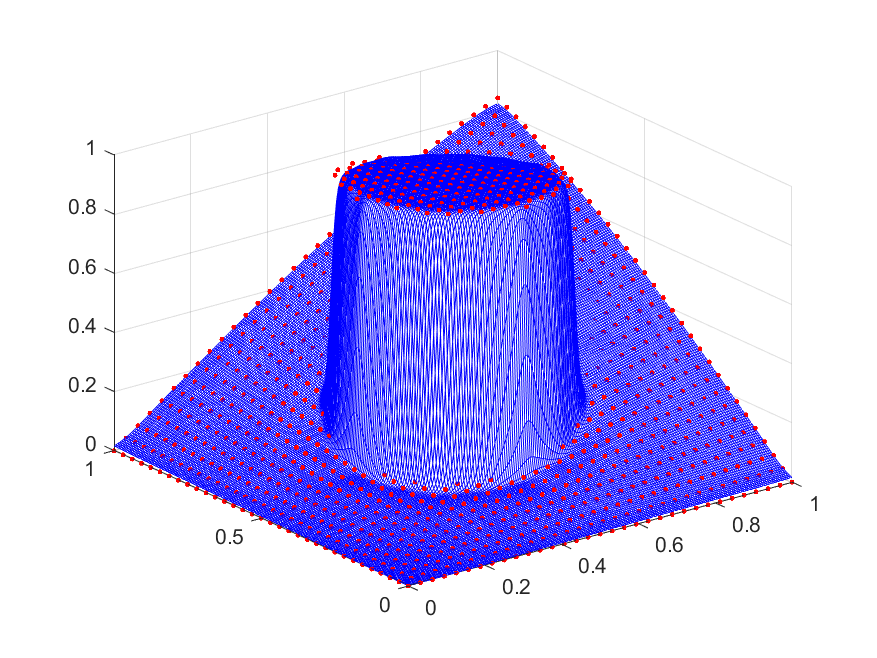} &
\hspace{-1cm}            \includegraphics[width=6.3cm]{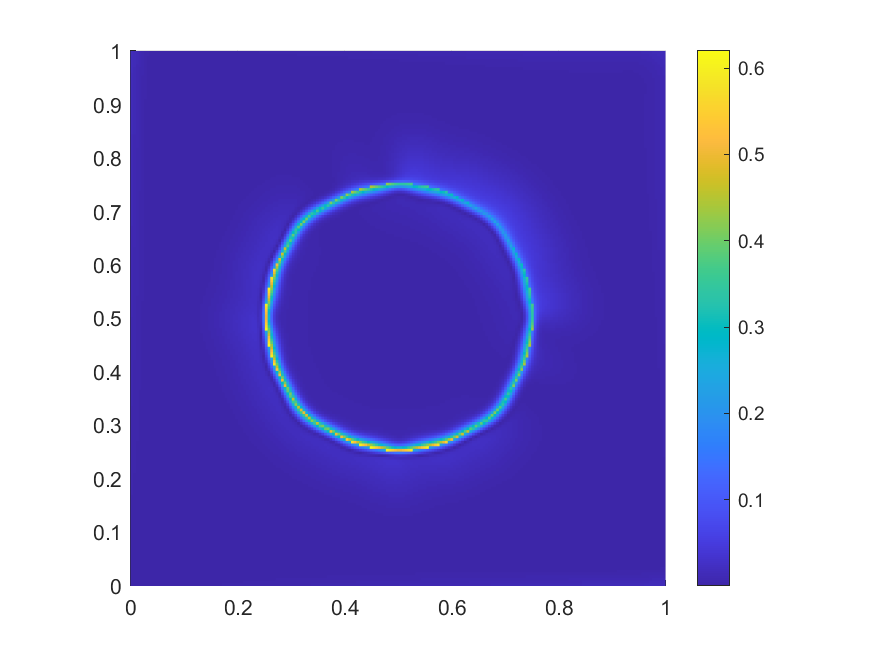} & 	
				\\
\hspace{-1cm}			\includegraphics[width=6.3cm]{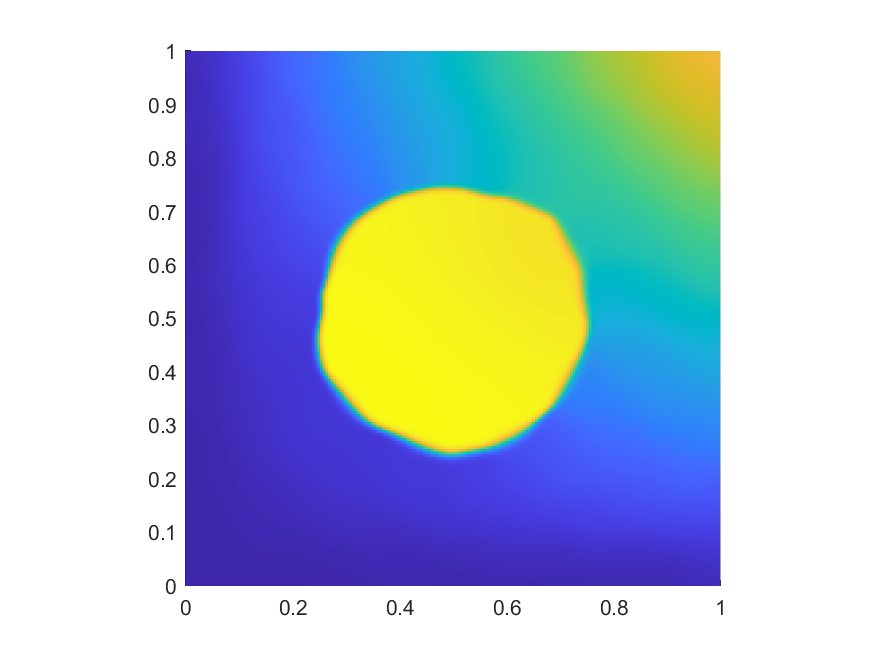}  &
\hspace{-1cm}			\includegraphics[width=6.3cm]{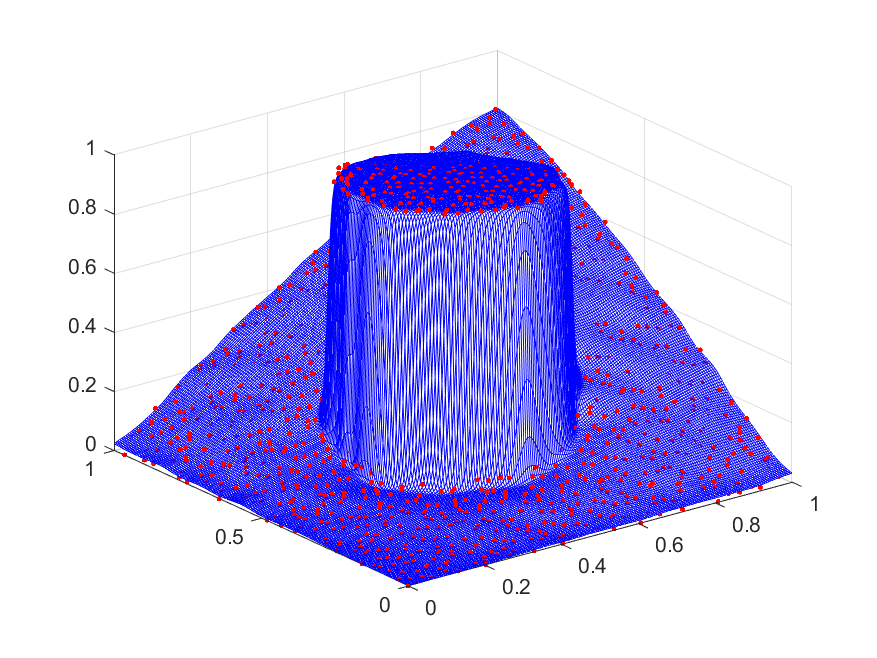} &
\hspace{-1cm}            \includegraphics[width=6.3cm]{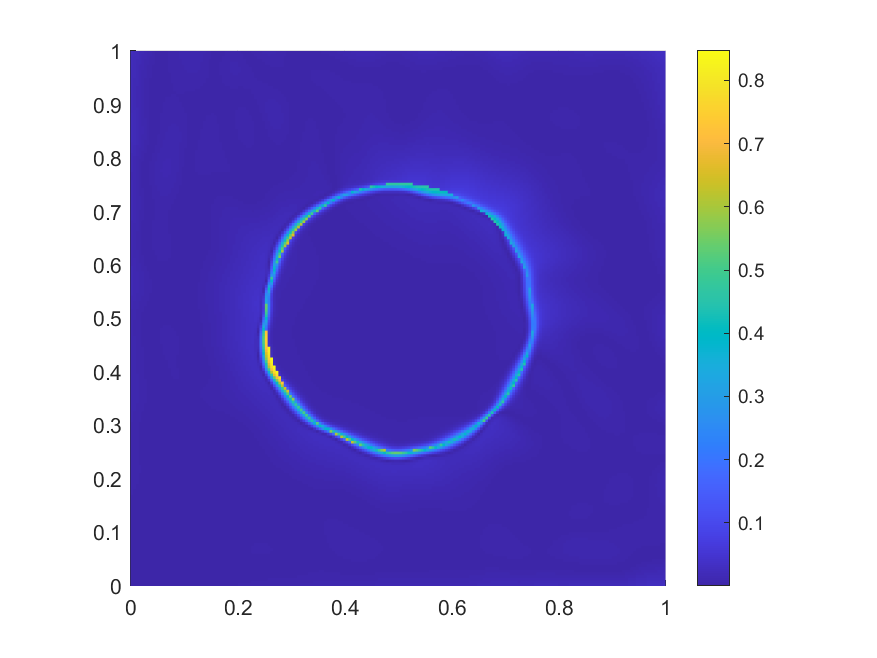} & 	

		\end{tabular}
\end{center}
				\caption{Approximation to function $z$, Eq. \eqref{ejemplonuma}, using W2, $d=0$ and $N=33^2$.}
		\label{figuraW2d0N33}
	\end{figure}

\begin{figure}[!ht]
\begin{center}
		\begin{tabular}{cccc}
\multicolumn{3}{c}{MLS$^0_{\text{W2}}$} \\
 \hspace{-1cm} 2D-plot& \hspace{-1cm} Approximation &  \hspace{-1cm} Error\\
\hspace{-1cm}			\includegraphics[width=6.3cm]{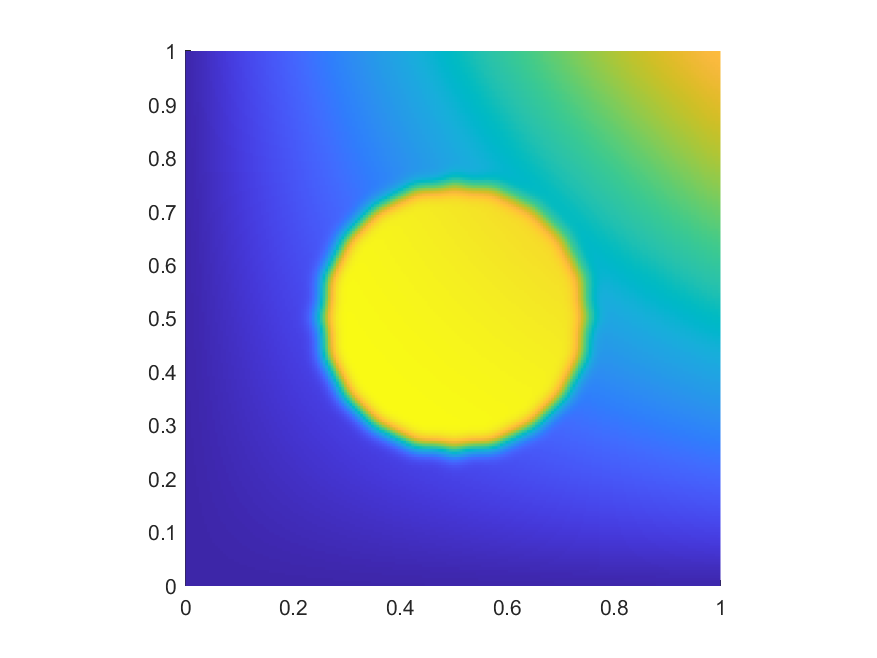}  &
\hspace{-1cm}		    \includegraphics[width=6.3cm]{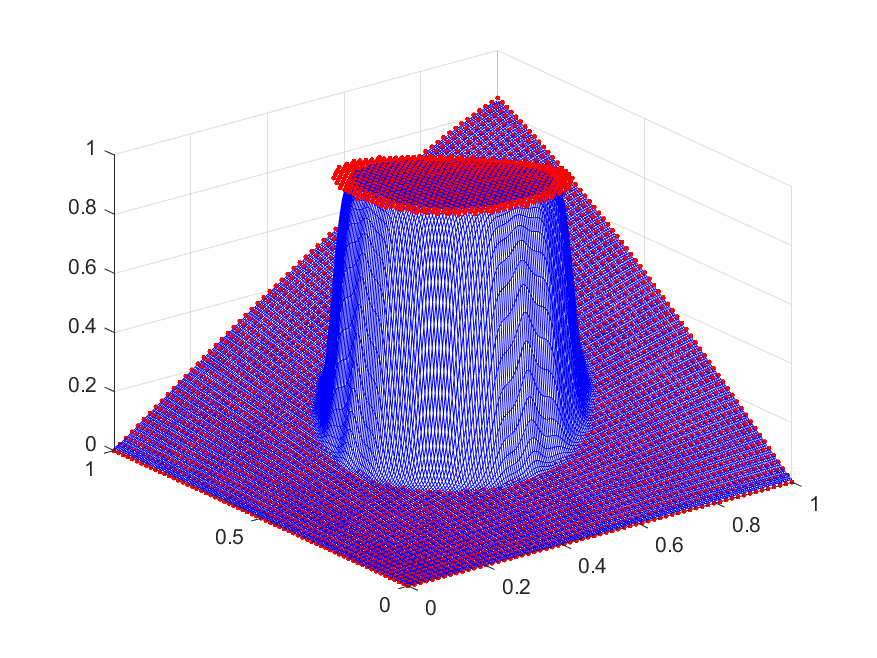} &	
\hspace{-1cm}            \includegraphics[width=6.3cm]{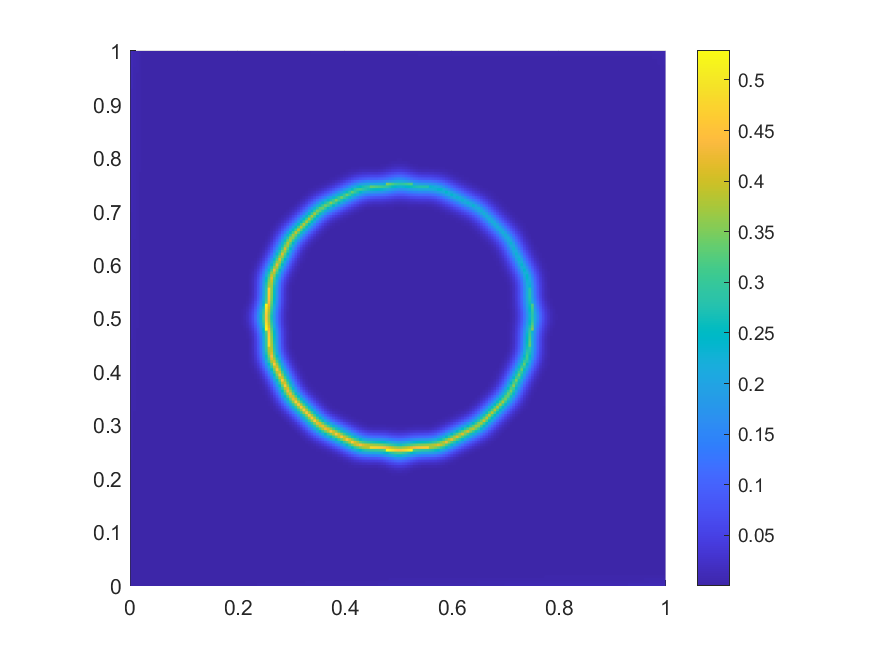} & 	
			\\
\hspace{-1cm}			\includegraphics[width=6.3cm]{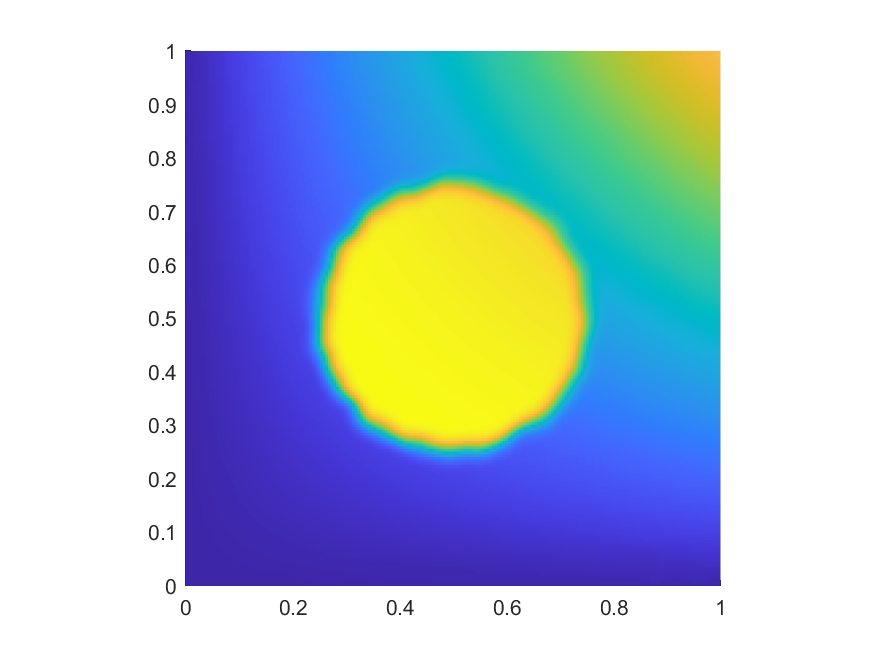}  &
\hspace{-1cm}		    \includegraphics[width=6.3cm]{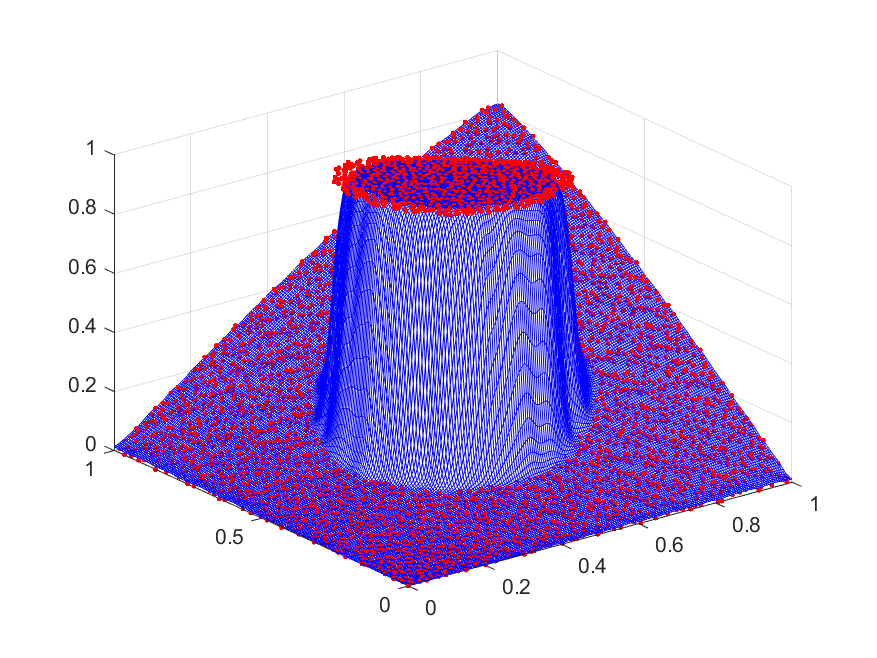} &	
\hspace{-1cm}            \includegraphics[width=6.3cm]{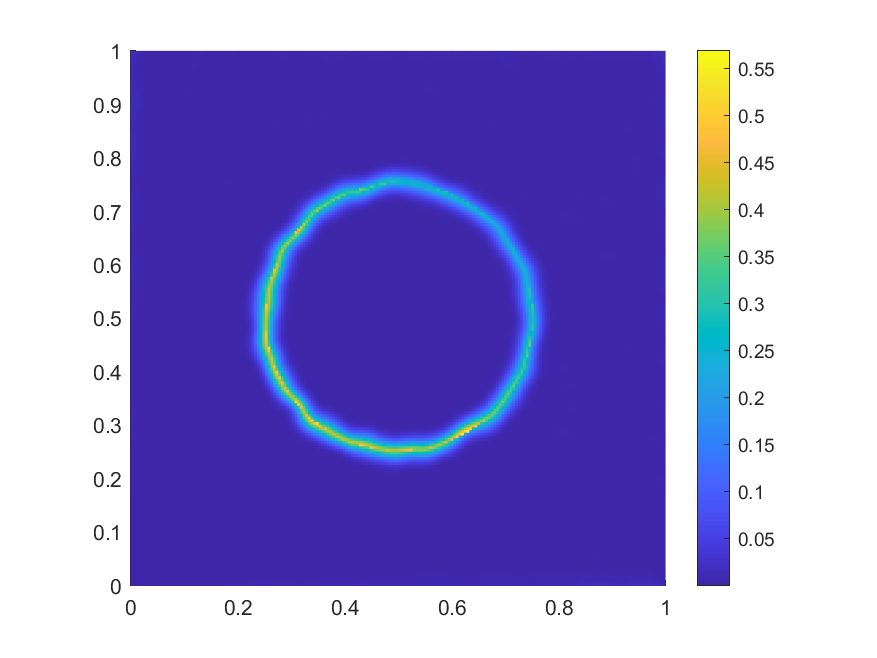} & 	
			\\
 \multicolumn{3}{c}{DD-MLS$^0_{\text{W2}}$}\\
\hspace{-1cm}			\includegraphics[width=6.3cm]{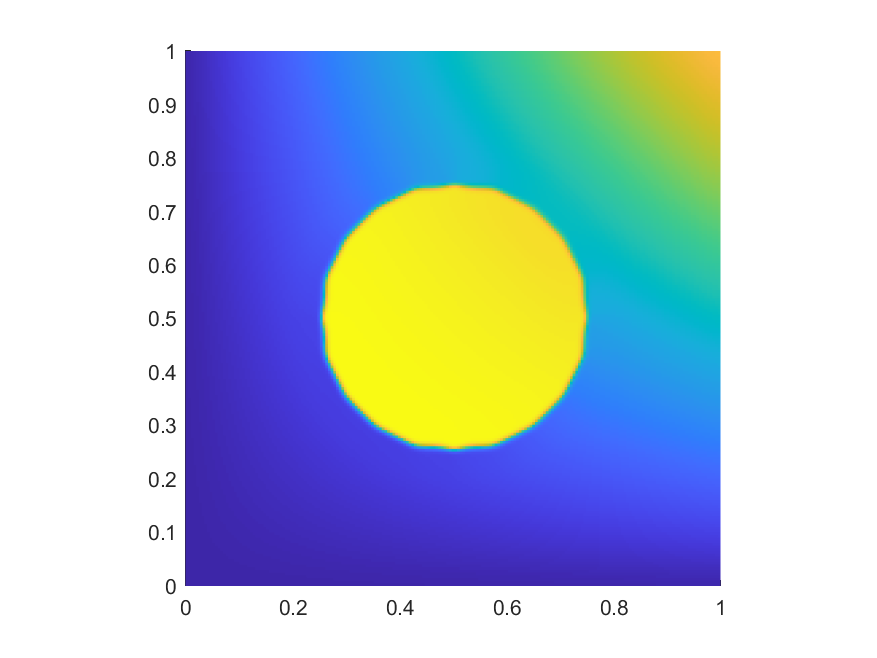}  &
\hspace{-1cm}			\includegraphics[width=6.3cm]{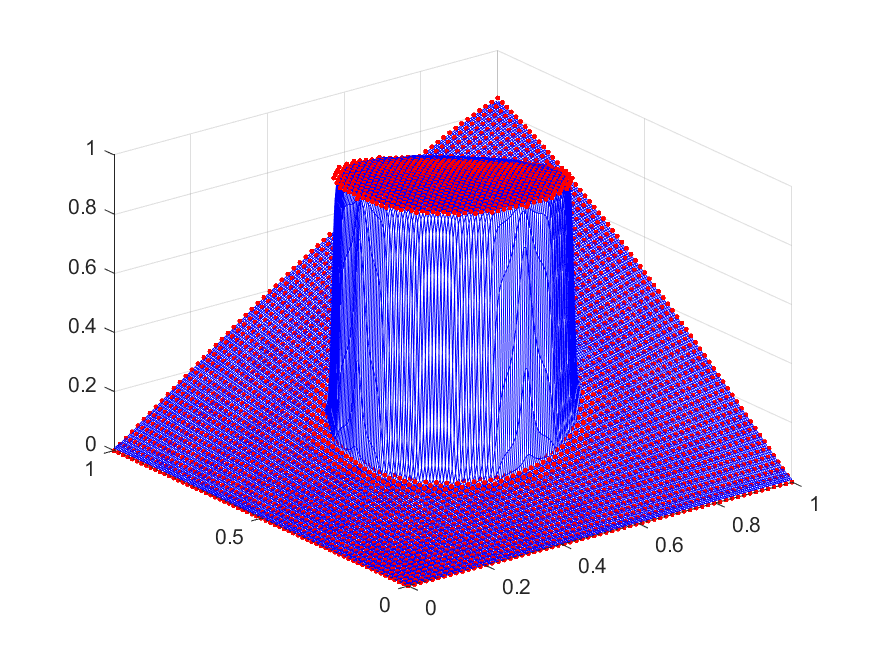} &
\hspace{-1cm}            \includegraphics[width=6.3cm]{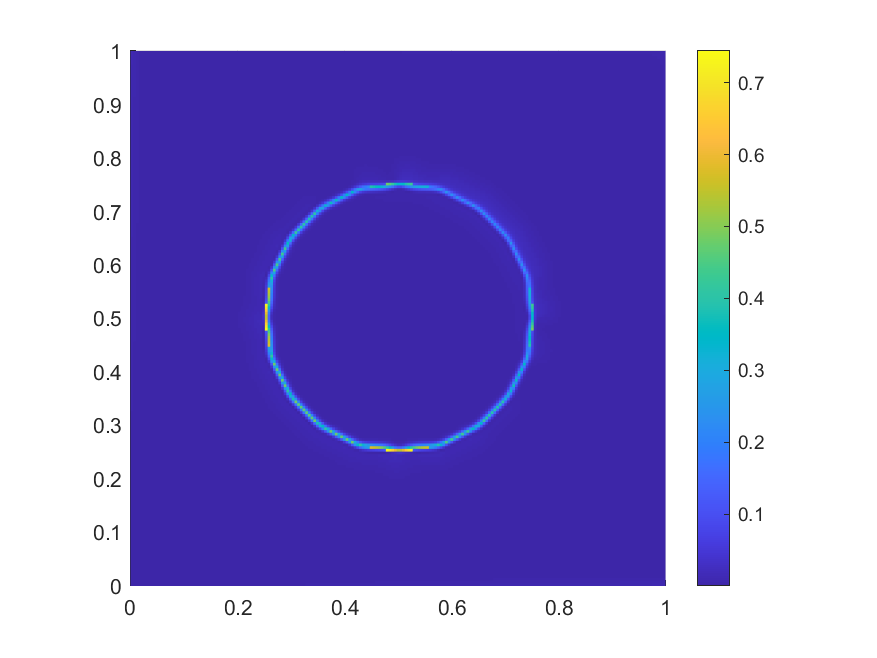} & 	
				\\
\hspace{-1cm}			\includegraphics[width=6.3cm]{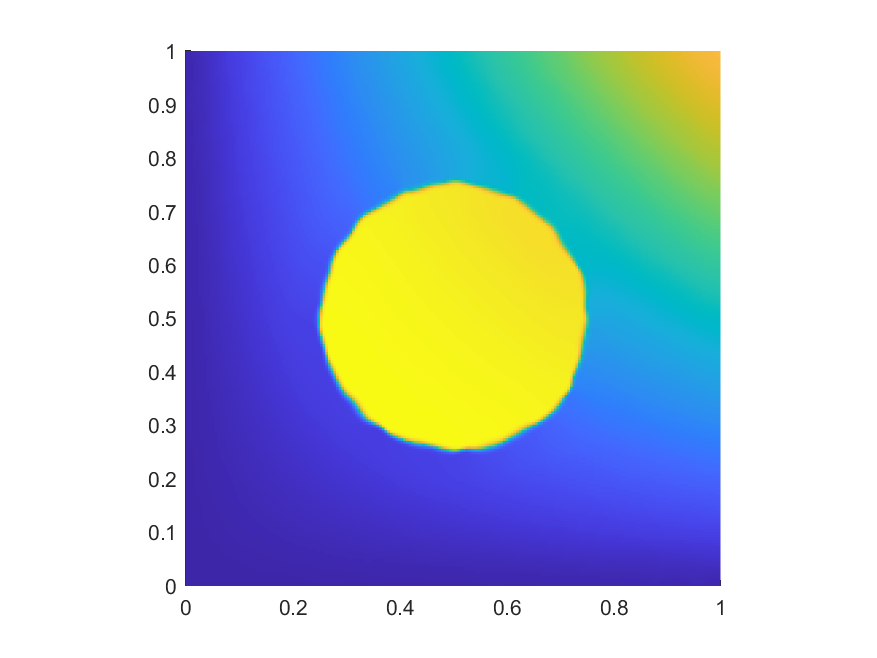}  &
\hspace{-1cm}			\includegraphics[width=6.3cm]{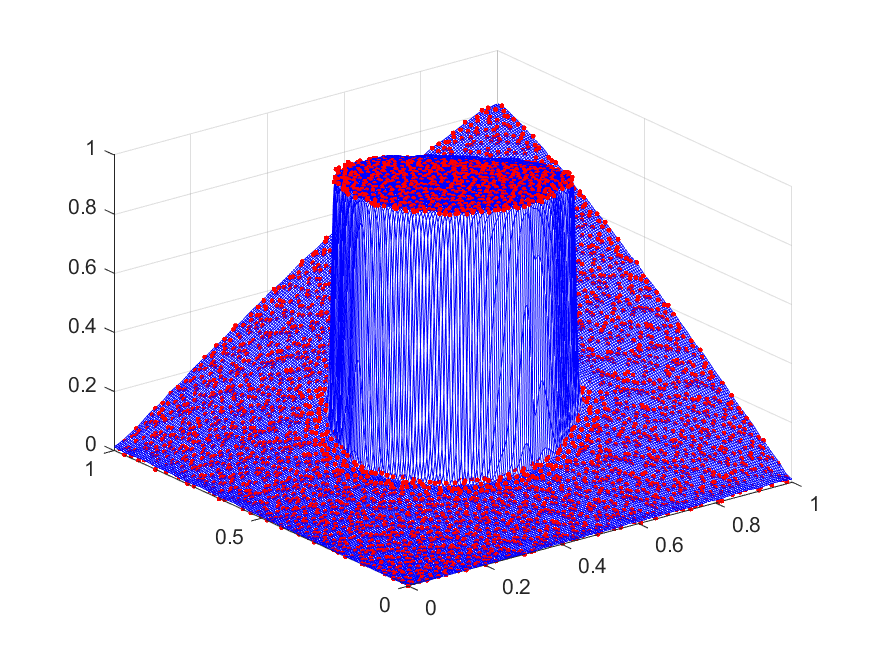} &
\hspace{-1cm}            \includegraphics[width=6.3cm]{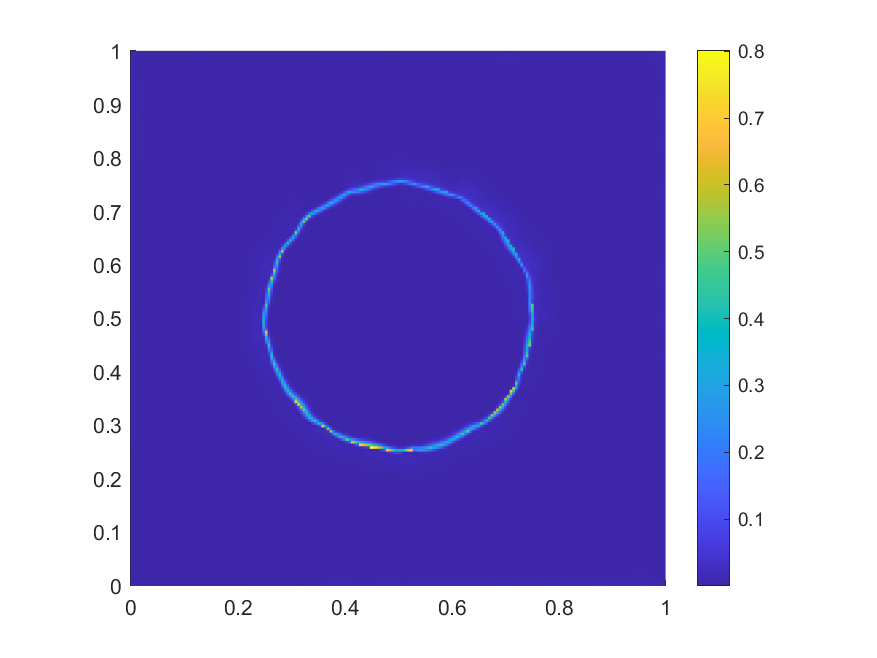} & 	

		\end{tabular}
\end{center}
				\caption{Approximation to function $z$, Eq. \eqref{ejemplonuma}, using W2, $d=0$ with $N=65^2$.}
		\label{figuraW2d0N65}
	\end{figure}

\begin{figure}[!ht]
\begin{center}
		\begin{tabular}{cccc}
\multicolumn{3}{c}{MLS$^0_{\text{G}}$} \\
 \hspace{-1cm} 2D-plot& \hspace{-1cm} Approximation &  \hspace{-1cm} Error\\
\hspace{-1cm}			\includegraphics[width=6.3cm]{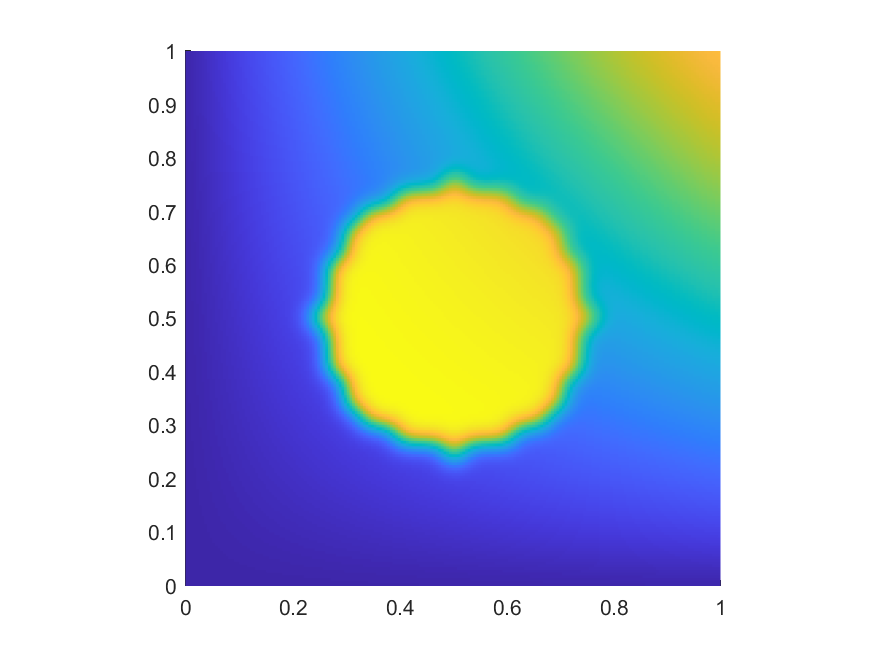}  &
\hspace{-1cm}		    \includegraphics[width=6.3cm]{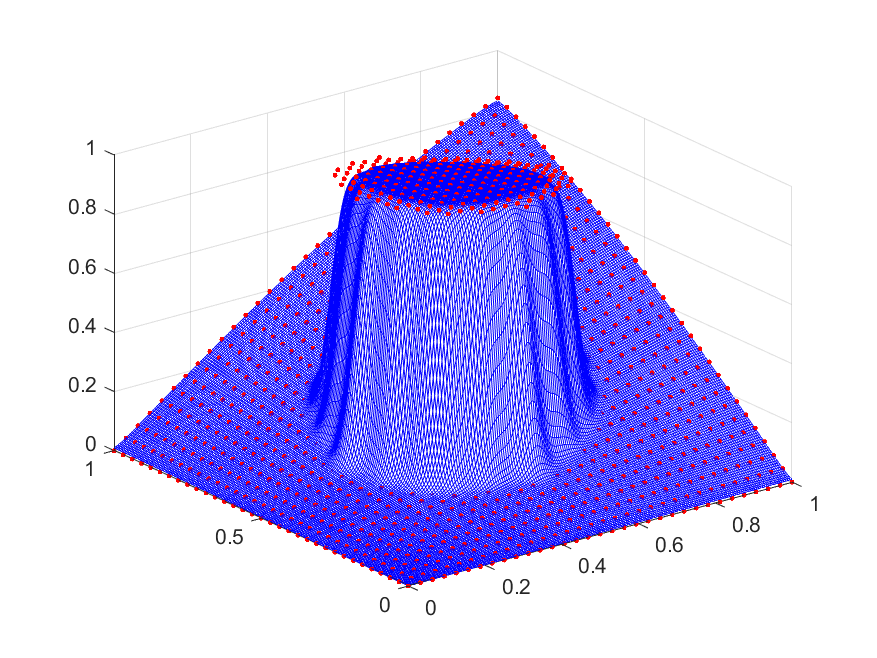} &	
\hspace{-1cm}            \includegraphics[width=6.3cm]{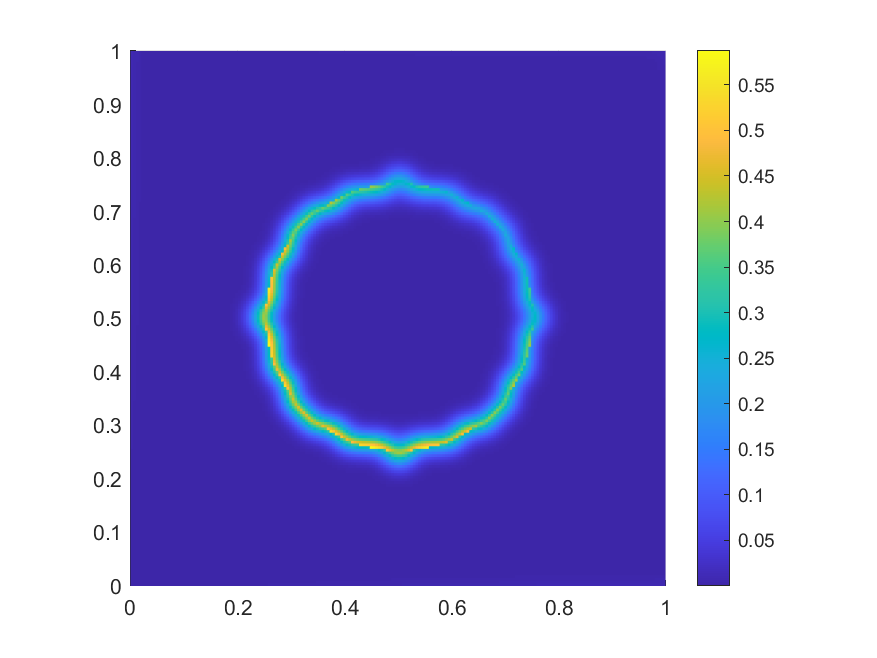} & 	
			\\
\hspace{-1cm}			\includegraphics[width=6.3cm]{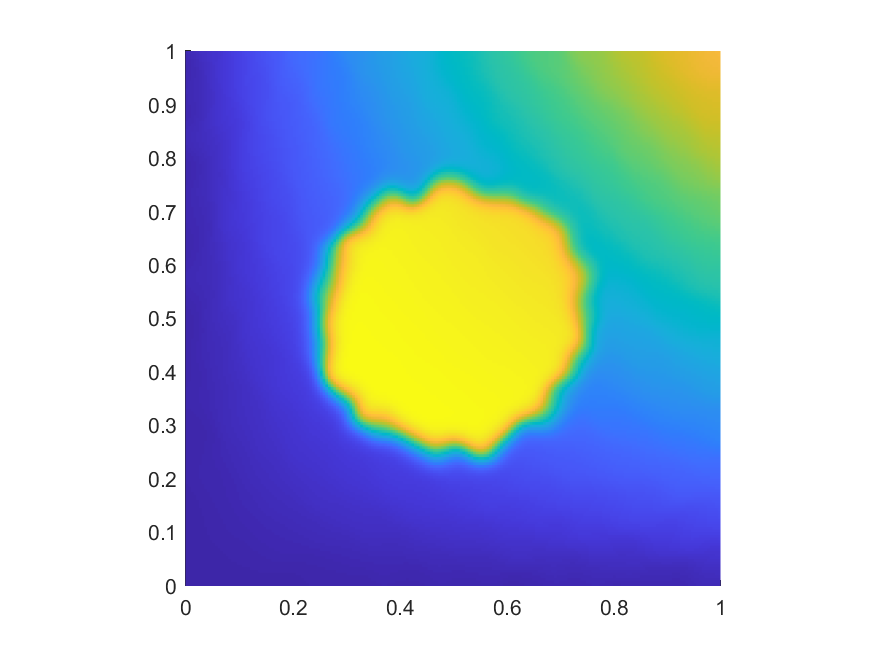}  &
\hspace{-1cm}		    \includegraphics[width=6.3cm]{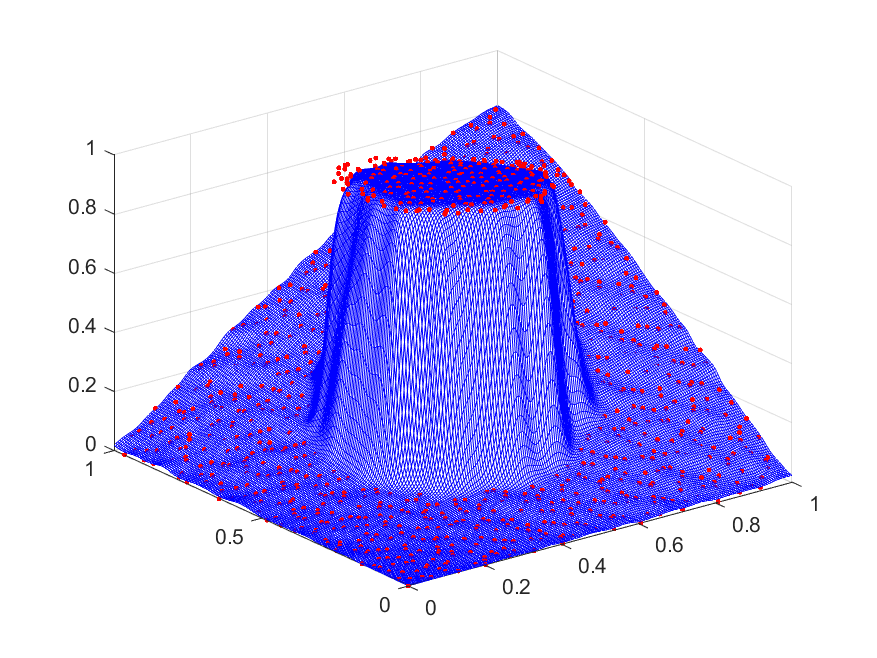} &	
\hspace{-1cm}            \includegraphics[width=6.3cm]{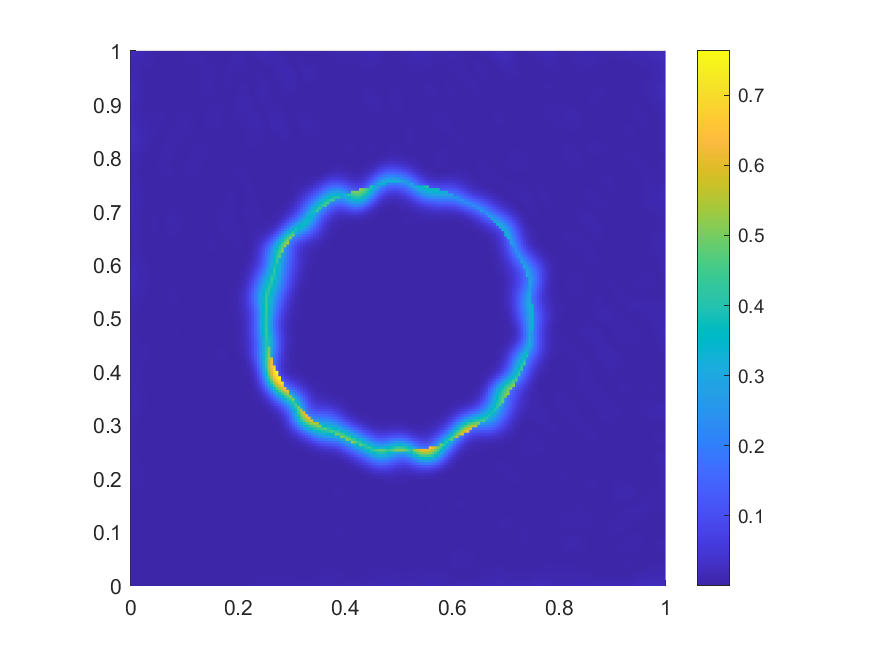} & 	
			\\
 \multicolumn{3}{c}{DD-MLS$^0_{\text{G}}$}\\
\hspace{-1cm}			\includegraphics[width=6.3cm]{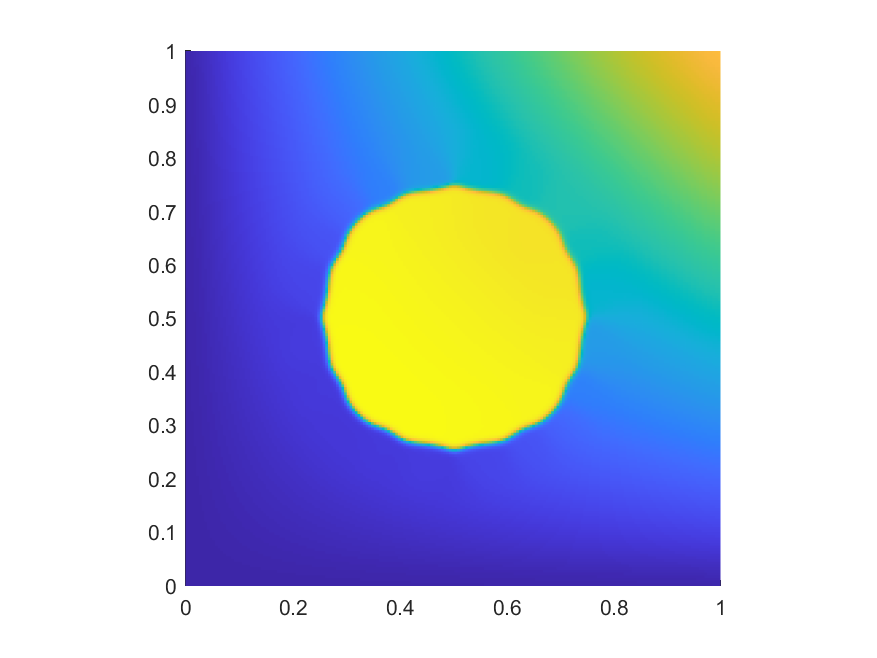}  &
\hspace{-1cm}			\includegraphics[width=6.3cm]{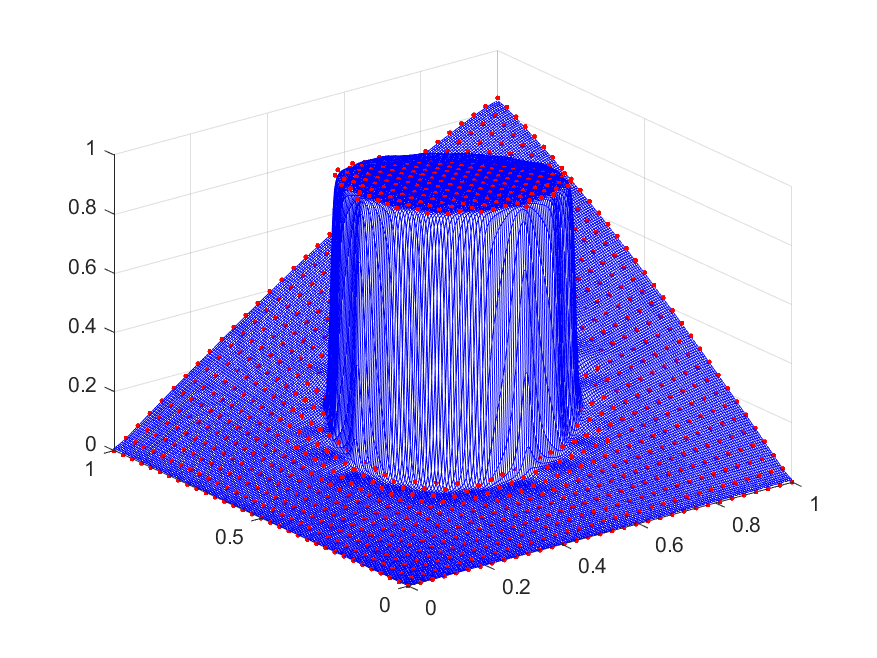} &
\hspace{-1cm}            \includegraphics[width=6.3cm]{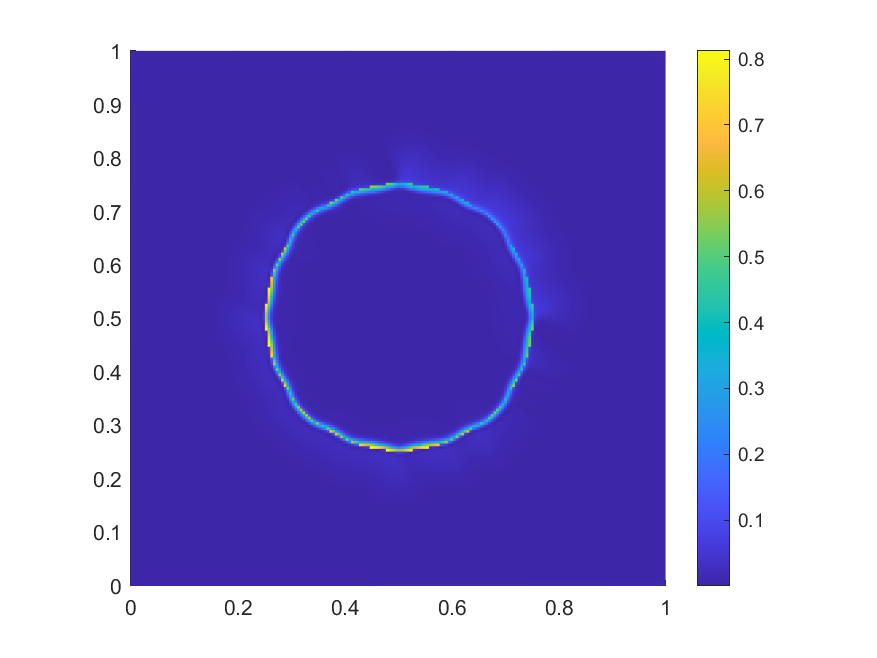} & 	
				\\
\hspace{-1cm}			\includegraphics[width=6.3cm]{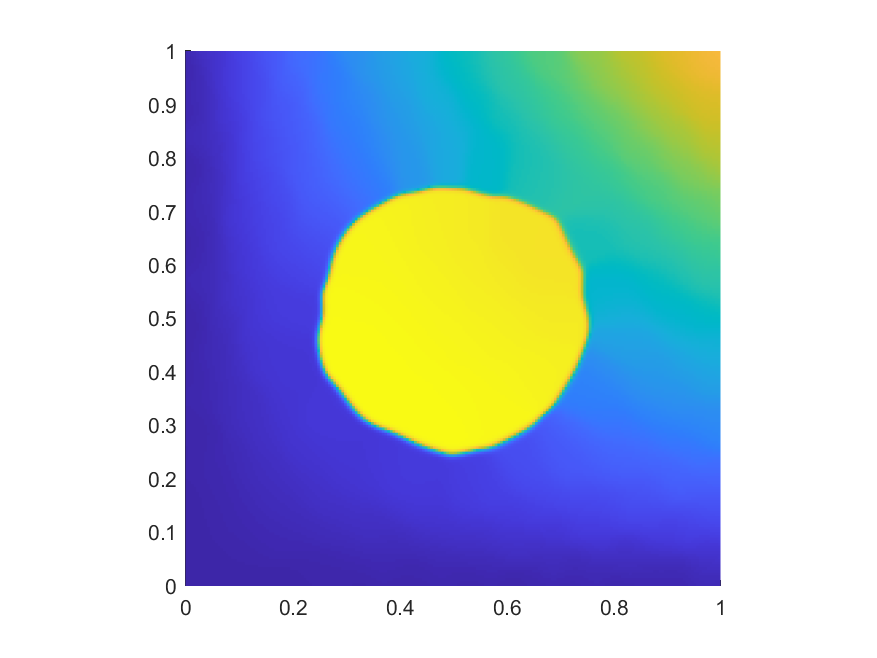}  &
\hspace{-1cm}			\includegraphics[width=6.3cm]{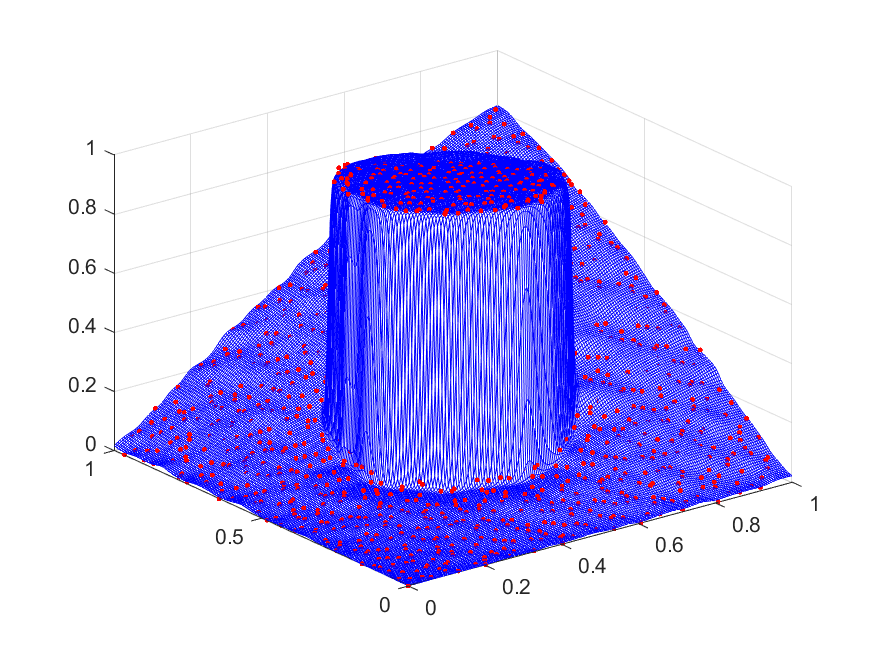} &
\hspace{-1cm}            \includegraphics[width=6.3cm]{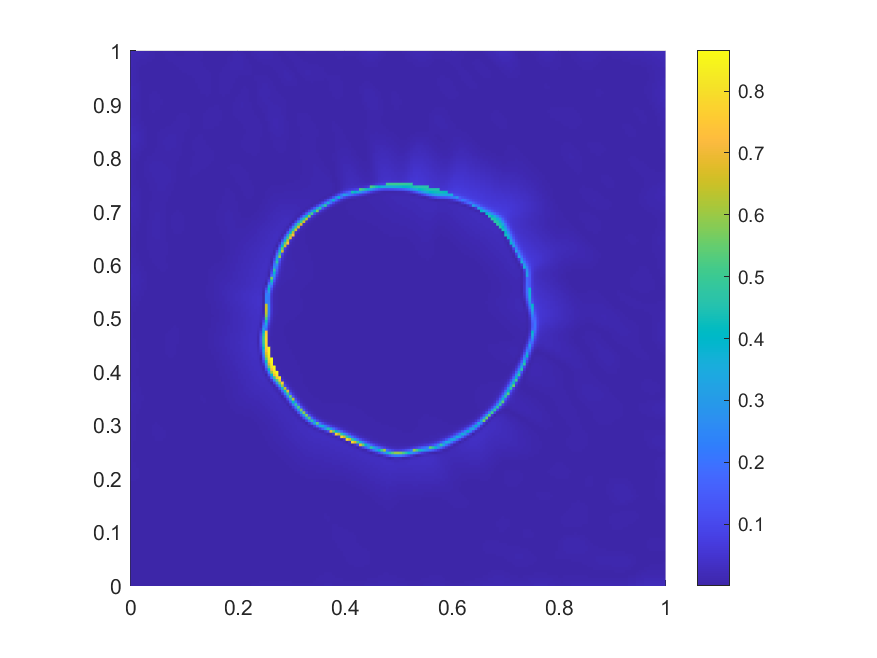} & 	

		\end{tabular}
\end{center}
				\caption{Approximation to function $z$, Eq. \eqref{ejemplonuma}, using G, $d=0$ with $N=33^2$.}
		\label{figuraGd0N33}
	\end{figure}

\begin{figure}[!ht]
\begin{center}
		\begin{tabular}{ccc}
\multicolumn{3}{c}{MLS$^0_{\text{G}}$} \\
 \hspace{-1cm} 2D-plot& \hspace{-1cm} Approximation &  \hspace{-1cm} Error\\
\hspace{-1cm}			\includegraphics[width=6.3cm]{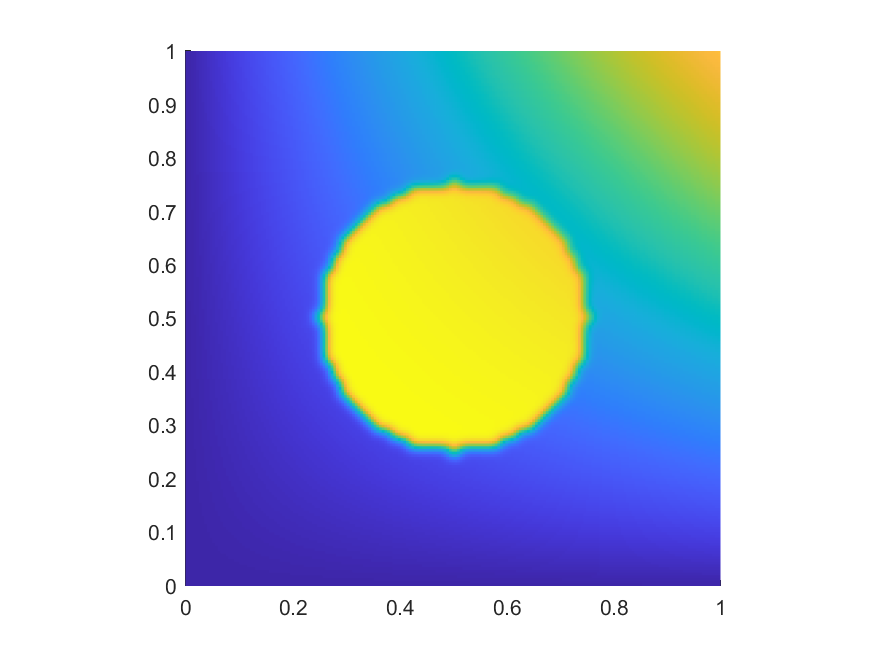}  &
\hspace{-1cm}		    \includegraphics[width=6.3cm]{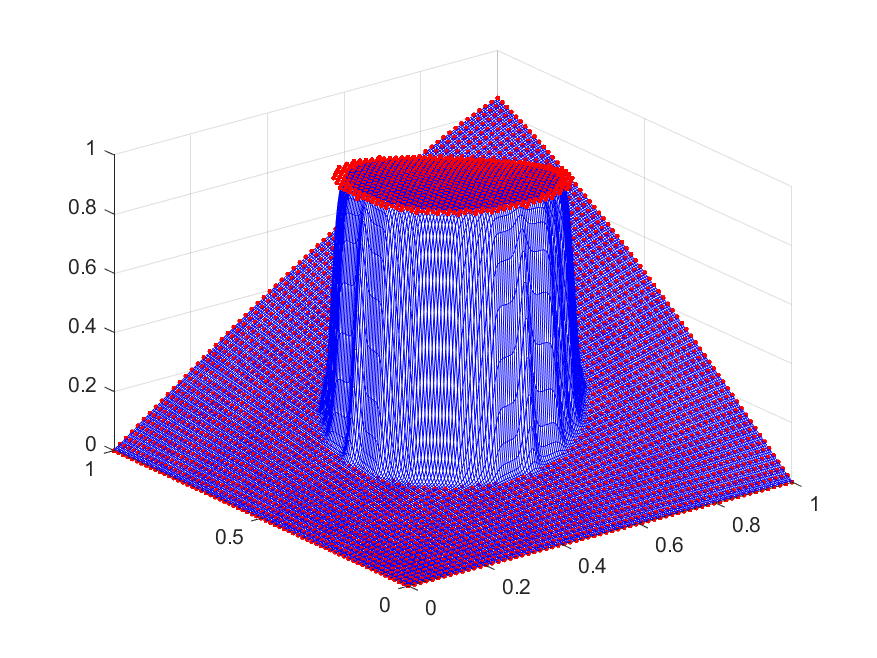} &	
\hspace{-1cm}            \includegraphics[width=6.3cm]{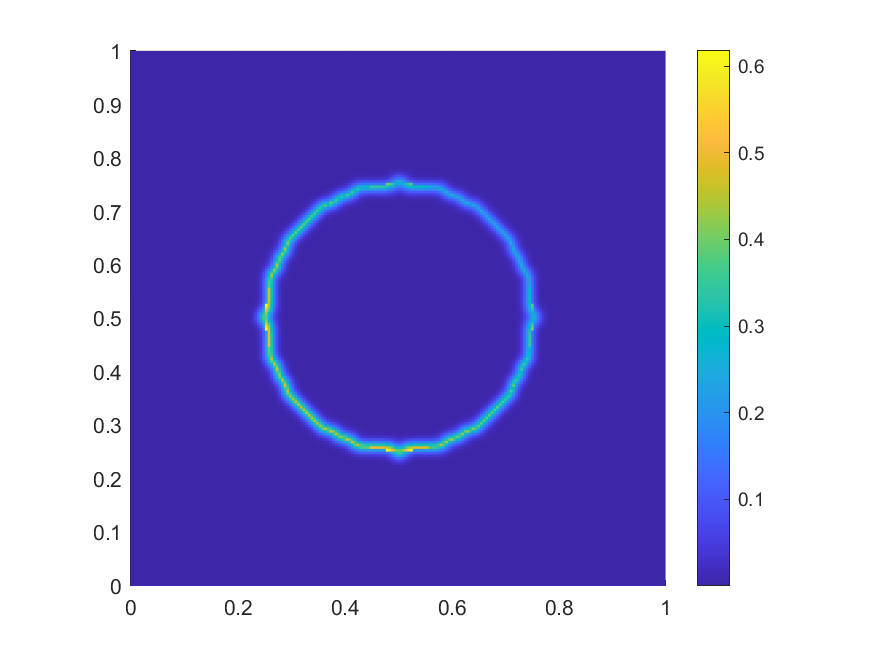}  	
			\\
\hspace{-1cm}			\includegraphics[width=6.3cm]{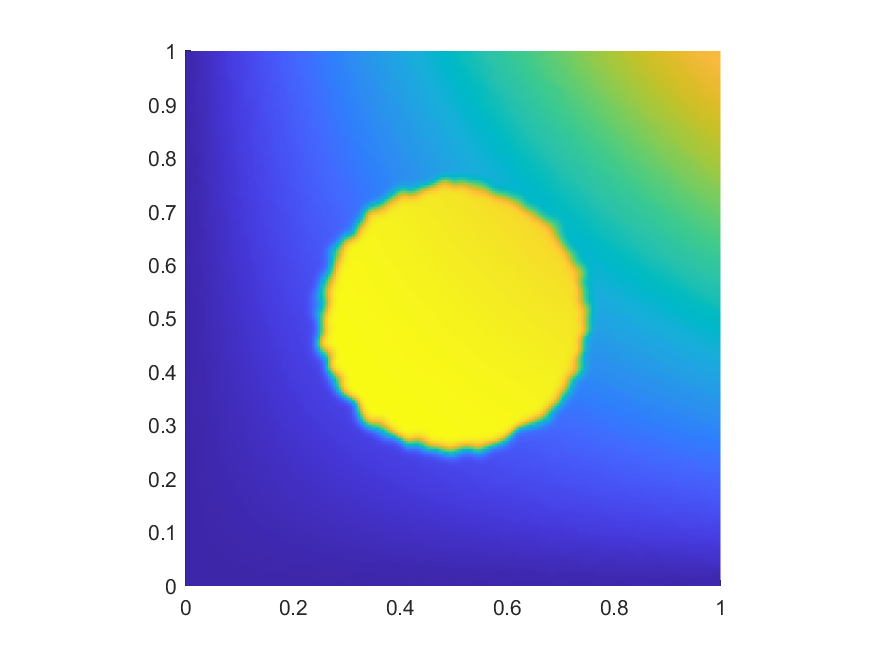}  &
\hspace{-1cm}		    \includegraphics[width=6.3cm]{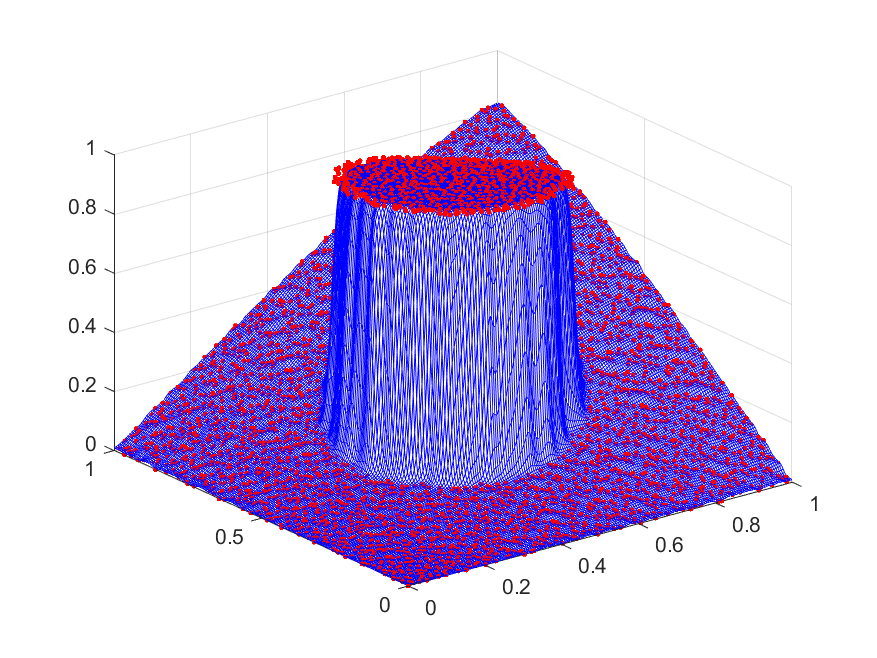} &	
\hspace{-1cm}            \includegraphics[width=6.3cm]{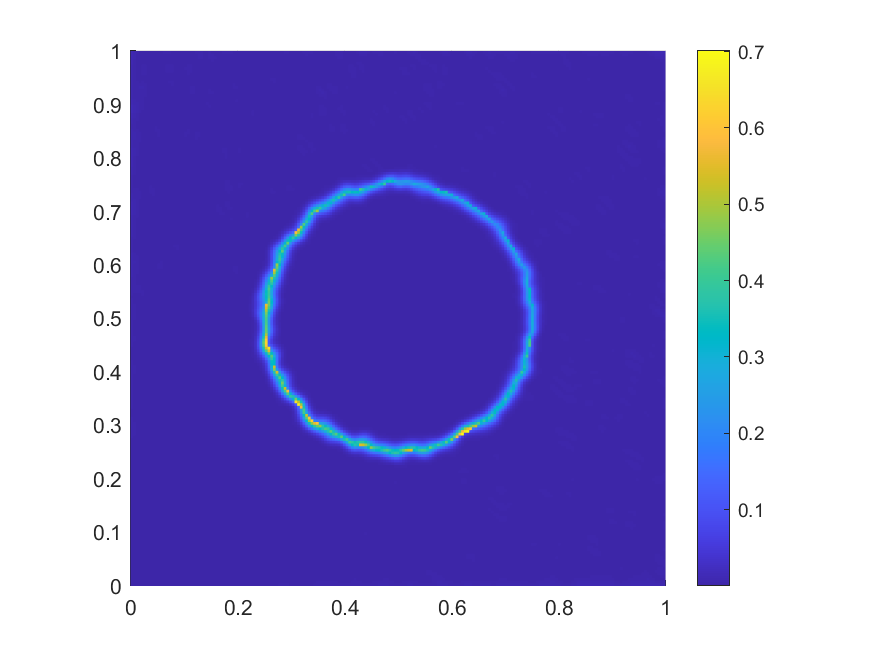} 	
			\\\multicolumn{3}{c}{DD-MLS$^0_{\text{G}}$}\\
\hspace{-1cm}			\includegraphics[width=6.3cm]{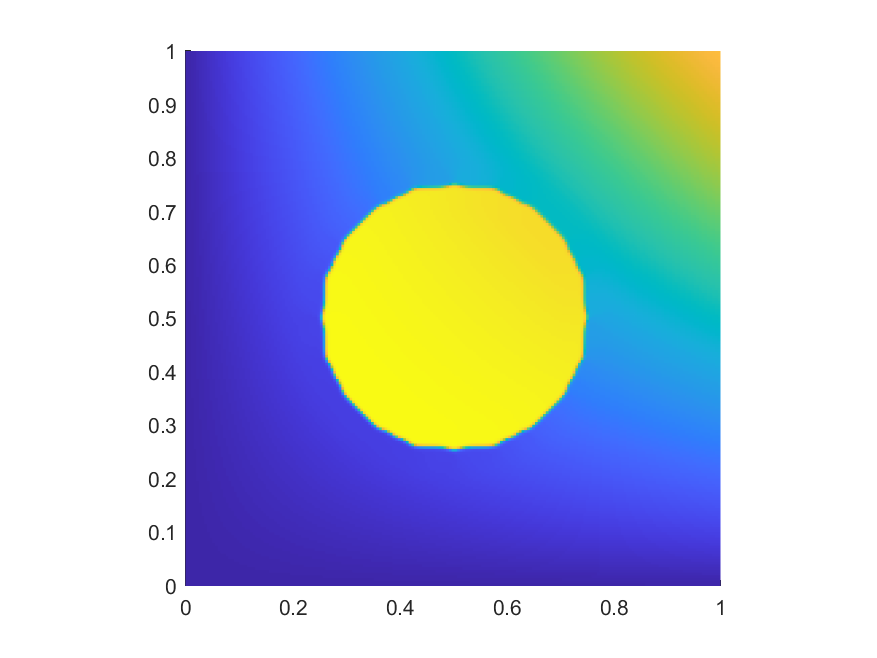}  &
\hspace{-1cm}			\includegraphics[width=6.3cm]{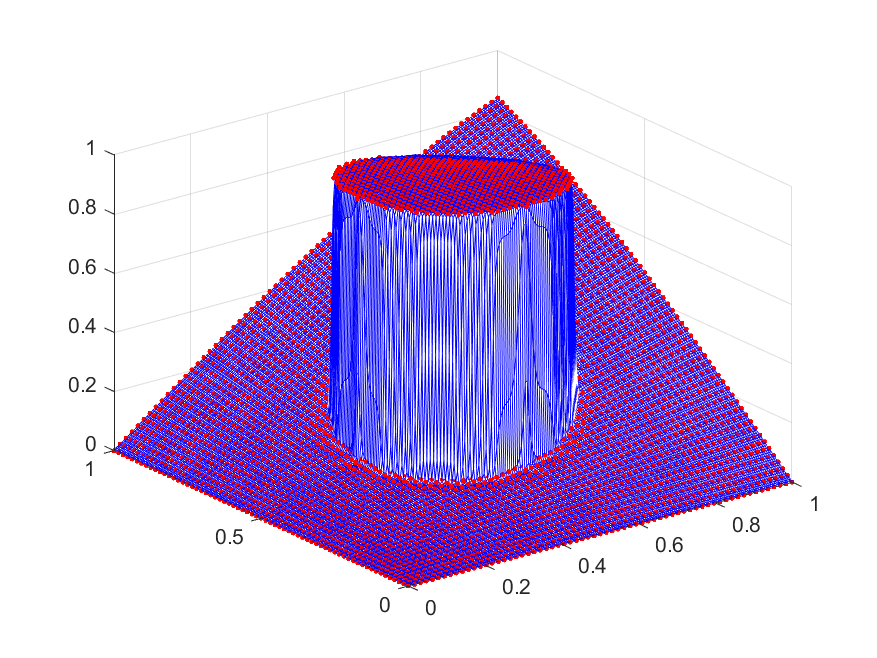} &
\hspace{-1cm}            \includegraphics[width=6.3cm]{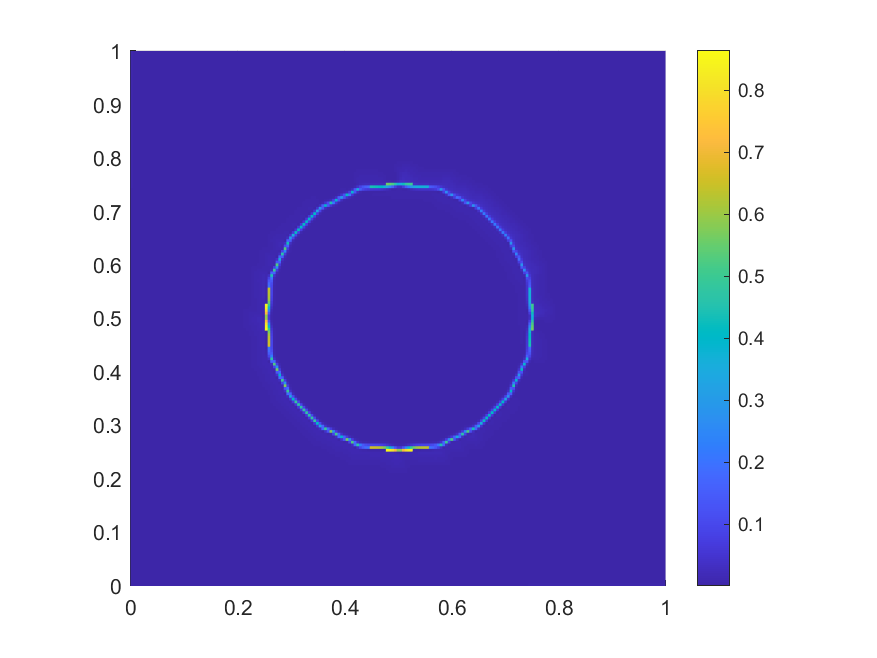}  	
				\\
\hspace{-1cm}			\includegraphics[width=6.3cm]{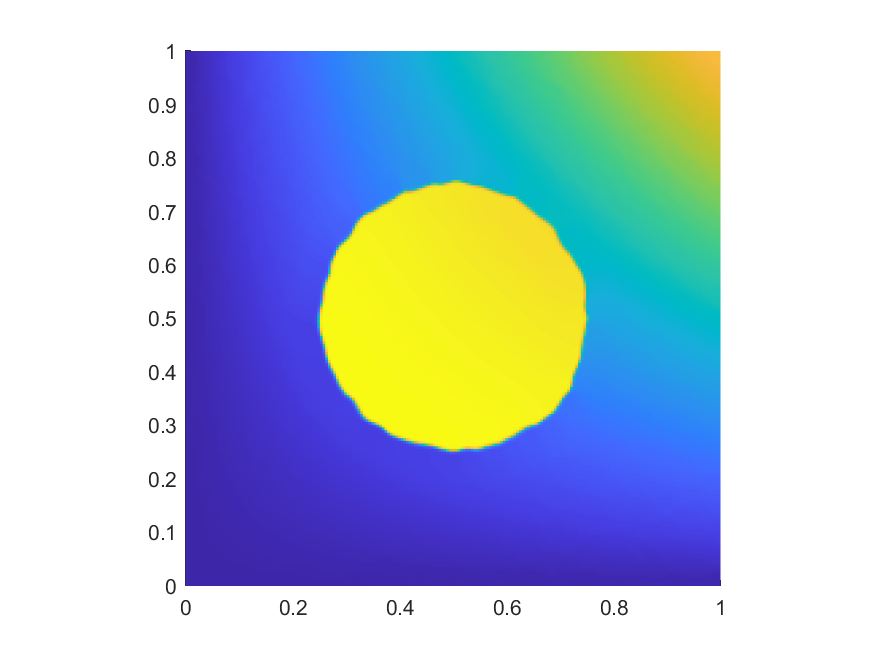}  &
\hspace{-1cm}			\includegraphics[width=6.3cm]{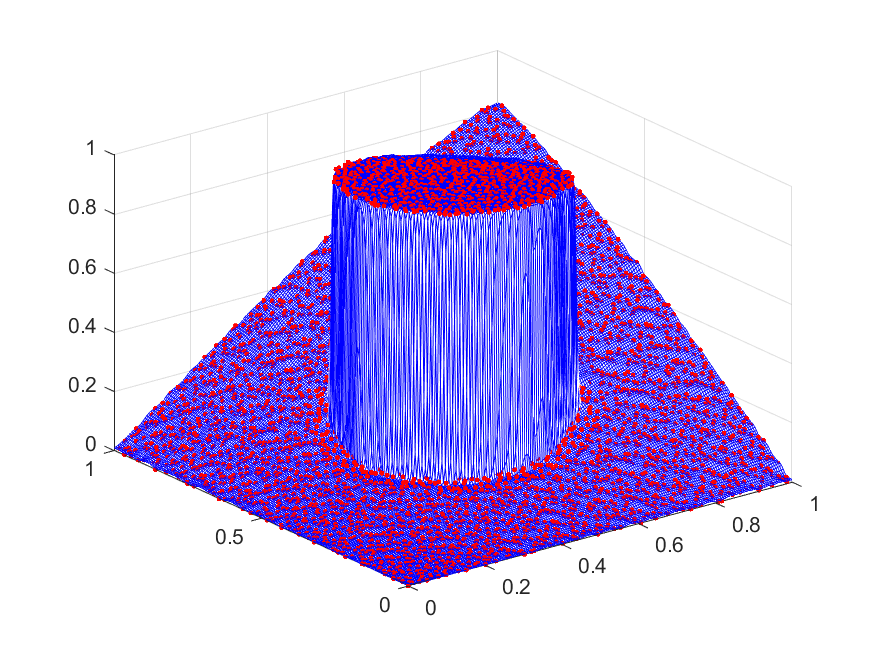} &
\hspace{-1cm}            \includegraphics[width=6.3cm]{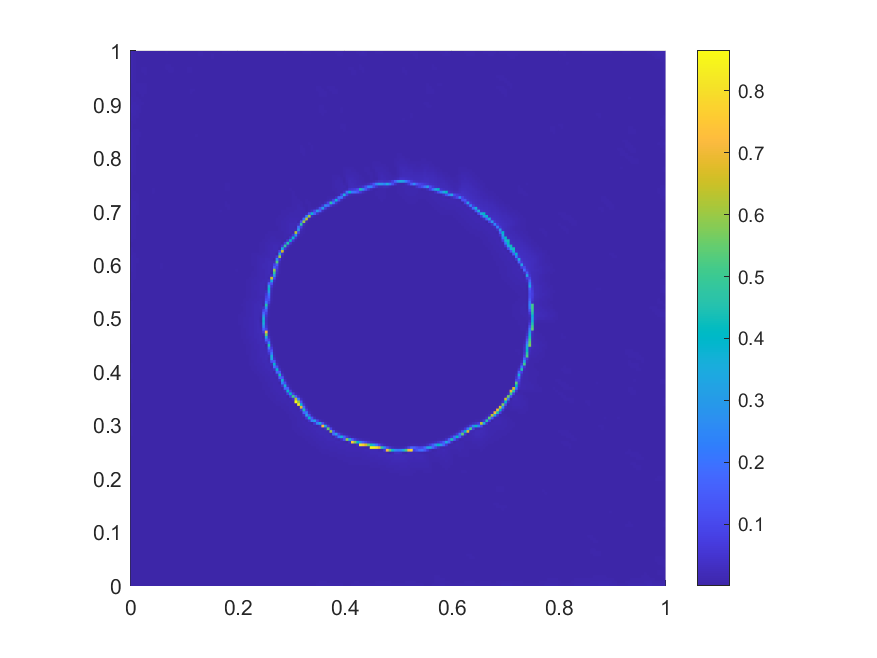} 	

		\end{tabular}
\end{center}
				\caption{Approximation to function $z$, Eq. \eqref{ejemplonuma}, using G, $d=0$ with $N=65^2$.}
		\label{figuraGd0N65}
	\end{figure}

In this section we study the results obtained when we employ the set of polynomials of degree $d=0$, general Shepard's method (see \cite{FASSHAUER}), or $d=1$. In these cases, the non-desired oscillations do not occur, but some diffusion effects appear when the linear
MLS algorithms are used. We compare these linear algorithms with the data-dependent ones using the same function, $z$, Eq. \eqref{ejemplonuma}. We perform experiments with $N=33^2$, Figures \ref{figuraW2d0N33} and \ref{figuraGd0N33} and $N=65^2$, Figures \ref{figuraW2d0N65} and \ref{figuraGd0N65}, using W2 and G functions, and for gridded and Halton's data points. In all the results, we can observe that the discontinuity is perfectly delineated in both the case of regular meshes and with Halton's points. Particularly striking is the result obtained for the function G, Figures \ref{figuraGd0N33} and \ref{figuraGd0N65}, where the band of diffusion around the discontinuity is reduced considerably.

\begin{figure}[!ht]
\begin{center}
		\begin{tabular}{cccc}
\multicolumn{2}{c}{MLS$^1_{\text{W2}}$}& \multicolumn{2}{c}{DD-MLS$^1_{\text{W2}}$}\\
\hspace{-1cm}		\includegraphics[width=5.7cm]{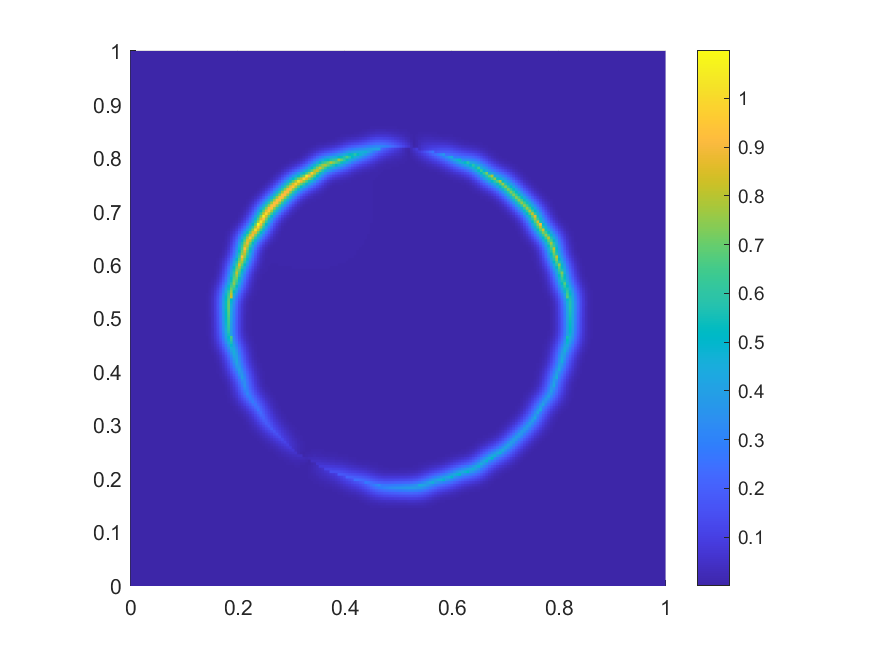} & \hspace{-1.25cm}		\includegraphics[width=5.7cm]{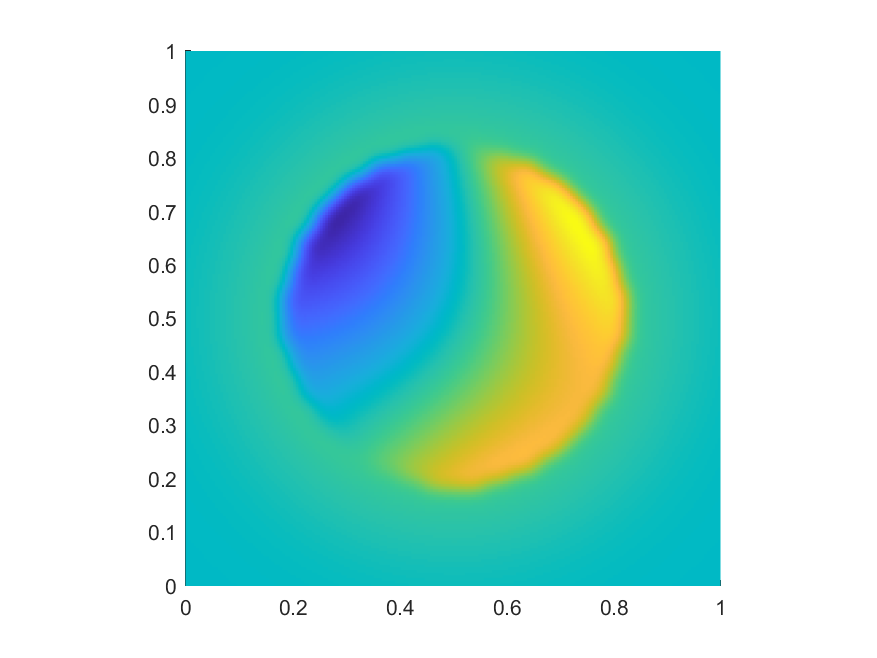} & \hspace{-1.5cm}	
			\includegraphics[width=5.7cm]{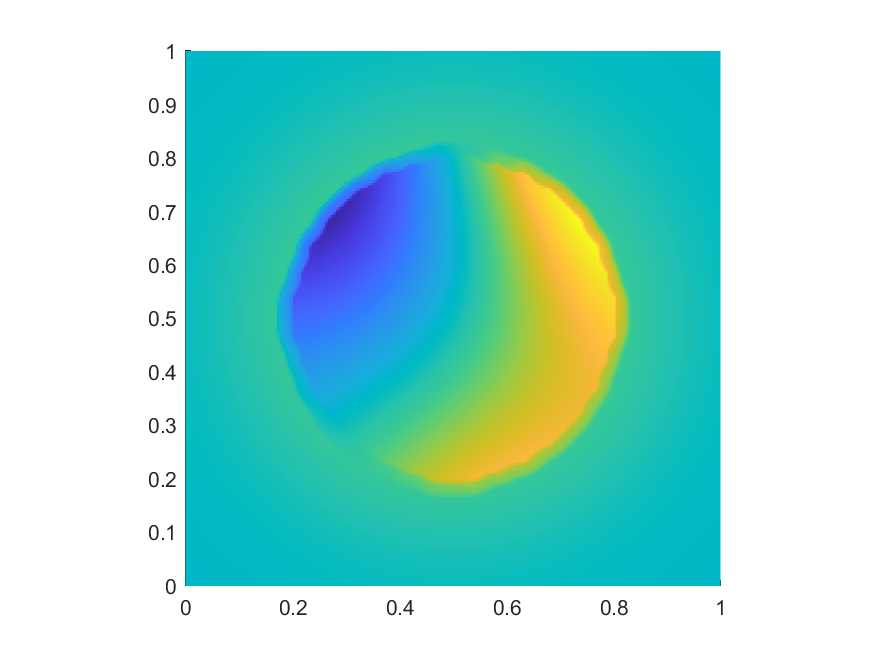}  & \hspace{-1.5cm}	
			\includegraphics[width=5.7cm]{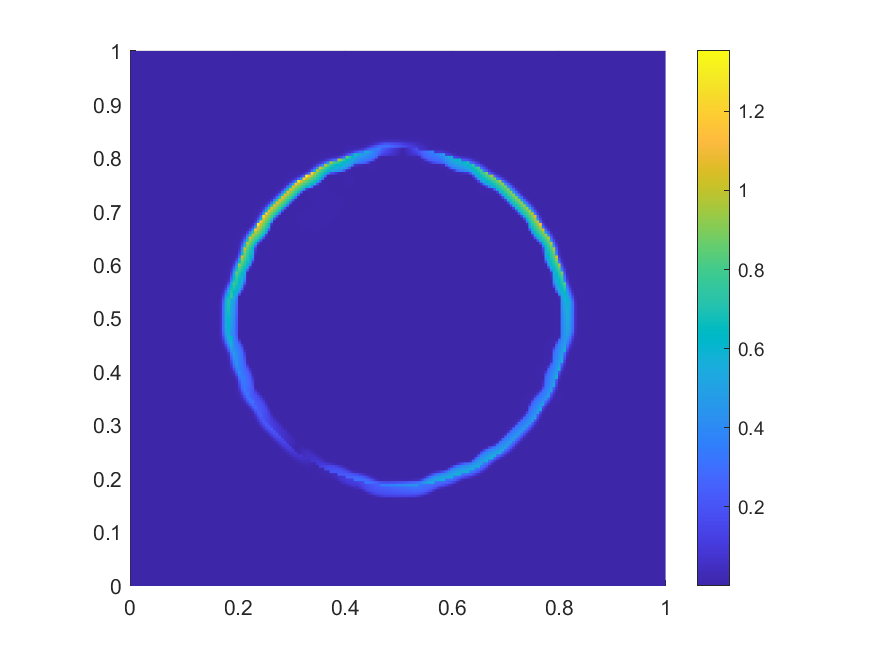} \\
\hspace{-1cm}		\includegraphics[width=5.7cm]{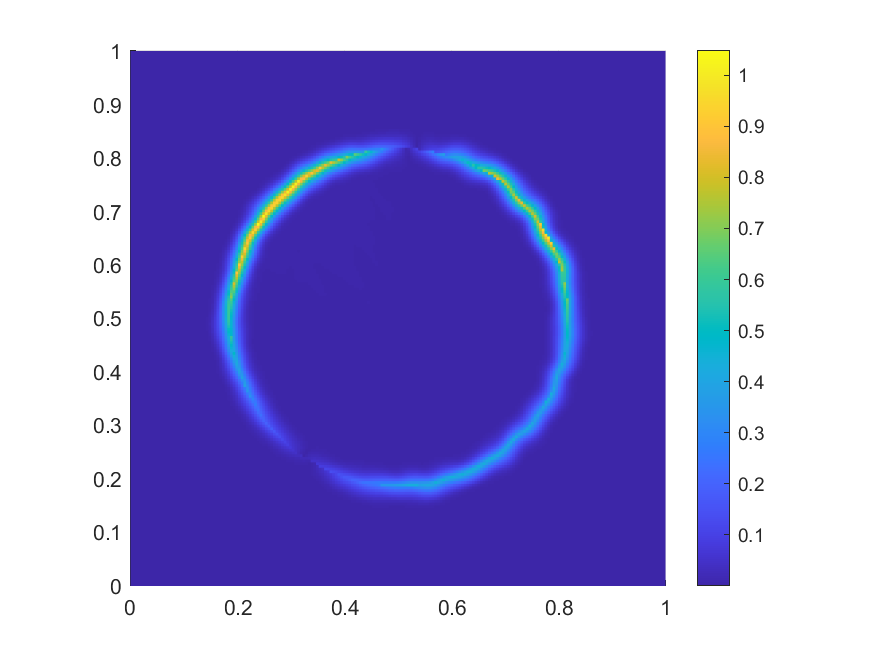} & \hspace{-1.25cm}		\includegraphics[width=5.7cm]{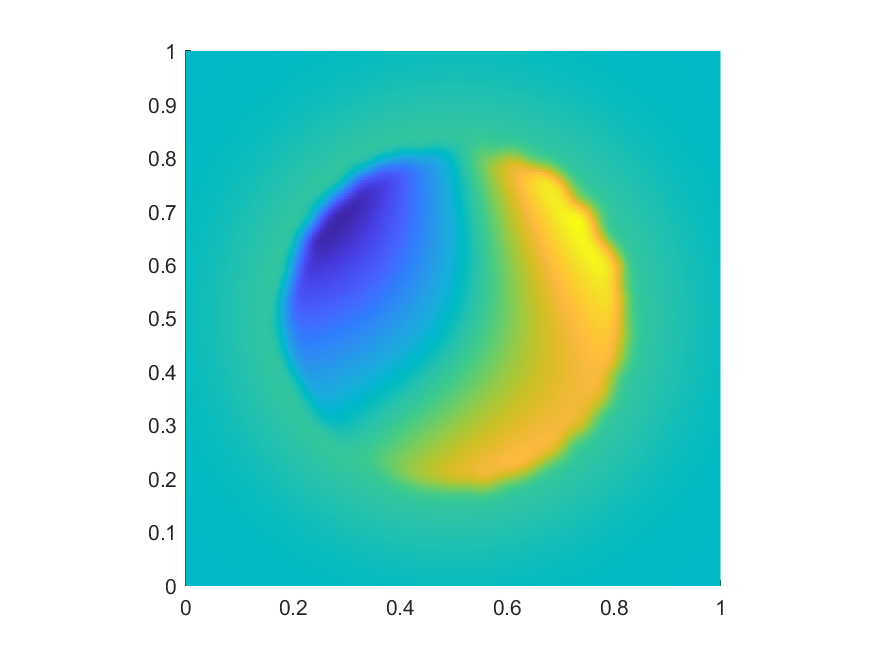} & \hspace{-1.5cm}	
			\includegraphics[width=5.7cm]{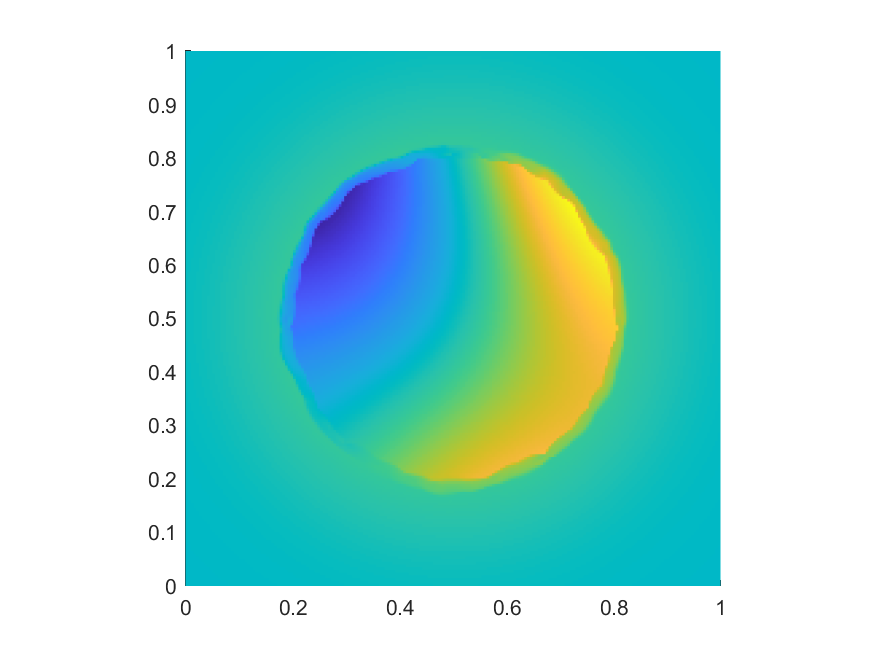}  & \hspace{-1.5cm}	
			\includegraphics[width=5.7cm]{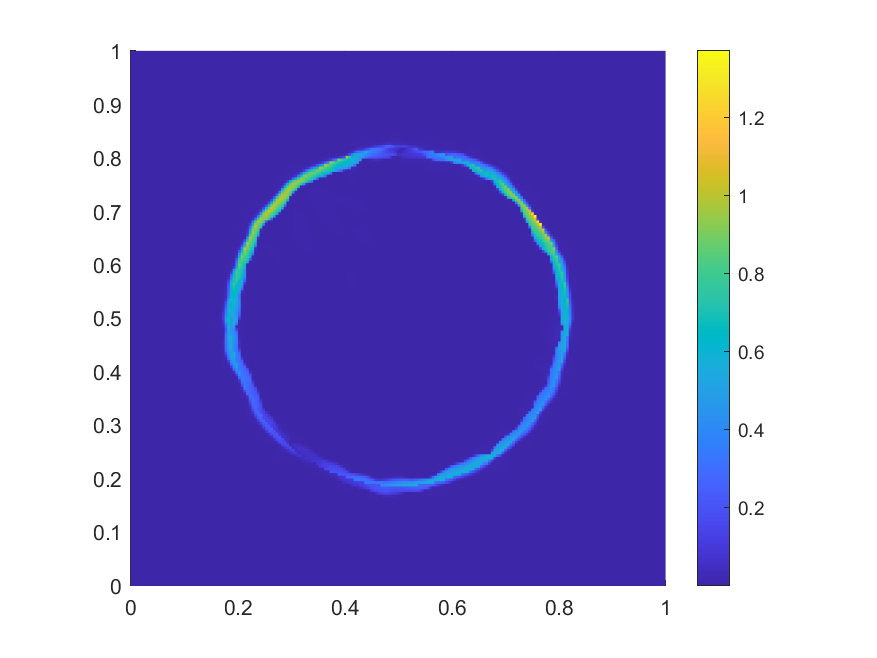} \\
\hspace{-1cm}		\includegraphics[width=5.7cm]{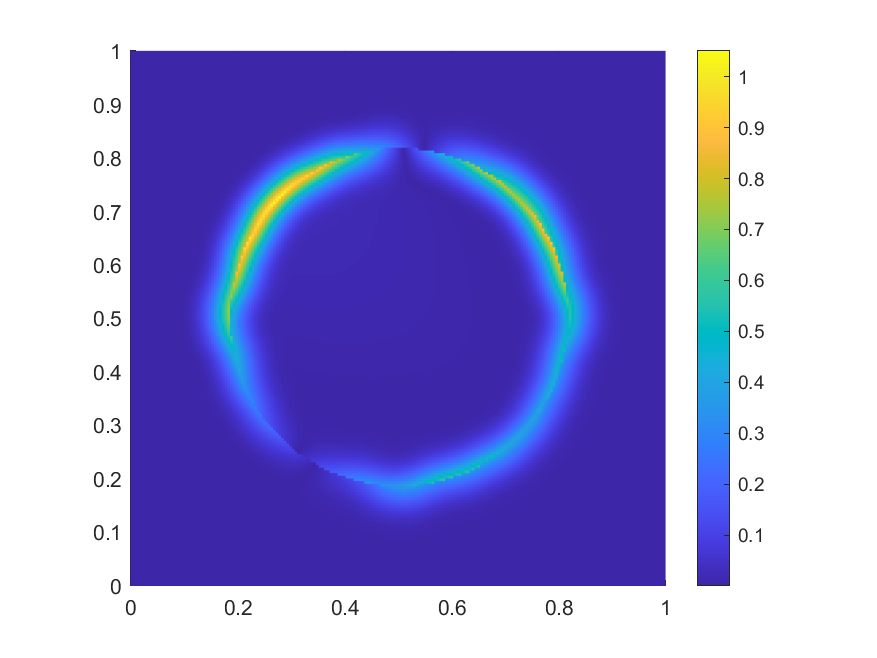} & \hspace{-1.25cm}		\includegraphics[width=5.7cm]{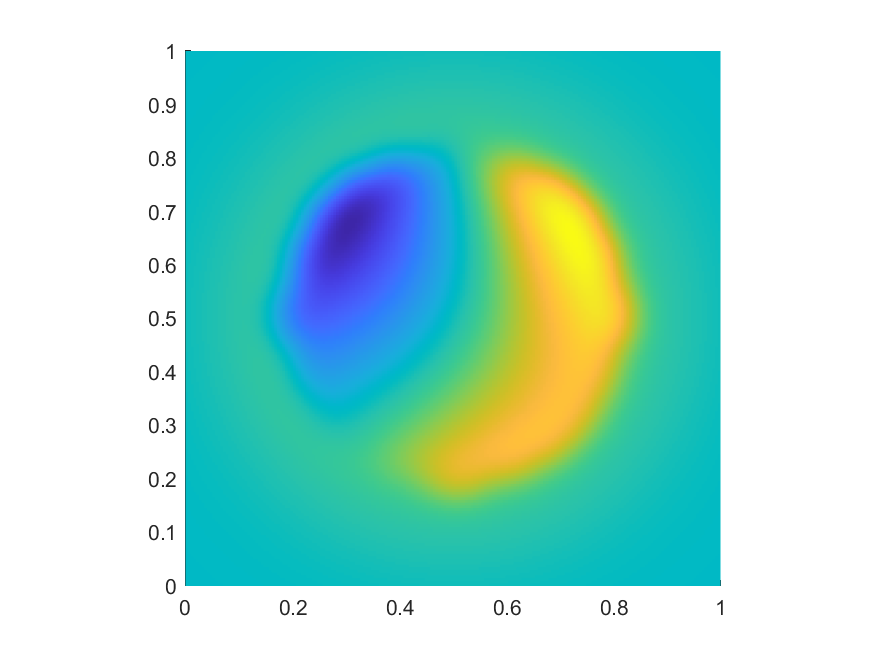} & \hspace{-1.5cm}	
			\includegraphics[width=5.7cm]{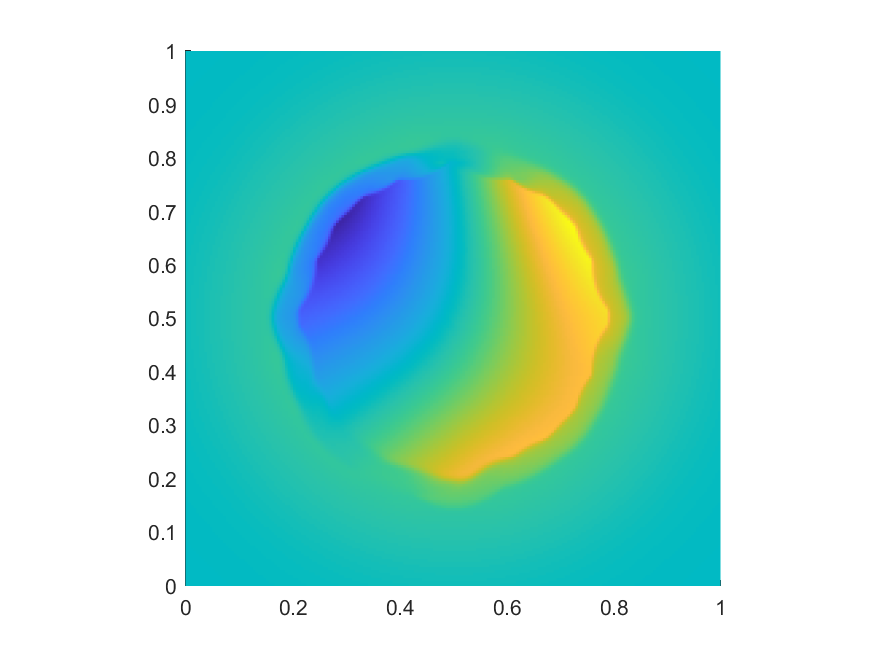}  & \hspace{-1.5cm}	
			\includegraphics[width=5.7cm]{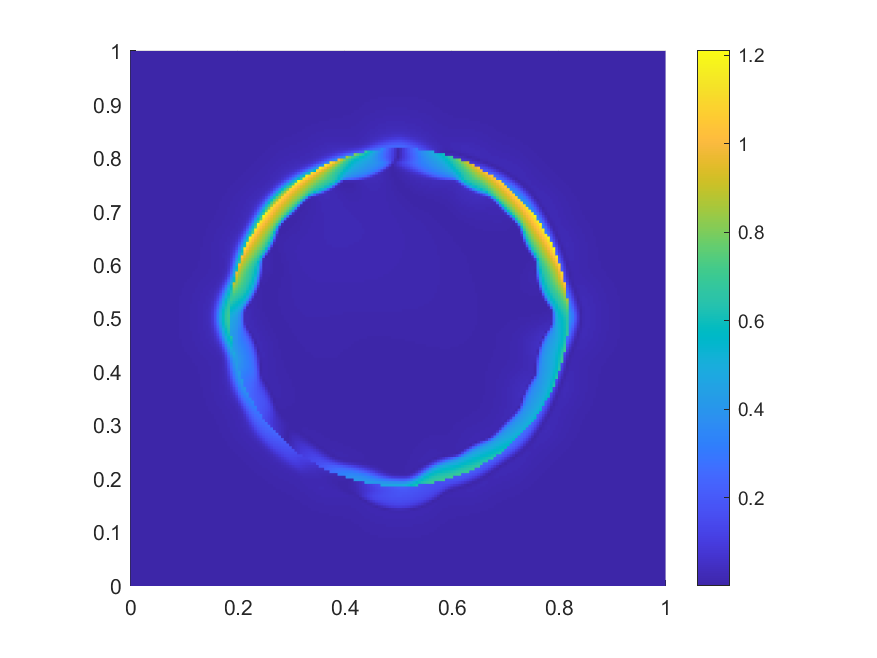} \\
\hspace{-1cm}		\includegraphics[width=5.7cm]{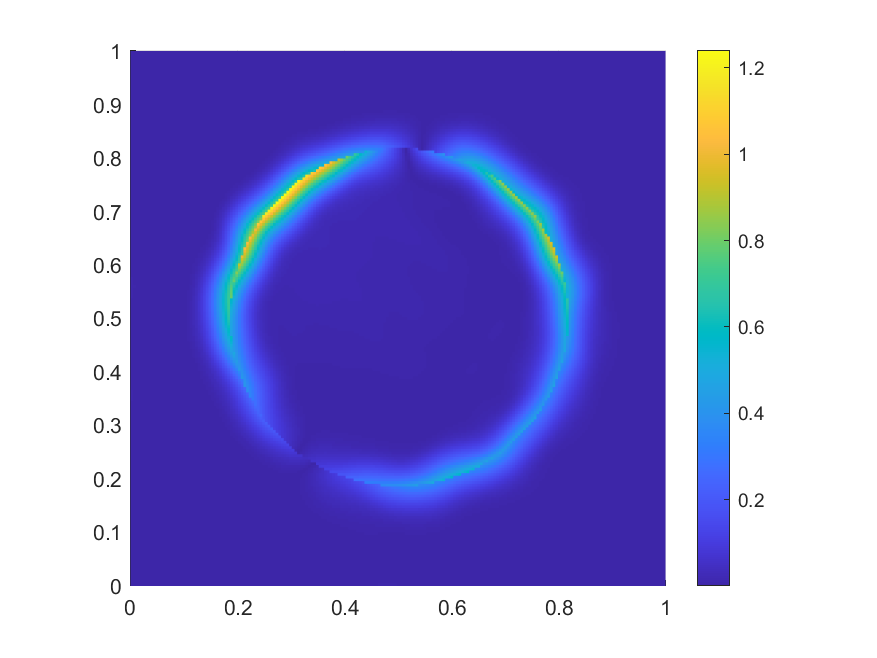} & \hspace{-1.25cm}		\includegraphics[width=5.7cm]{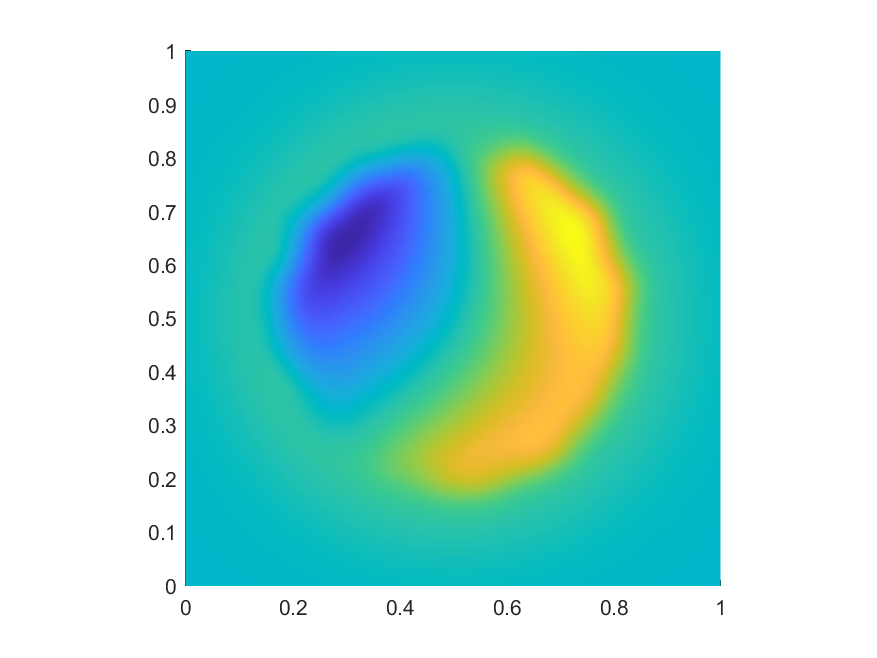} & \hspace{-1.5cm}	
			\includegraphics[width=5.7cm]{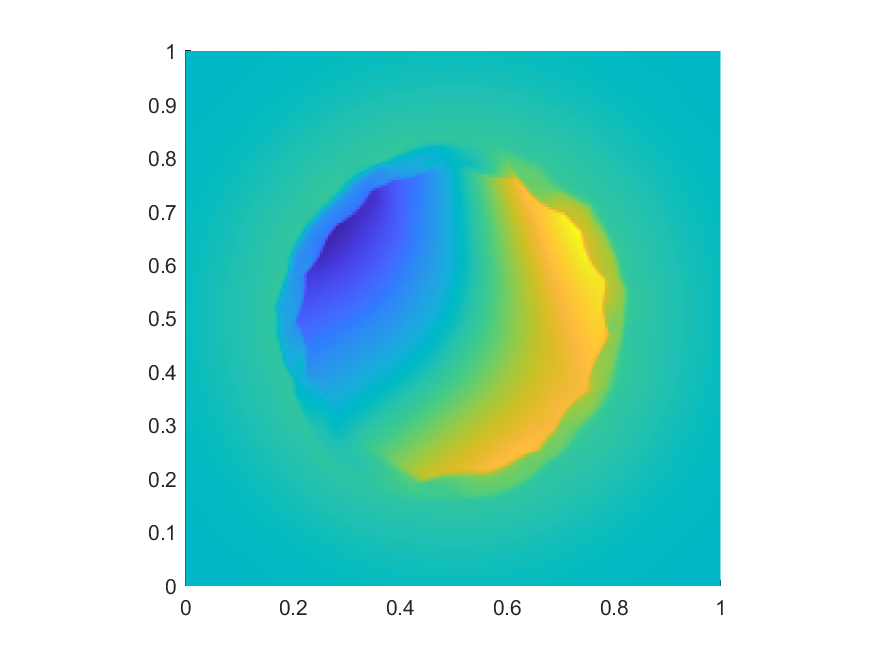}  & \hspace{-1.5cm}	
			\includegraphics[width=5.7cm]{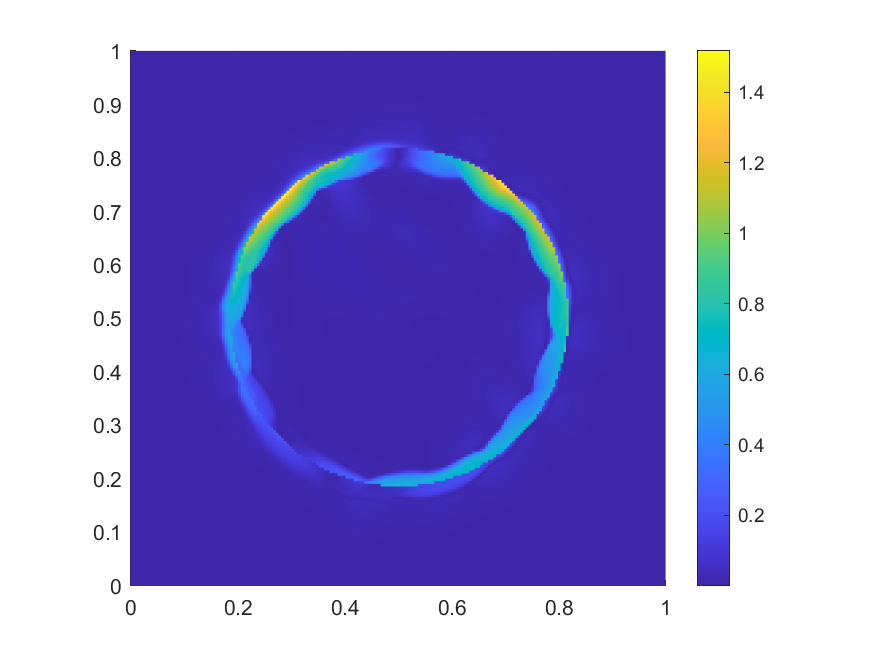}
		\end{tabular}
\end{center}
				\caption{Approximation to function $f$, Eq. \eqref{ejemplolevin}, and errors in a 2-D plot using regular grid data points with $N=65^2, 33^2$.}
		\label{figurad11}
	\end{figure}

\begin{figure}[!ht]
\begin{center}
		\begin{tabular}{cccc}
\multicolumn{2}{c}{MLS$^1_{\text{W2}}$}& \multicolumn{2}{c}{DD-MLS$^1_{\text{W2}}$}\\
\hspace{-1cm}		\includegraphics[width=5.7cm]{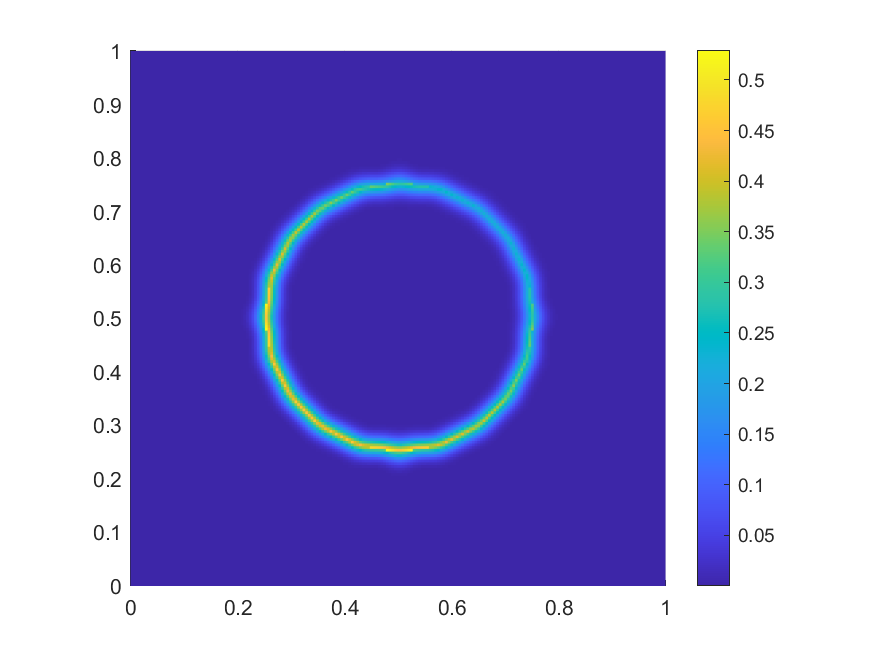} & \hspace{-1.25cm}		\includegraphics[width=5.7cm]{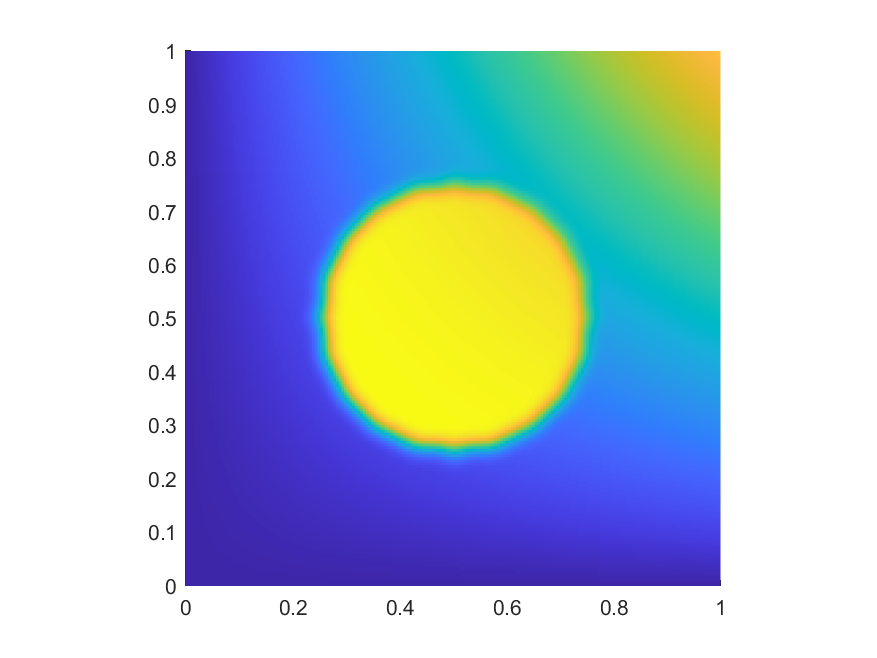} & \hspace{-1.5cm}	
			\includegraphics[width=5.7cm]{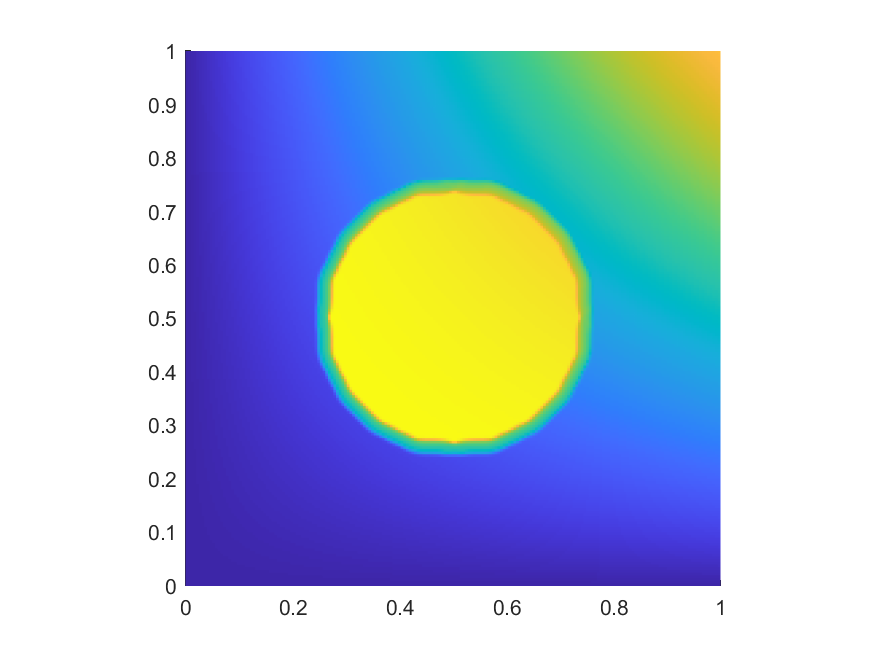}  & \hspace{-1.5cm}	
			\includegraphics[width=5.7cm]{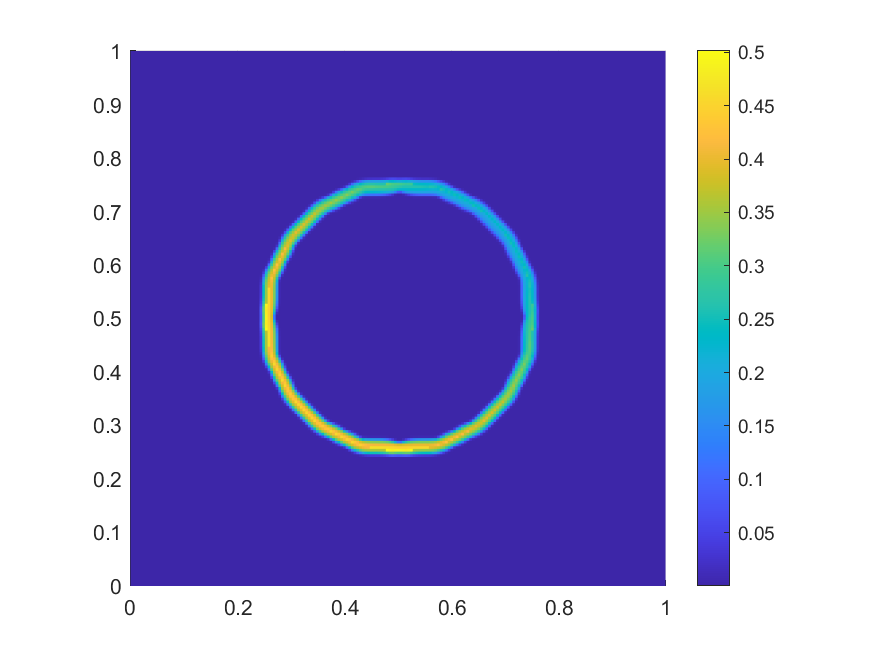} \\
\hspace{-1cm}		\includegraphics[width=5.7cm]{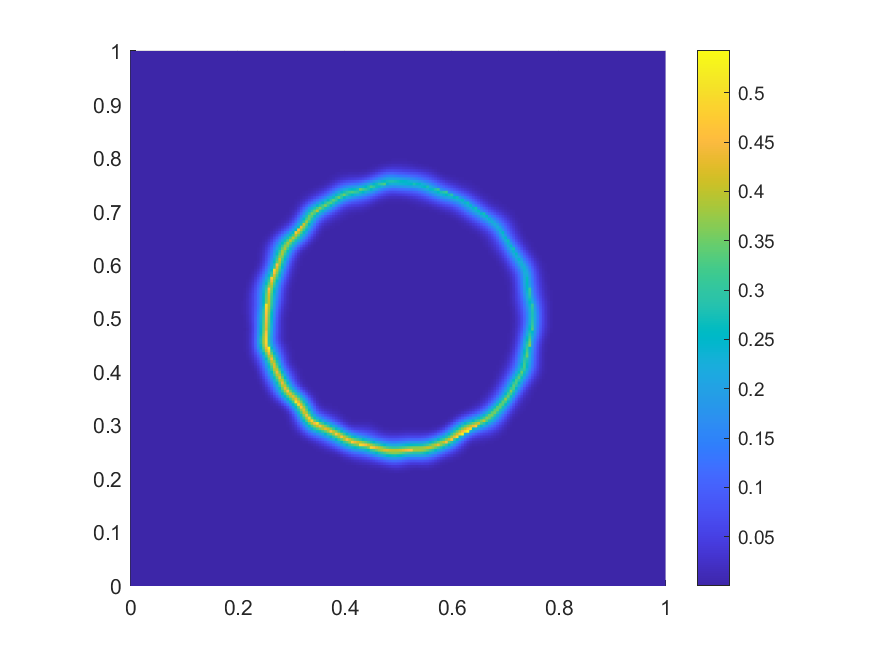} & \hspace{-1.25cm}		\includegraphics[width=5.7cm]{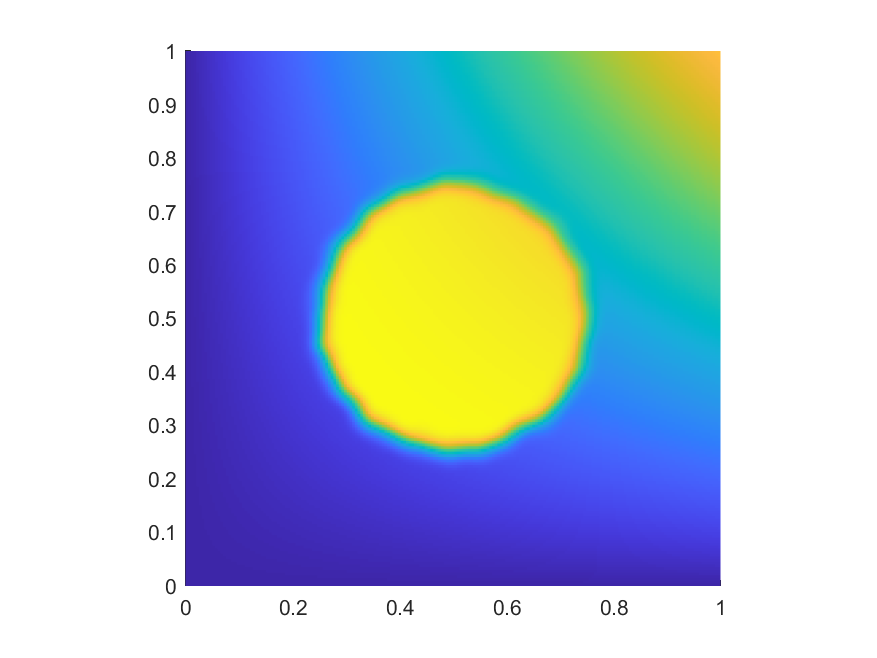} & \hspace{-1.5cm}	
			\includegraphics[width=5.7cm]{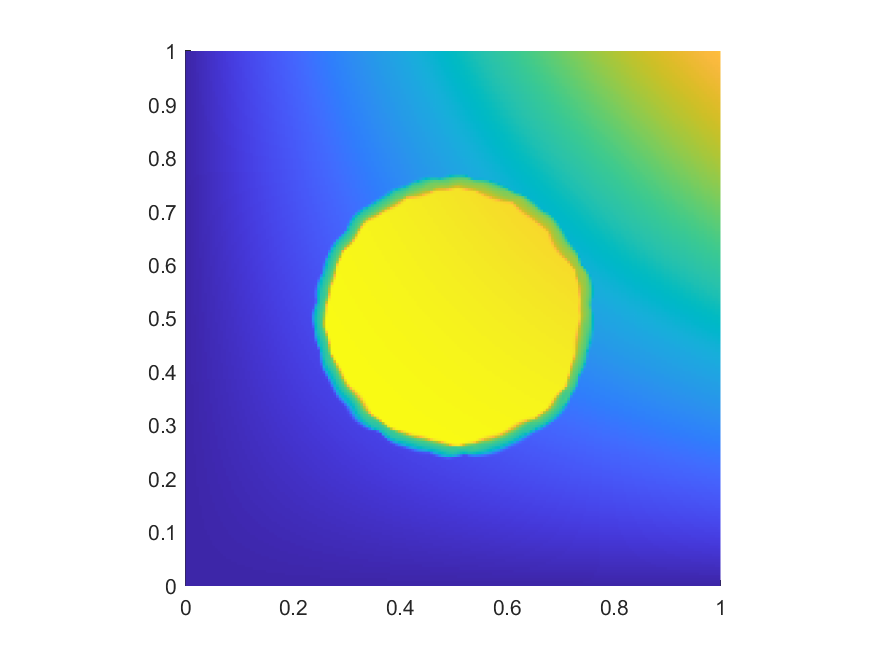}  & \hspace{-1.5cm}	
			\includegraphics[width=5.7cm]{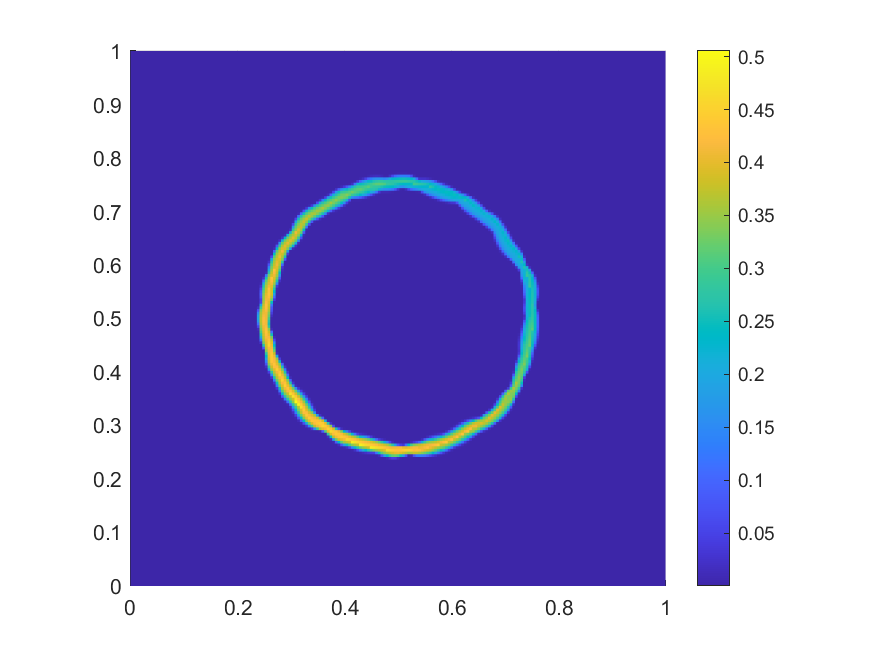} \\
\hspace{-1cm}		\includegraphics[width=5.7cm]{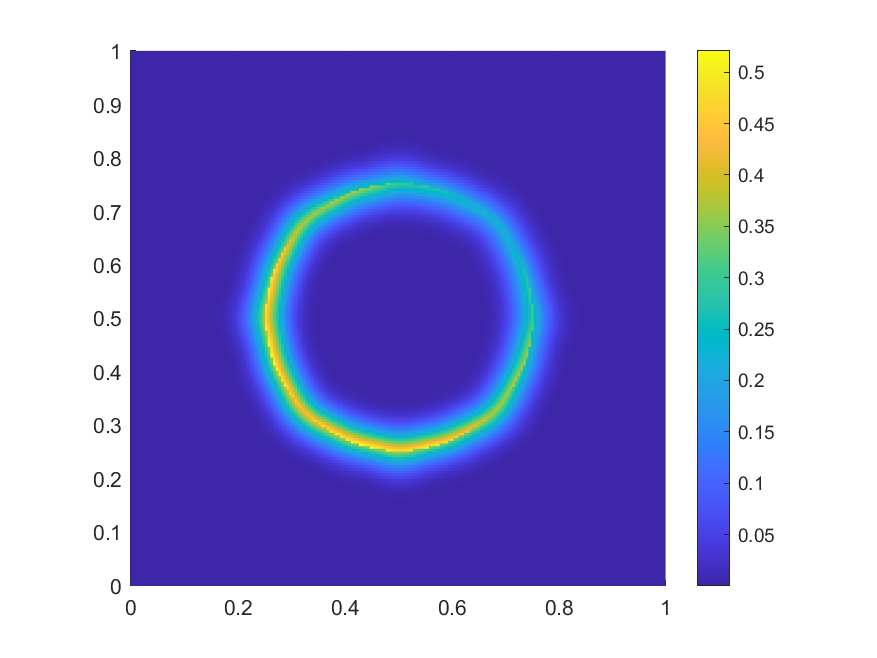} & \hspace{-1.25cm}		\includegraphics[width=5.7cm]{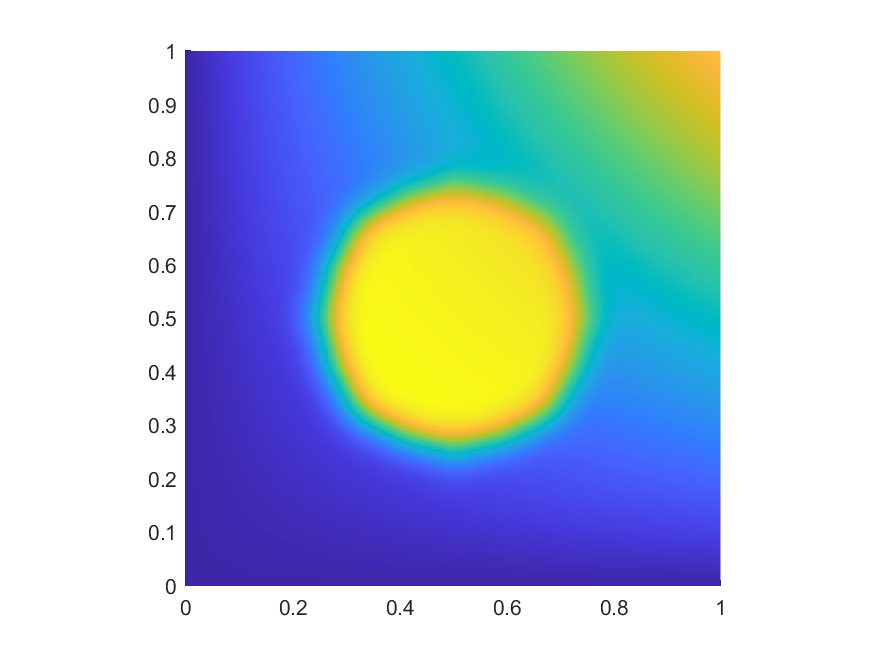} & \hspace{-1.5cm}	
			\includegraphics[width=5.7cm]{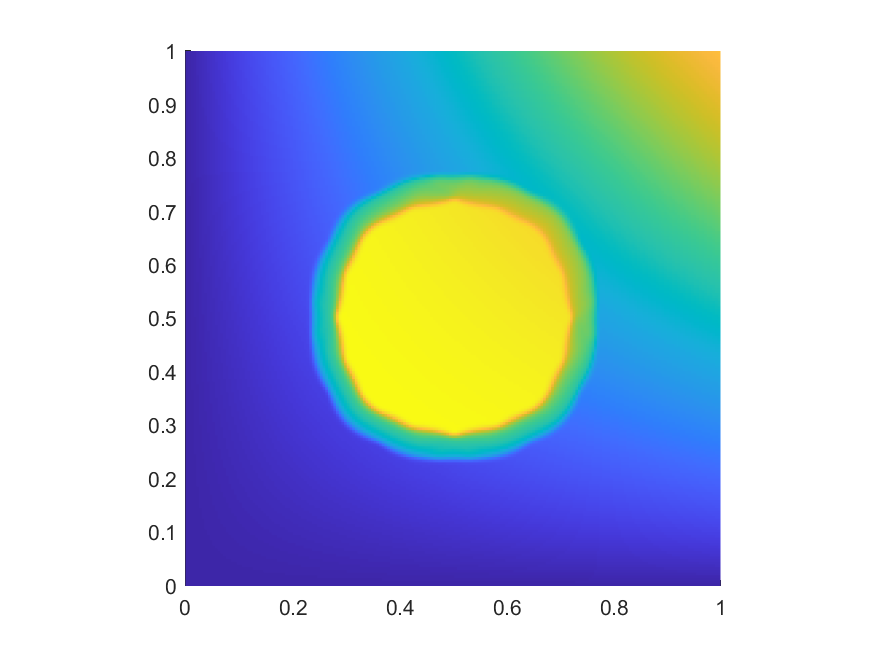}  & \hspace{-1.5cm}	
			\includegraphics[width=5.7cm]{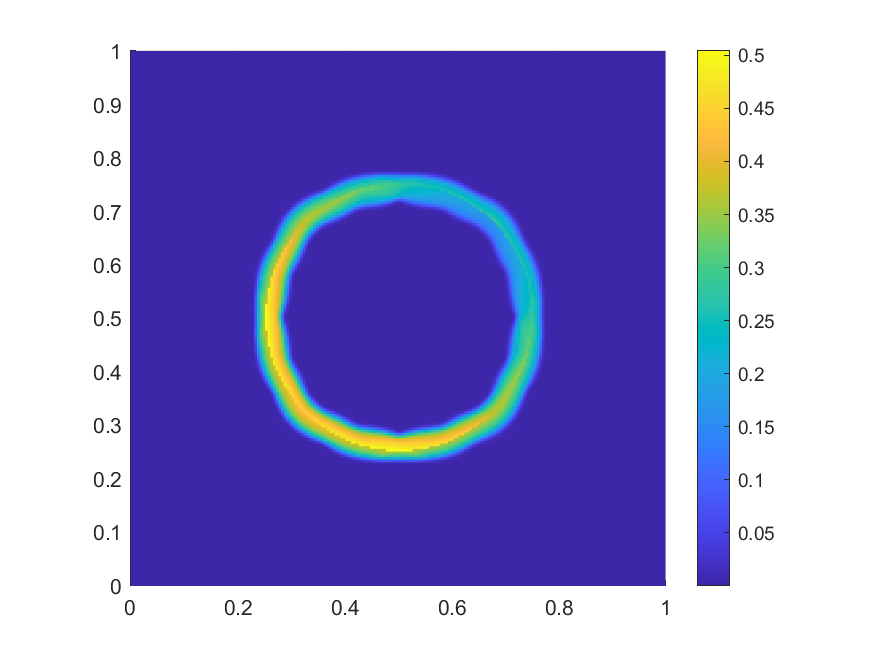} \\
  \\
\hspace{-1cm}		\includegraphics[width=5.7cm]{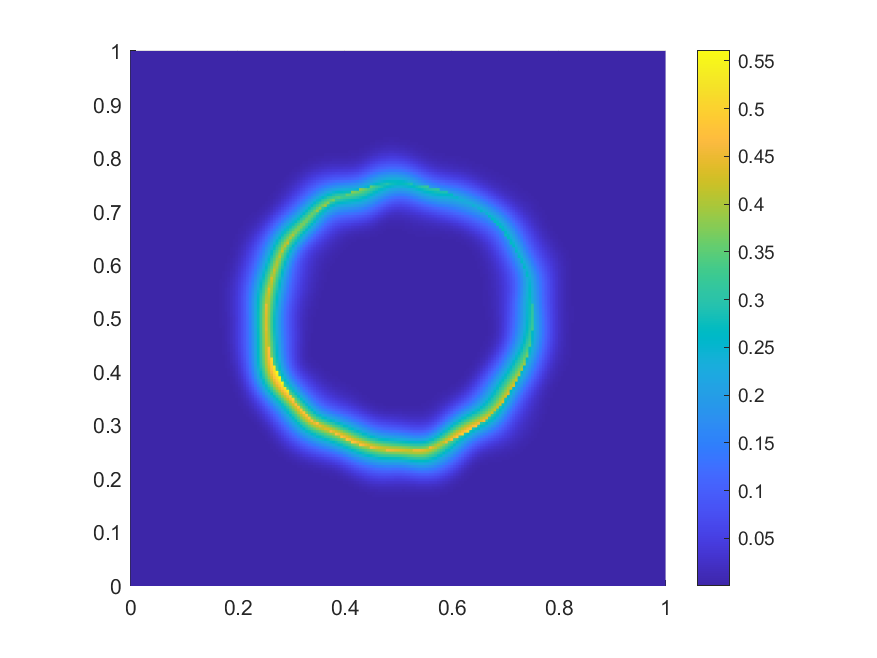} & \hspace{-1.25cm}		\includegraphics[width=5.7cm]{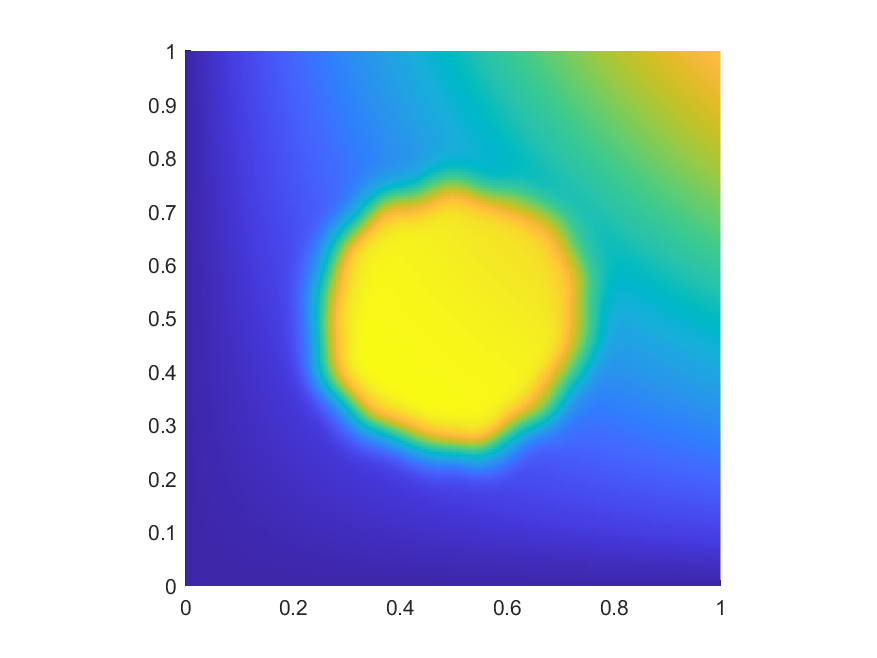} & \hspace{-1.5cm}	
			\includegraphics[width=5.7cm]{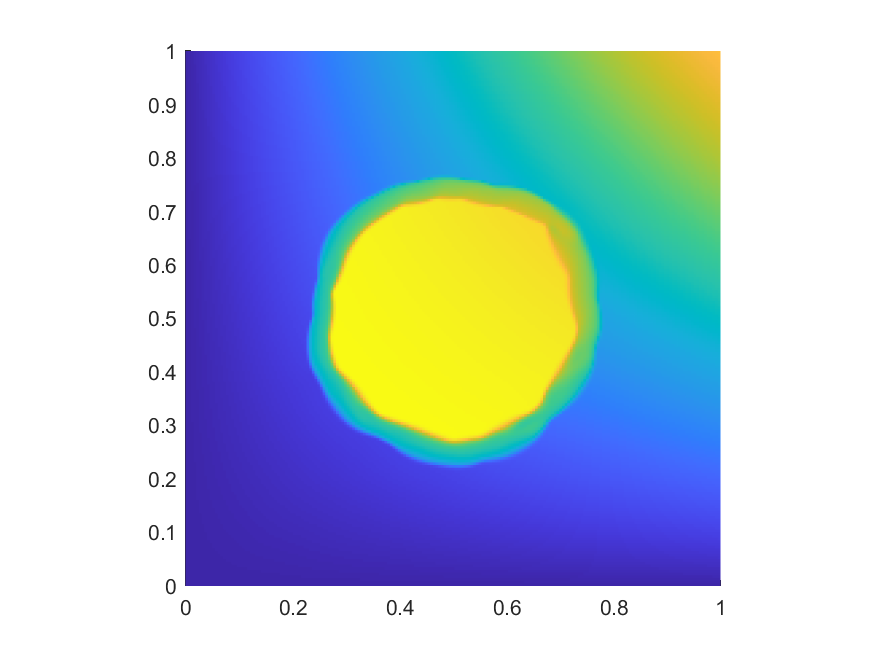}  & \hspace{-1.5cm}	
			\includegraphics[width=5.7cm]{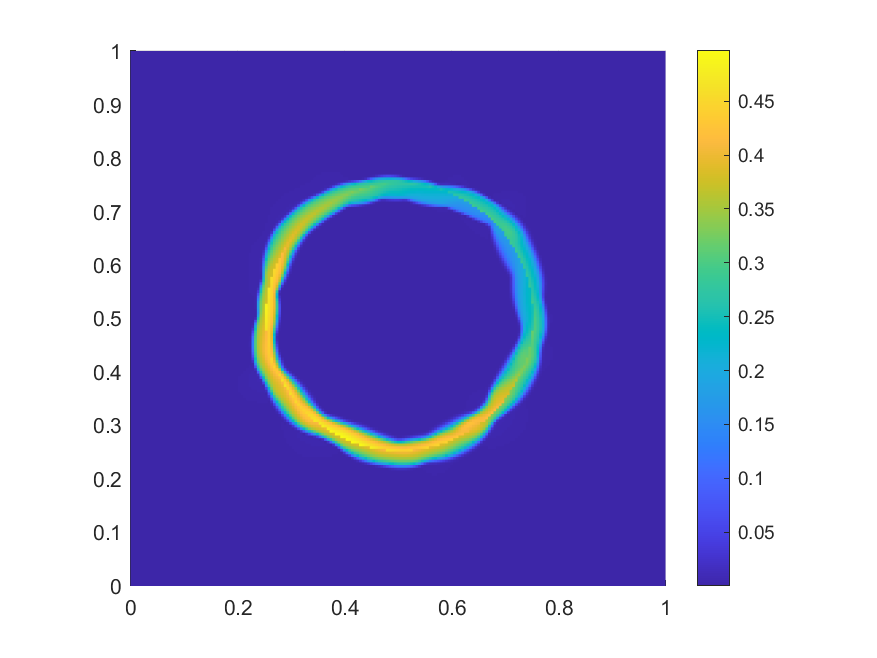}
		\end{tabular}
\end{center}
				\caption{Approximation to function $z$, Eq. \eqref{ejemplonuma}, and errors in a 2-D plot using regular grid data points  with $N=65^2, 33^2$.}
		\label{figurad12}
	\end{figure}

Finally, some examples with $d=1$ are shown in Figures \ref{figurad11} and \ref{figurad12}. Again, the diffusion effects decrease using DD-MLS.

\section{Conclusions and Future work}

In this study, we have introduced a novel approach to the moving least squares (MLS) problem (\ref{MLSproblem}) by replacing the traditional weight functions with new functions that assign greater weight to nodes farther from discontinuities, while still assigning smaller weights to nodes far from the point of approximation. This adjustment effectively mitigates the Gibbs phenomenon and reduces the smearing of discontinuities in the final approximation of the original data.

Our method uses smoothness indicators to accurately identify {\it infected} nodes, i.e. those affected by the presence of discontinuities, in a way inspired by the WENO method. This results in a data-dependent weighted least squares problem where the weights are influenced by both the distances between nodes and the point of approximation, and the distances between isolated discontinuities and the nodes. We think that these criteria could be adapted to other requirements, such as point density or monotonicity, but we leave these ideas for future explorations.

Through the design and analysis of the new data dependent approximation, we have demonstrated its properties, including {polynomial reproduction}, accuracy, smoothness or neglecting of the Gibbs oscillations close to the discontinuities. Our numerical experiments validate the theoretical findings, showing the effectiveness of the proposed method.

{

}

\end{document}